\documentclass[12pt]{amsart}
\usepackage[left=2.90cm, right=2.90cm, top=3.30cm, bottom=3.30cm]{geometry}

\usepackage{amsfonts,graphics,amsmath,amsthm,amsfonts,amscd,amssymb,amsmath,latexsym,multicol,
 mathrsfs}
\usepackage{epsfig,url}
\usepackage{flafter}
\usepackage{fancyhdr}
\usepackage[hyperindex,breaklinks]{hyperref}
\hypersetup{colorlinks=true, linkcolor=black, breaklinks=true}

\usepackage{graphicx}
\usepackage{xcolor}
\usepackage{float}
\usepackage{bm}
\usepackage{enumitem}
\usepackage{multirow}
\usepackage{scalerel,stackengine}
\usepackage{stmaryrd}
\usepackage{datetime}
\usepackage[title]{appendix}
\usepackage{url}
\usepackage{tabulary}
\usepackage{booktabs}
\usepackage{tikz}
\usetikzlibrary{positioning}
\usepackage{verbatim}

\usepackage[utf8]{inputenc}
\usepackage[english]{babel}
\usepackage{csquotes}

\usepackage{comment}

\usepackage{CJKutf8}

\usepackage{quiver}

\theoremstyle{cupplain}
\newtheorem{theorem}{Theorem}[section]
\newtheorem{lemma}[theorem]{Lemma}
\newtheorem{corollary}[theorem]{Corollary}

\newtheorem{proposition}[theorem]{Proposition}
\newtheorem{conjecture}[theorem]{Conjecture}

\newtheorem{assumption}[theorem]{Assumptions}

\theoremstyle{cupremark}
\newtheorem{remark}[theorem]{Remark}

\theoremstyle{cupdefinition}
\newtheorem{definition}[theorem]{Definition} 

\theoremstyle{definition}
\newtheorem*{Ack}{Acknowledgment}

\title[On the cone conjecture for pairs of large Kodaira dimension]{On the cone conjecture for pairs of large Kodaira dimension}
\author{Fulin Xu
}
\date{
}

\begin{document}
\begin{abstract}
In this paper, by running the minimal model program and considering the anti-canonical fibration, we prove the Morrison-Kawamata cone conjecture for klt Calabi-Yau pairs $(X,\Delta)$ such that the Iitaka dimension $\kappa(X,-K_X)$ is at least $\dim X - 2$. 
\end{abstract}
\maketitle


\tableofcontents

\section{Introduction}

One of the central ideas of birational geometry is to study all 
possible contractions from a fixed projective variety $X$ by looking at the cone of curves in $N_1(X)$, or the dual cone, the nef cone in $N^1(X)$. This leads to the study of the behavior of the cones. 

Also, minimal model program predicts that there are three types of building blocks for algebraic varieties, i.e. Fano varieties, Calabi-Yau varieties, and varieties of general type. For (log) Fano varieties, the cone is rational polyhedral by cone theorem. For
varieties of general type, the cone can be far from rational polyhedral, while the birational automorphism group is always finite.

The Morrison-Kawamata cone conjecture  would give a complete picture in this story. 
The conjecture predicts that the effective nef cone of a Calabi-Yau variety admits a rational polyhedral
fundamental domain under the action of the automorphism group. This means that Calabi-Yau varieties lie precisely in the middle of general type varieties and Fano varieties: the
cones can be wild, but the failure of becoming rational polyhedral always should come from an
infinite automorphism group. It is actually not clear where these automorphisms should
come from, and all known cases follow from some specific kind of geometry.

The first formulation of the cone conjecture is due to Morrison and Kawamata, see \cite{Mor96} and \cite{Kaw97}. 
We consider the following version of the cone conjecture, formulated by Totaro in \cite{Tot09}:

\begin{conjecture}
Let $(X/S,\Delta)$ be a projective klt Calabi-Yau pair, where $X$ is $\mathbb{Q}$-factorial with $\Delta$ an $\mathbb{R}$-boundary. Then 

(a) The number of ${\rm Aut}(X/S, \Delta)$-equivalence classes of faces of the cone $\mathcal{A}^e(X/S)$
corresponding to birational contractions or fiber space structures is finite. Moreover, the action of ${\rm Aut}(X/S,\Delta)$ on the effective nef cone $\mathcal{A}^e(X/S)$ admits a rational polyhedral fundamental domain $\Pi$, in the sense that: 

(1) $\mathcal{A}^e(X/S) = \bigcup_{g\in {\rm Aut}(X/S,\Delta)}g_*\Pi$, 

(2) $\text{Int }\Pi \cap g_*\text{Int }\Pi= \emptyset$ unless $g_* = 1$.

(b) The number of ${\rm PsAut}(X/S, \Delta)$-equivalence classes of chambers $\mathcal{A}^e
(X_0 /S)$
in the cone $\mathcal{M}^e
(X/S)$ corresponding to marked small $\mathbb{Q}$-factorial modifications $X_0 \to S$ of $X \to S$ is finite. Moreover, the action of ${\rm PsAut}(X/S,\Delta)$ on the effective movable cone $\mathcal{M}^e(X/S)$ admits a rational polyhedral fundamental domain $\Pi^\prime$, in the sense that: 

(1) $\mathcal{M}^e(X/S) = \bigcup_{g\in {\rm PsAut}(X/S,\Delta)}g_*\Pi^\prime$, 

(2) $\text{Int }\Pi^\prime \cap g_*\text{Int }\Pi^\prime= \emptyset$ unless $g_* = 1$. 

\end{conjecture}

\begin{remark}\label{remark cone}
The first (resp. second) part of (a) is usually refered as the weak nef cone conjecture (resp. nef cone conjecture), and the first (resp. second) part of (b) is usually refered as the weak movable cone conjecture (resp. movable cone conjecture). Recall that a small $\mathbb{Q}$-factorial modification (SQM) of $X$ over $S$ means a pseudo-isomorphism over $S$ from $X$ to some other $\mathbb{Q}$-factorial variety with a projective
morphism to $S$.

In this paper, we mainly consider the absolute case, i.e. when $S$ is a single point. 

Assume Assumption \ref{Assumption main} holds for $X$, existence of a rational polyhedral fundamental domain for $\mathcal{M}^e(X)$ \emph{implies all the other statements} by Section \ref{Relations between cone conjectures}. So we mainly focus on existence of rational polyhedral fundamental domain for $\mathcal{M}^e(X)$. 

It is known that it makes no difference to assume $\Delta$ to be a $\mathbb{Q}$-divisor. Indeed, we may approximate $\Delta$ by $\mathbb{Q}$-divisors $\Delta^\prime$ such that $\text{Supp} \Delta = \text{Supp} \Delta^\prime$, $K_X+\Delta^\prime \equiv 0$, and the differences of coefficients are arbitrarily small. 

It is worth noting that the conjecture is no longer true for lc pairs. There exist counterexamples when $X$ is a lc surface with $K_X\equiv 0$ and $\Delta=0$ by \cite{Tot09}. 
\end{remark}

We briefly summarise some known results: 

In the relative setting, the weak cone conjecture for all 3-dimensional terminal Calabi-Yau fiber spaces over a nontrivial base is proved by Kawamata in \cite{Kaw97}. Finiteness of ${\rm PsAut}(X/S)$-equivalence classes for elliptic fibrations is proved by Filipazzi-Hacon-Svaldi in \cite{FHS21}. Recently, finiteness of ${\rm PsAut}(X/S)$-equivalence classes and existence of a rational polyhedral fundamental domain for $\mathcal{M}(X/S)_+$ for surface fibrations are proved by Li and Zhao in \cite{LZ22}, \cite{Li23}. 

When $S$ is a point, the cone conjecture is
known in dimension $2$ by Totaro in \cite{Tot09}, for abelian varieties by Prendergast-Smith in \cite{PS12a}. The cone conjecture is
known for Calabi-Yau manifolds with Picard number $2$ with infinite birational automorphism group by Lazi\'c-Peternell \cite{LP13} and Oguiso \cite{Ogu14}. Recently, the cone conjecture is proved for generalised hyperelliptic varieties by Monti-Quedo in \cite{MQ24}. 

Over an arbitrary field of characteristic $\neq 2$, the cone conjecture is known for K3 surfaces by Bright-Logan-van Luijk, see \cite{BLvL19}.

In the statement of the cone conjecture, sometimes $\mathcal{M}^e(X/S)$ is replaced by the rational hull $\mathcal{M}(X/S)_+$. In this case, the cone conjecture is also known for projective hyperk\"ahler manifolds over $\mathbb{C}$, see \cite{Mar11}, and it is generalised to any field of characteristic $0$ by Takamatsu in \cite{Tak25}. Recently, it is known for a special class of Enriques manifolds with prime index by Pacienza-Sarti, see \cite{PS23}, and the more general results on the cone conjecture of finite quotient is obtained by Gachet, see \cite{Gac25}. It is proved by Lutz in \cite{Lut24} that the cone conjecture for Calabi-Yau threefolds is preserved under deformation. 

Recently, equivalent formulations of the cone conjecture are studied systematically by Gachet-Lin-Stenger-Wang in \cite{GLSW24}. 

We recommand a good survey by Lazic–Oguiso–Peternell, see \cite{LOP18}.

In this paper, we mainly consider the case when the base $S$ is a point. We first prove the following reduction, see Section \ref{Section reduction}: 

\begin{theorem}\label{main thm} \label{main theorem 2}
Let $f: (X,\Delta_X) \to (Y,\Delta_Y)$ be a crepant birational morphism between $\mathbb{Q}$-factorial klt Calabi-Yau pairs over a field $k$ of characteristic $0$, such that the support of $\Delta_X$ contains the exceptional divisors. Suppose that Assumption \ref{Assumption main} or \ref{Assumption GLM} holds for $X$. 

Then the movable cone conjecture holds for $(X,\Delta_X)$ if and only if the movable cone conjecture holds for $(Y,\Delta_Y)$.

\end{theorem}

\begin{remark}
In this case, we can actually prove that the pseudo-automorphism groups ${\rm PsAut}(X,\Delta_X)$ and ${\rm PsAut}(Y,\Delta_Y)$ have a common finite index subgroup, and these cones are of polyhedral type with respect to this common finite index subgroup. This relation is no longer true when we drop the condition that the support of $\Delta_X$ contains the exceptional divisors. 
\end{remark}

In order to prove the cone conjecture via fibration structures, it turns out that we need the following variant of the cone conjecture in dimension $d\leq 2$, see Subsection \ref{Subsection generalised cone}: 

\begin{theorem}\label{Thm generalised cone}
Let $X$ be a projective klt Calabi-Yau variety of dimension $d\leq 2$ over a finitely generated field $k$ over $\mathbb{Q}$. The Picard group ${\rm Pic}(X)$ is known to be finitely generated. Let ${\rm Nef}^e(X)$ be the cone in ${\rm Pic}(X)\otimes_\mathbb{Z} \mathbb{R}$ given by $\{D \in {\rm Pic}(X) \otimes_\mathbb{Z} \mathbb{R}| D \text{ is nef and effective}\}$. 

Then there exists a closed rational polyhedral cone $P\subset {\rm Pic}(X)\otimes_\mathbb{Z} \mathbb{R}$ such that ${\rm Aut}(X) P = {\rm Nef}^e(X)$. 
\end{theorem}

\begin{remark}
    In general, ${\rm Nef}^e(X)$ is different from the usual nef cone $\mathcal{A}^e(X)$. See Remark \ref{remark difference} for an example. 
\end{remark}

Then we prove that this variant of the cone conjecture behaves well with respect to fibrations under certain conditions, see Subsection \ref{Sec fibration and cone}: 

\begin{theorem}\label{Lem Main}

Let $(X, \Delta)$ be a projective klt Calabi-Yau pair over a finitely generated extension of $\mathbb{Q}$, say $k$, where $X$ is $\mathbb{Q}$-factorial of dimenison $d$, and Assumption \ref{Assumption main} or \ref{Assumption GLM} holds for $X$. Suppose $-K_X$ is semiample, inducing a contraction $f: X \to S$, with $-K_X = f^*A$ for some ample $\mathbb{Q}$-divisor $A$ on $S$. Let $\eta$ be the generic point of $S$. 

Assume further that: 

(1) There exists a closed rational polyhedral cone $P\subset {\rm Pic}(X_\eta)\otimes_\mathbb{Z} \mathbb{R}$ such that ${\rm PsAut}(X_\eta)P = {\rm Mov}^e(X_\eta)$. 

(2) $G = {\rm PsAut}(X_\eta) \cap {\rm PsAut}(X,\Delta)$ is of essentially finite index in ${\rm PsAut}(X_\eta)$.

Then there exists a closed rational polyhedral cone $\Pi\subset {\rm Pic}(X)\otimes_\mathbb{Z} \mathbb{R}$ such that ${\rm PsAut}(X,\Delta)\Pi = {\rm Mov}^e(X)$. 

Furthermore, the cone conjecture holds for $(X,\Delta)$. 
\end{theorem}

\begin{remark}
See Definition \ref{Def cones in Pic} for the definition of ${\rm Mov}^e(X_\eta)$, and see Definition \ref{Def finiteness condition} for the definition of essentially finite index subgroups. 
\end{remark}

Combining these results, by running MMP and considering the anti-canonical fibration,  we prove our main result, see Subsection \ref{Subsection conclusion}: 

\begin{theorem}\label{The Main Theorem}\label{Thm main dimension 3}

Let $(X,\Delta)$ be a $\mathbb{Q}$-factorial klt Calabi-Yau pair of dimension $d$ over a field $k$ of characteristic $0$, and the Iitaka dimension $\kappa(X,-K_X)\geq d-2$. Then the cone conjecture holds for $(X,\Delta)$. 
\end{theorem}

\begin{remark}
    When $d = 2$, this is a generalisation of \cite{Tot09} to general base fields of characteristic $0$. 

    When $d = 3$, this generalises the results of \cite{PS12}, \cite{PS15} and \cite{Li23a}. 
\end{remark}

\begin{remark}
Our main strategy in this paper is running MMP and considering the anti-canonical fibration, and then we try to prove that the cone conjecture is preserved under these operations. We should note that such ideas already appeared in Totaro's result, see \cite{Tot09}. 

For running MMP, we can actually show that in each step the pseudo-automorphism groups are commensurable, and the complexity of the cones are essentially the same. So it suffices to consider the minimal model. 

In this part, we mainly use Shokurov polytope. We note that Shokurov polytope is known to be useful for this conjecture by Lazi\'c-Oguiso-Peternell in \cite{LOP20}, and it is used to prove the relative cone conjecture by Li-Zhao in \cite{LZ22}. The difference between \cite{LZ22} and our situation is that \cite{LZ22} works on the relative case, but we work on the absolute case. So we have to extract positivity from the boundary divisor, which is not needed in the relative case. 

For fibrations, we can still apply similar arguments, using Shokurov polytope. So the main ideas are parallel to the previous case, but we need to overcome some new difficulties. 

One of the new difficulties is that the cone conjecture over $\mathbb{C}$ does not behave well with respect to anti-canonical fibrations when $X$ is not rationally connected. In Totaro's paper \cite{Tot09}, this difficulty is avoided by a classification result. Our solution is to pass to non-algebraically closed fields (In fact, we finally work on finitely generated fields over $\mathbb{Q}$.) and to consider ${\rm Pic}(X)\otimes_\mathbb{Z} \mathbb{R}$ instead of $N^1(X)$. In this case, ${\rm Pic}(X)\otimes_\mathbb{Z} \mathbb{R}$ is still finite dimensional, and the kernel of the restriction map to the generic fiber $\ker({\rm Pic}(X)\otimes_\mathbb{Z} \mathbb{R} \to {\rm Pic}(X_\eta)\otimes_\mathbb{Z} \mathbb{R})$ is generated by vertical divisors by Lemma \ref{Lem Pic restriction}. 

Also, the lifting of pseudo-automorphisms from the generic fiber is not easy for klt singularities. Our solution is to consider the global index $1$ cover and we finally reduce to the terminal case, see Lemma \ref{Lem lifting 1}, Lemma \ref{Lem lifting 2}. 
\end{remark}

\begin{Ack}
The author would like to thank Professor Caucher Birkar for his constant support and guidance. He thanks Professor Mao Sheng for providing a family of examples, which motivates the arguments in this paper. He thanks Bingyi Chen for reading a draft of this paper and giving a lot of useful comments. He thanks Jia Jia, Xintong Jiang for useful discussions. He also thanks Long Wang for some suggestions on the organization of this paper. 

The author thanks the anonymous referee for useful suggestions that make this paper better organized. 
\end{Ack}

\section{Preliminaries}

\subsection{Notations and conventions}

Our conventions mainly follow from \cite{Kaw97}, but we allow any base field $k$ of characteristic $0$ (possibly non-algebraically closed). 

By a \emph{variety}, we always mean a geometrically connected integral scheme of finite type over a field, unless specified. 

Let $f:X \to S$ be a proper surjective
morphism of normal varieties over $k$, with geometrically connected fibers. (Essentially, throughout this article, we mainly consider the case when $S$ is a point.)

Let $D$ be a Cartier divisor on $X$. We recall some standard definitions for different kinds of divisors. It is said to be \emph{$f$-nef} if $(D \cdot C) \geq 0$ holds for any curve $C$ (curve means integral one-dimensional closed subscheme) on $X$ which is
contracted by $f$. It is said to be \emph{$f$-movable} if $\dim \text{Supp } \text{Coker}(f^*f_*\mathcal{O}_X(D) \to \mathcal{O}_X(D)) \geq 2$. It is said to be \emph{$f$-effective} if $f_*\mathcal{O}_X(D) \neq 0$. And it is said to be \emph{$f$-big} if its restriction to the generic fiber has maximal Kodaira dimension. 

The \emph{N\'eron-Severi group}
$$NS
(X/S) = \{\text{Cartier divisor on }X\}/(\text{numerical equivalence over }S)$$
is known to be finitely generated. Let $N^1
(X/S) = NS
(X/S)\otimes_{\mathbb{Z}} \mathbb{R}$. We set $\rho(X/S) = \dim N^1
(X/S)$ to be the \emph{relative Picard number}. The numerical class of an $\mathbb{R}$-Cartier
divisor $D$ in $N^1
(X/S)$ is denoted by $[D]$.

The \emph{$f$-pseudo-effective cone} is denoted by $\bar{\mathcal{B}}(X/S)$, the \emph{$f$-nef cone} is denoted by $\bar{\mathcal{A}}(X/S)$, and the \emph{closed $f$-movable cone} is denoted by $\bar{\mathcal{M}}(X/S)$. 
They are defined to be the closed convex cones in $N^1
(X/S)$ generated
by the numerical classes of $f$-effective
divisors, $f$-nef divisors, and $f$-movable divisors respectively. 


If $f$ is projective, the interior $\mathcal{B}(X/S) \subset \bar{\mathcal{B}}(X/S)$  (resp. $\mathcal{A}(X/S) \subset \bar{\mathcal{A}}(X/S)$) is the open convex subcone generated by the numerical classes of $f$-big divisors ($f$-ample divisors) and it is called
an \emph{$f$-big cone} (resp. \emph{$f$-ample cone}). Let $\mathcal{M}(X/S)\subset \bar{\mathcal{M}}(X/S)$ be the \emph{$f$-movable cone}, i.e. the subcone of $\bar{\mathcal{M}}(X/S)$ generated by $f$-movable divisors. 

We denote the \emph{$f$-effective cone} by $\mathcal{B}^e 
(X/S)$, i.e. the convex cone
generated by $f$-effective Cartier divisors. We call $\mathcal{A}^e
(X/S) = \bar{\mathcal{A}}(X/S) \cap \mathcal{B}^e
(X/S)$
and $\mathcal{M}^e
(X/S) = \bar{\mathcal{M}}(X/S) \cap \mathcal{B}^e
(X/S)$ the \emph{$f$-effective $f$-nef cone} and \emph{$f$-effective
$f$-movable cone}, respectively. 

We also introduce some notations for pairs. 

A \emph{pair} is a normal variety $X$ together with an $\mathbb{R}$-divisor $\Delta = \sum_i a_i\Delta_i$ such that $K_X+\Delta$ is $\mathbb{R}$-Cartier and $a_i\in [0,1]$. We adapt usual definition of singularities of pairs (lc, klt, canonical, terminal), see \cite[Definition 2.8]{Ko13}. 

Let $(X,\Delta)$ be a pair, $f:X\to S$ be a surjective morphism. Sometimes we use the notation $(X/S,\Delta)$ for a pair together with a surjective morphism $f:X\to S$. Let $g:T\to S$ be another morphism. We define the base change $(X,\Delta)\times_S T : = (X^\prime,\Delta^\prime)$, where $X^\prime = X\times_S T$, $\Delta^\prime_i = \Delta_i\times_S T$, $\Delta^\prime = \sum a_i\Delta^\prime_i$. 

We say a pair $(X/S,\Delta)$ is \emph{Calabi-Yau} if $K_X+\Delta \equiv_S 0$. When $S$ is a point, we often omit $S$ from the notation. When $S$ is not a point, $(X/S,\Delta)$ is called a \emph{Calabi-Yau fiber space}. 

Let $(X,\Delta_X)$, $(Y,\Delta_Y)$ be two pairs. An \emph{isomorphism} from $(X,\Delta_X)$ to $(Y,\Delta_Y)$ is an isomorphism $\phi: X\to Y$ of varieties such that $\phi_* \Delta_X= \Delta_Y$. A \emph{pseudo-isomorphism} from $(X,\Delta_X)$ to $(Y,\Delta_Y)$ is a small birational map  $\phi: X\dashrightarrow Y$ such that $\phi_* \Delta_X= \Delta_Y$. In a similar way we define \emph{automorphism}s and \emph{pseudo-automorphism}s. We denote the \emph{automorphism group} by ${\rm Aut}(X,\Delta)$, and the \emph{pseudo-automorphism group} by ${\rm PsAut}(X,\Delta)$. 

\subsection{Finite generation of Picard group}

We recall Lang and N\'eron's result on rational points in abelian varieties \cite[Theorem 1]{LN59}: 

\begin{theorem}\label{thm lang-neron}
Let $K$ be a finitely generated regular extension of a field $k$. 
Let $A$ be an abelian variety defined over $K$, and let $(B,\tau)$ be its $K/k$-trace.
Then $A(K)/\tau B(k)$ is a finitely generated abelian group. 
\end{theorem}

We have the following corollary, which is mentioned after \cite[Theorem 2]{LN59}: 

\begin{corollary}
Let $k$ be a finitely generated extension of $\mathbb{Q}$. 
Let $A$ be an abelian variety defined over $k$, then $A(k)$ is a finitely generated abelian group. 
\end{corollary}

\begin{proof}
Let $l$ be the algebraic closure of $\mathbb{Q}$ in $k$, i.e. the field consisting of all element in $k$ that are algebraic over $\mathbb{Q}$. Then $k/l$ is a regular extension. By Mordell-Weil theorem, see also \cite[Theorem 2]{LN59}, $B(l)$ is finitely generated, so we conclude that $A(k)$ is finitely generated for any abelian variety $A$ over $k$ by Theorem \ref{thm lang-neron}. 

\end{proof}

Let $X$ be a normal projective variety over $k$. Recall that we always assume $X$ is geometrically connected. Then the identity component of the Picard scheme $\mathcal{P}ic^0(X)$ is a proper group scheme by \cite[Theorem 9.5.4]{FAG}, hence an abelian variety. We always denote the Picard group by ${\rm Pic}(X)$, and denote the Picard scheme by $\mathcal{P}ic(X)$.

\begin{proposition}
Let $X$ be a normal projective variety over a finitely generated field $k$ over $\mathbb{Q}$. Then ${\rm Pic}^0(X)$ is finitely generated. 
\end{proposition}

\begin{proof}
Note that the Picard scheme $\mathcal{P}ic^0(X)$ represents the sheafification of the Picard functor in \'etale topology, so we do not have the equality $\mathcal{P}ic^0(X)(k) = {\rm Pic}^0(X)$. 

But we have a sheafification map ${\rm Pic}^0(X) \to \mathcal{P}ic^0(X)(k)$, which is injective by faithfully flat descent, see \cite[Theorem 9.2.5]{FAG}. So ${\rm Pic}^0(X)$ is finitely generated since it is a subgroup of a finitely generated abelian group. 

\end{proof}

\begin{corollary}
Let $X$ be a normal projective variety over a finitely generated field $k$ over $\mathbb{Q}$. Then ${\rm Pic}(X)$ is finitely generated. 
\end{corollary}

\begin{proof}
Recall that we have an exact sequence $0 \to {\rm Pic}^0(X) \to {\rm Pic}(X) \to NS(X) \to 0$. The N\'eron-Severi group $NS(X)$ is known to be finitely generated by theorem of the base. 

So ${\rm Pic}^0(X)$ is finitely generated since it is an extension of finitely generated abelian groups. 

\end{proof}

With this corollary, ${\rm Pic}(X)\otimes_\mathbb{Z} \mathbb{R}$ is known to be a finite dimensional vector space. 

We define some cones in ${\rm Pic}(X)\otimes_\mathbb{Z} \mathbb{R}$, which will be used in our main results: 

\begin{definition}\label{Def cones in Pic}
Let $X$ be a normal projective variety over a finitely generated field $k$ over $\mathbb{Q}$. 

We define \emph{the effective nef cone in ${\rm Pic}(X)\otimes_\mathbb{Z} \mathbb{R}$}, denoted by ${\rm Nef}^e(X)$, to be the cone in ${\rm Pic}(X)\otimes_\mathbb{Z} \mathbb{R}$ generated by effective nef $\mathbb{R}$-divisors. 

Similarly, we define \emph{the effective movable cone in ${\rm Pic}(X)\otimes_\mathbb{Z} \mathbb{R}$}, denoted by ${\rm Mov}^e(X)$, to be the cone in ${\rm Pic}(X)\otimes_\mathbb{Z} \mathbb{R}$ generated by effective $\mathbb{R}$-divisors $D$ such that its numerical class lies in $\mathcal{M}^e(X)$. 
    
\end{definition}

\begin{remark}\label{remark difference}
${\rm Nef}^e(X)$ (resp. ${\rm Mov}^e(X)$) is the same as $\mathcal{A}^e(X)$ (resp. $\mathcal{M}^e(X)$) when the identity component of the Picard group admits no non-torsion $k$-rational points. In particular, this holds when $h^1(X,\mathcal{O}_X) = 0$. 

In general, they can be different. A typical example is an elliptic curve $E$ with a zero point $O$ and a non-torsion $k$-rational point $P$. Then $[O]$ and $[P]$ give linearly independent elements in ${\rm Nef}^e(E)$ (resp. ${\rm Mov}^e(E)$), while $\mathcal{A}^e(E)$ (resp. $\mathcal{M}^e(E)$) is of dimension $1$. 
\end{remark}

\subsection{Some standard conjectures in MMP}

We recall some standard conjectures in MMP, and package an assumption which will be used for our main results. 

\begin{conjecture}
(Existence of minimal models) 
Let $(X/Z,\Delta)$ be a projective lc pair. If $K_X+\Delta$ is pseudo-effective over $Z$, then there exists a log minimal model $(X^\prime/Z,\Delta^\prime)$ of $(X/Z,\Delta)$. 
\end{conjecture}

\begin{remark}
    We refer to \cite[Subsection 2.1]{HH19} for a definition of log minimal model. 
\end{remark}
    
\begin{conjecture}(Non-vanishing) 
Let $(X/Z,\Delta)$ be a projective lc pair. If $K_X+\Delta$ is pseudo-effective, then $K_X+\Delta \sim_{\mathbb{R},Z} D$ for some effective $\mathbb{R}$-divisor $D$. 
\end{conjecture}

Throughout this paper, the following assumption will be used repeatedly. 

\begin{assumption}\label{Assumption main}
Let $X$ be a normal projective variety of dimension $d$ over a field $k$ of characteristic $0$. 

We say Assumption \ref{Assumption main} holds for $X$ if for any $\mathbb{R}$-divisor $B$ on $X$ such that $(X,B)$ is lc, 

(a) existence of minimal models holds for $(X,B)$, and 

(b) non-vanishing holds for $(X,B)$. 
\end{assumption}

\begin{proposition}\label{Prop Assumption}
Let $(X,\Delta)$ be a klt Calabi-Yau pair of dimension $d$ over a field $k$ of characteristic $0$. Suppose the Iitaka dimension $\kappa(X,-K_X)\geq d-3$.  

Then for any $\mathbb{R}$-divisor $B$ on $X$ such that $(X,B)$ is lc, existence of good minimal models holds for $(X,B)$. In particular, Assumption \ref{Assumption main} holds for $X$. 
\end{proposition}

\begin{proof}
We may assume $K_X+B$ is pseudo-effective, otherwise there is nothing to prove. Replacing $(X,B)$ by $(X,\frac{1}{2}(B+\Delta))$, we may assume $(X,B)$ is klt and ${\rm Supp}(\Delta) \subset {\rm Supp}(B)$. 

By \cite[Lemma 2.9.3]{HMX18}, we may assume $k$ is an algebraically closed field. 

Let $\epsilon$ be a small enough positive real number such that $(X,(1+\epsilon)\Delta)$ is klt. Let $\delta$ be a small enough positive real number such that we have $(1-\delta)(1+\epsilon)>1$. We may consider $B^\prime = \delta B+(1-\delta) \Delta$ and $B^{\prime\prime} = \delta B+(1-\delta)(1+\epsilon) \Delta$. 

Then we have $B^\prime\leq B^{\prime\prime}$ ${\rm Supp}(B^\prime) = {\rm Supp}(B^{\prime\prime})$.  
Moreover, we have $K_X+B^{\prime\prime} \sim_{\mathbb{R}} \delta B +((1-\delta)(1+\epsilon)-1)\Delta$, so $\kappa(X,K_X+B^{\prime\prime})\geq \kappa(X,-K_X)\geq d-3$. 
By \cite[Proposition 3.2]{CHL25}, $(X,B^{\prime\prime})$ admits a good minimal model. 
By \cite[Theorem 1.4]{MZ23}, $(X,B^\prime)$ admits a good minimal model. 

Since $K_X+B^\prime  = K_X+\delta B+(1-\delta) \Delta \sim_{\mathbb{R}}\delta(K_X+B)$, $(X,B)$ also admits a good minimal model. 

\end{proof}

\subsection{The geography of log models}

In this subsection, we work over an arbitrary field $k$ of characteristic $0$, unless specified. The main result in this subsection come from appendix $A$, and the ideas come from Choi's paper \cite{Choi08} and Shokurov-Choi's paper \cite{SC11}. We first state the conjecture of geography of log models (GLM). 

\begin{conjecture}\label{Conj GLM}
(\cite[Conjecture 3.2.3]{Choi08}, GLM) Fix a reduced divisor $F =
\sum_i D_i$ on a normal variety $X/Z$ of dimension $d$, where $Z$ is a normal affine variety. Let $B_F =
\bigoplus_i [0, 1]D_i$ and let $N_F$ be the subset
of divisors $B$ in $B_F$ such that $(X/Z, B)$ admits a wlc model.
Then the following hold:

(1) The set $N_F$ is closed convex rational polyhedral.

(2) There exists a finite decomposition $N_F = \bigsqcup_i \Pi_i^\circ$ such that each $\Pi_i$ is closed convex rational polyhedral, $\Pi_i^\circ$ is the relative interior of $\Pi_i$, and any divisors $D_1$ and $D_2$ in $\Pi_i$ admit the same wlc models. 
\end{conjecture}

\begin{remark}
    We refer to \cite[Definition 1.2.6]{Choi08} for a definition of weak log canonical models (i.e. wlc models). 
    We note that log minimal models are wlc models. 
\end{remark}

\cite{Choi08} proved that MMP implies GLM: 

\begin{theorem}(\cite[Main Theorem 1]{Choi08})
MMP in dimension $d$ implies GLM in dimension $d$. More
precisely, we can relax MMP to the following statement: for any $\mathbb{Q}$-factorial lc pair $(X/Z, B)$, where $X\to Z$ is projective and $Z$ is affine,
there exists a terminated sequence of extremal log flips, i.e. the sequence is finite and the last pair in the
sequence is either a $\mathbb{Q}$-factorial log minimal model or a Mori log fibration.

\end{theorem}

By Theorem \ref{Thm Galois descent}, Theorem \ref{Theorem A.4}, and the fact that Assumption \ref{Assumption main}(a) is preserved under change of base fields, we have the following more refined result: 

\begin{theorem}\label{mmp-GLM}
Let $X$ be a normal projective variety over a field $k$ of characteristic $0$.  Suppose Assumption \ref{Assumption main}(a) holds for $X$. Then Conjecture \ref{Conj GLM} holds for $X$. 
\end{theorem}

But Conjecture \ref{Conj GLM} is not enough for our purpose. We need the following version of GLM, which includes non-vanishing conjecture, and we call it ``strong GLM": 

\begin{conjecture}(Strong GLM)\label{Conj strong GLM}
Fix a reduced divisor $F =
\sum_i D_i$ on a normal variety $X/Z$, where $Z$ is a normal affine variety. Let $B_F =
\bigoplus_i [0, 1]D_i$ and let $N_F$ be the subset
of divisors $B$ in $B_F$ such that $(X,B)$ is lc and $K_X+B$ is $\mathbb{R}$-linearly equivalent to an effective $\mathbb{R}$-divisor over $Z$. 
Then the following hold:

(1) The set $N_F$ is closed convex rational polyhedral.

(2) There exists a finite decomposition $N_F = \bigsqcup_i \Pi_i^\circ$ such that each $\Pi_i$ is closed convex rational polyhedral, $\Pi_i^\circ$ is the relative interior of $\Pi_i$, and any divisors $D_1$ and $D_2$ in $\Pi_i$ admit the same wlc models. 
\end{conjecture}

We package strong GLM into an assumption: 

\begin{assumption}\label{Assumption GLM}
We say that Assumption \ref{Assumption GLM} holds for $X$ if strong GLM holds for $X$ with any reduced divisor $F$ on $X$. 
\end{assumption}

\begin{corollary}\label{strong GLM}
Fix a normal projective variety $X$. If Assumption \ref{Assumption main} holds for $X$, then Assumption \ref{Assumption GLM} holds for $X$.  
\end{corollary}

\begin{proof}

Let $B$ be any boundary divisor such that $(X,B)$ is lc. Since non-vanishing and existence of minimal models hold for $(X,B)$ by Assumption \ref{Assumption main}, we know that $K_X+B$ is $\mathbb{R}$-linearly equivalent to an effective divisor if and only if $(X,B)$ admits a weak log canonical model. 

Now since GLM holds for $X$ by Theorem \ref{mmp-GLM}, Assumption \ref{Assumption GLM} holds for $X$.

\end{proof}

\subsection{Polyhedral type cones}

By a polyhedral cone, we mean a closed cone generated by finitely many vectors, unless specified. 

Li-Zhao's result \cite[Theorem 1.3]{LZ22} essentially showed that for varieties over $\mathbb{C}$, the cone conjecture is equivalent to the existence of a (not necessarily rational) polyhedral cone $P\subset \mathcal{B}^e(X)$ such that $\mathcal{M}^e(X) \subset \bigcup_{g\in {\rm PsAut}(X,\Delta)}g_*P$. In this section, we generalize their result to more general fields, but restrict to rational polyhedral cones for convenience; see
 Theorem \ref{Thm equivalent cone}. 

As we need to work on more general fields, we have to trace back to Looijenga's result on polyhedral type cones \cite{Loo14}. 

More precisely, we recall

\begin{theorem}\label{Prop-Def}
(\cite[Proposition-Definition 4.1]{Loo14}) 
Let $V$ denote a real finite dimensional vector space
equipped with a rational structure $V(\mathbb{Q}) \subset V$ and $C$ is an open nondegenerate convex cone in $V$. Let $C_+$ be the rational hull of $C$, i.e. the cone generated by rational points in the closure of $C$. 

Let $\Gamma$ be a subgroup of ${\rm GL}(V)$ which stabilizes
$C$ and some lattice in $V(\mathbb{Q})$. Then the following conditions are equivalent:

$(i)$ There exists a polyhedral cone $\Pi$ in $C_+$ with $\Gamma \cdot \Pi = C_+$.

$(ii)$ There exists a polyhedral cone $\Pi$ in $C_+$ with $\Gamma \cdot \Pi \supset C$. 

(i)$^*$,(ii)$^*$ The corresponding property for the contragradient action of $\Gamma$ on the open dual cone ${\rm Int}(C^*)$.

If one of these
equivalent conditions is fulfilled, we say that $(V(\mathbb{Q}), C, \Gamma)$ is a polyhedral
triple or simply, that $(C_+, \Gamma )$ is of polyhedral type. 
\end{theorem}

\begin{remark}
Theorem \ref{Prop-Def} will be applied in the following two cases: 

(1) When $V = N^1(X)$ for some klt pair $(X,\Delta)$ such that $X$ is $\mathbb{Q}$-factorial, and $\Gamma$ is a subgroup of ${\rm PsAut}(X,\Delta)$. The cones will be $C = \mathcal{M}^e(X)$ or $\mathcal{A}^e(X)$. Then we may take the lattice $L$ to be the subgroup of $V$ generated by Weil divisors. The lattice $L$ is clearly preserved by ${\rm PsAut}(X,\Delta)$. 

(2) When $V = Pic(X) \otimes_{\mathbb{Z}}\mathbb{R}$ for some klt pair $(X,\Delta)$ over a finitely generated field $k$ over $\mathbb{Q}$ such that $X$ is $\mathbb{Q}$-factorial, and $\Gamma$ is a subgroup of ${\rm PsAut}(X,\Delta)$. The cones will be $C = {\rm Mov}^e(X)$ or ${\rm Nef}^e(X)$, as defined in Definition \ref{Def cones in Pic}. Then we may take the lattice $L$ to be the image of Weil divisors. The lattice $L$ is clearly preserved by ${\rm PsAut}(X,\Delta)$. 

We will omit the lattice in applications. 
\end{remark}

To see $L$ is indeed a lattice in the above two cases, we need the following lemma, which is suggested by Jia Jia: 

\begin{lemma}\label{Lemma finite gen of Weil div}
For a fixed klt pair $(X,\Delta)$ over a field $k$ of characteristic $0$, a fixed multiple of a $\mathbb{Q}$-Cartier Weil divisor on $X$ is Cartier. 
\end{lemma}

\begin{proof}

This is known 
when $k=\mathbb{C}$ by \cite[Theorem 1.9, Remark 1.10]{GKP16}, when $k=\bar{k}$ by \cite[Theorem 1.10]{HLQ20} (The latter in fact gives a uniform bound in any bounded family with bounded singularities.) For general base field, we note that a Weil divisor is Cartier if and only if
the ideal sheaf is invertible. So such a property satisfies faithfully flat descent, and we
conclude the general case. 

\end{proof}

\begin{theorem}\label{Thm equivalent cone}
Let $(X, \Delta)$ be a projective klt Calabi-Yau pair of over a field $k$ of characteristic $0$. Suppose Assumption \ref{Assumption main} or \ref{Assumption GLM} holds for $X$. Then the following are equivalent: 

(1) The cone conjecture holds for $(X,\Delta)$. 

(2) There exists a rational polyhedral cone $P \subset \mathcal{B}^e(X)$ such that
$$\bigcup_{g\in {\rm PsAut}(X,\Delta)}
g_* P \supset \mathcal{M}^e(X)$$.
\end{theorem}

\begin{proof}
It suffices to show (2) $\Rightarrow$ (1). By Corollary \ref{strong GLM}, the cone $\Pi = P\cap \mathcal{M}^e(X)$ is rational polyhedral. So $\bigcup_{g\in {\rm PsAut}(X,\Delta)}
g_* \Pi = \mathcal{M}^e(X)$. 

Observe that $\Pi = P\cap \mathcal{M}^e(X)\subset \mathcal{M}_+(X)$. Then by Theorem \ref{Prop-Def}, we necessarily have $\bigcup_{g\in {\rm PsAut}(X,\Delta)}
g_* \Pi = \mathcal{M}_+(X)$. But we have $\bigcup_{g\in {\rm PsAut}(X,\Delta)}
g_* \Pi =  \mathcal{M}^e(X)$, which implies $\bigcup_{g\in {\rm PsAut}(X,\Delta)}
g_* \Pi = \mathcal{M}^e(X)= \mathcal{M}_+(X)$. 

Then the cone conjecture holds by Lemma \ref{lemma existence of trivial stablizer} and \cite[Application 4.14]{Loo14}, see also \cite[Lemma 3.5]{LZ22} for an explanation. 

\end{proof}

\begin{lemma}\label{lemma existence of trivial stablizer}
Let $(V(\mathbb{Q}),C,\Gamma)$ be a polyhedral triple. Then there exists $\xi \in {\rm Int}(C^*) \cap  V(\mathbb{Q})^*$, such that the stabilizer of $\xi$ is trivial. 
    
\end{lemma}

\begin{proof}
By Theorem \ref{Prop-Def}, the definition of the polyhedral triple is preserved under taking dual. So it suffices to show that there exists some $\zeta \in C \cap  V(\mathbb{Q})$, such that the stabilizer of $\zeta$ is trivial. 

By \cite[Theorem 3.8]{Loo14}, $\Gamma$ acts properly discontinously on $C$. So for any polyhedral cone $P \subset C$ with $\dim P = \dim C$, the set $S = \{g\in G|g \neq id, gP\cap P \neq \emptyset\}$ is finite. Let $v_1, \dots, v_n \in P$ be a basis of $V(\mathbb{Q})$. Consider the vectors in $C\cap V(\mathbb{Q})$: $$\zeta_k = \sum_{i=1}^n k^iv_i, \text{ } k\geq 1$$

We claim that $\zeta_k$ has trivial stabilizer for some $k\geq 1$. If not, for each $k$, there exists some $g\in S$ such that $g\zeta_k = \zeta_k$. Since $S$ is finite, there exists some $g\in S$ such that $g\zeta_k = \zeta_k$ for infinitely many $k$. But this means $g=id$ since any $n$ different vertors among $\{\zeta_k\}_{k\geq 1}$ is linearly independent. This leads to a contradiction since $id \notin S$.

\end{proof}

\subsection{The behavior of the cones under the change of the base field}

We collect some results related to the behavior of the cones under the change of the base field.

\begin{lemma}\label{Lem curve}
Suppose we have the following diagram: 
\[\begin{tikzcd}
	{C_1} & {C_2} \\
	K & k
	\arrow[from=2-1, to=2-2]
	\arrow[from=1-1, to=2-1]
	\arrow["f", from=1-1, to=1-2]
	\arrow[from=1-2, to=2-2]
\end{tikzcd}\]
where $f: C_1 \to C_2$ is a finite surjective morphism of proper integral curves over  $K$ and $k$, such that $K/k$ is a finite extension, and $L$ is a line bundle on $C_2$. Then $L$ has degree $0$ on $C_2/k$ if and only if $f^*L$  has degree $0$ on $C_1/K$. 
\end{lemma}

\begin{proof}
By \cite[Chapter 7, Proposition 3.7, 3.8]{Liu02}, $\deg_K f^*L = ({\deg f}/{[K:k]})\deg_k L$.

\end{proof}

\begin{lemma}\label{Lemma cone of algebraic closure}
Let $X$ be a normal projective variety over a field $k$ of characteristic $0$, and $K$ be the algebraic closure of $k$. Let $X_K$ be the base change of $X$ to $K$. Then: 

(1) The pullback map $\tau^*: N^1(X) \to N^1(X_K)$ is well-defined and injective. 

(2) The image of $\tau^*$ is $N^1(X_K)^{{\rm Gal}(K/k)}$. 
\end{lemma}

\begin{proof}
(1) To show the map is well-defined, it suffices to prove that for any line bundle $L$ on $X$, if $L$ has zero intersection number with all integral curves (i.e. one dimensional integral closed subschemes) $C$ on $X$, then its pullback to $X_K$ has zero intersection number with all curves on $X_K$. 

Consider the commutative diagram: 
\[\begin{tikzcd}
	C & {X_K} & X \\
	& K & k
	\arrow[from=1-3, to=2-3]
	\arrow[from=2-2, to=2-3]
	\arrow[from=1-2, to=1-3]
	\arrow[from=1-2, to=2-2]
	\arrow["\ulcorner"{anchor=center, pos=0.125}, draw=none, from=1-2, to=2-3]
	\arrow[from=1-1, to=1-2]
	\arrow[from=1-1, to=2-2]
\end{tikzcd}\]
where $C$ is a curve on $X_K$. By passage to limit, see \cite[Proposition 1.10.9]{Fu11}, there exists a finite extension $k_0$ over $k$ such that we have the following Cartisian diagram (by taking the Galois closure, we may assume $k_0/k$ is Galois): 
\[\begin{tikzcd}
	C & {C_0} \\
	{X_K} & {X_{k_0}} & X \\
	K & {k_0} & k
	\arrow[from=2-3, to=3-3]
	\arrow[from=3-2, to=3-3]
	\arrow[from=2-2, to=2-3]
	\arrow[from=2-2, to=3-2]
	\arrow["\ulcorner"{anchor=center, pos=0.125}, draw=none, from=2-2, to=3-3]
	\arrow[from=2-1, to=2-2]
	\arrow[from=3-1, to=3-2]
	\arrow[from=2-1, to=3-1]
	\arrow[from=1-1, to=2-1]
	\arrow[from=1-2, to=2-2]
	\arrow[from=1-1, to=1-2]
	\arrow["\ulcorner"{anchor=center, pos=0.125}, draw=none, from=1-1, to=2-2]
	\arrow["\ulcorner"{anchor=center, pos=0.125}, draw=none, from=2-1, to=3-2]
\end{tikzcd}\]

Let $Z$ be the scheme-theoretic image of $C_0$ on $X$, then $Z$ is integral of dimension $1$ by universal properties. By our condition, $L$ has degree $0$ on $Z$. Observe that $C_0 \to X$ is finite, so $C_0 \to Z$ is finite. Then by Lemma \ref{Lem curve}, $L$ has degree $0$ on $C_0/k_0$. By \cite[Chapter 7, Proposition 3.7]{Liu02}, $\deg_{C/K}L=\deg_{C_0/k_0}L = 0$. 

To show injectivity, let $C$ be an integral curve on $X$. Then consider the base change diagram 
\[\begin{tikzcd}
	{C_K} & C \\
	{X_K} & X \\
	K & k
	\arrow[from=2-2, to=3-2]
	\arrow[from=3-1, to=3-2]
	\arrow[from=2-1, to=3-1]
	\arrow[from=2-1, to=2-2]
	\arrow[from=1-2, to=2-2]
	\arrow[from=1-1, to=1-2]
	\arrow[from=1-1, to=2-1]
	\arrow["\ulcorner"{anchor=center, pos=0.125}, draw=none, from=1-1, to=2-2]
	\arrow["\ulcorner"{anchor=center, pos=0.125}, draw=none, from=2-1, to=3-2]
\end{tikzcd}\]
By \cite[Chapter 7, Proposition 3.7]{Liu02}, $\deg_{C_0/k_0}L= \deg_{C/K}L = 0$.

(2) Let $[L]\in N^1(X_K)^{{\rm Gal}(K/k)}$, where $L\in {\rm Pic}(X_K)$. Then by passage to limit, see \cite[Proposition 1.10.2]{Fu11}, there exists a finite extension $k_0/k$ such that $L = f^*L_0$, where $L_0$ is a line bundle on $X_{k_0}$ and $f:X_K\to X_{k_0}$ is the natural map. By taking Galois closure, we may assume $k_0/k$ is Galois. Consider $L_0^{\prime} = \otimes_{g\in {\rm Gal}(k_0/k)}g^*L_0$. Then ${\rm Gal}(k_0/k)$ acts on $L_0^\prime$ and by faithfully flat descent $L_0^\prime = g^*M$, where $M$ is a line bundle on $X$ and $g:X_{k_0}\to X$ is the natural map. Now $[k_0:k][L] = [(g\circ f)^*M] = \tau^*[M]$.

\end{proof}

\begin{lemma} \label{Lem change of closed fields}
Let $X$ be a normal proper variety over a field $k$ of characteristic $0$, and $K$ be an extension of $k$. Assume $k$ and $K$ are both algebraically closed. Let $X_K$ be the base change of $X$ to $K$. Then the pullback map $\tau^*: N^1(X) \to N^1(X_K)$ is an isomorphism. 

\end{lemma}

\begin{proof}
By \cite[Theorem 9.6.3]{FAG}, the N\'eron-Severi group can be written as a quotient ${\rm Pic}(X)/{\rm Pic}^\tau(X)$, where ${\rm Pic}^\tau(X)$ denotes the quotient components. By base change properties of the Picard scheme, this is invariant under base change of algebraically closed fields. So $\tau^*$ is an isomorphism.

\end{proof}

\begin{lemma} \label{Lem change of closed fields}
Let $X$ be a normal proper variety over a field $k$ of characteristic $0$, and $K$ be an extension of $k$. Let $X_K$ be the base change of $X$ to $K$. Then the pullback map $\tau^*: N^1(X) \to N^1(X_K)$ is well-defined. 
\end{lemma}

\begin{proof}
We have the pullback for Picard groups $\tau^*_{\rm Pic}: {\rm Pic}(X) \to {\rm Pic}(X_K)$. 
To show that the induced map $\tau^*: N^1(X) \to N^1(X_K)$ well-defined, it suffices to show that if $L$ is a numerically trivial line bundle on $X$, then the base change $L_K$ is numerically trivial. 

Let $\bar{k}$ (resp. $\bar{K}$) be the algebraic closure of $k$ (resp. $K$). 

By Lemma \ref{Lemma cone of algebraic closure}(1), $L_{\bar{k}}$ is numerically trivial. By Lemma \ref{Lem change of closed fields},  $L_{\bar{K}}$ is numerically trivial. By Lemma \ref{Lemma cone of algebraic closure}(1) again, $L_K$ is numerically trivial. 

\end{proof}

\begin{lemma}\label{Lem cone up}
Let $(X,\Delta)$ be a projective pair over a field $k$ of characteristic $0$. Let $K$ be an extension of $k$, and $(X_K,\Delta_K) = (X,\Delta)\times_k K$. 
Suppose the natural map $\tau^*: N^1(X) \to N^1(X_K)$ is an isomorphism. 
Then:

(1) $\tau^*$ induces an isomorphism on the effective, ample, closed movable cones. 

(2) The image of ${\rm Aut}(X,\Delta)$ (resp. ${\rm PsAut}(X,\Delta)$) on ${\rm GL}(N^1(X))$ is contained in the image of ${\rm Aut}(X_K,\Delta_K)$ (resp. ${\rm PsAut}(X_K,\Delta_K)$) on ${\rm GL}(N^1(X_K))$ under the identification $\tau^*$. 

(3) The nef (resp. movable) cone conjecture holds for $(X_K,\Delta_K)$ if the nef (resp. movable) cone conjecture holds for $(X,\Delta)$. 

\end{lemma}

\begin{proof}
(1) Let $\tau: X_K \to X$ be the natural map. 
Let $L$ be a line bundle on $X$, $L_K$ be the base change to $X_K$. Note that in all these cases injectivity is clear. 

\emph{Effective cone}: By flat base change, $L$ admits a section if and only if $L_K$ admits a section. So the induced map is well-defined. 

On the other hand, let $E_K$ be an effective divisor on $X_K$
. By passage to limit, there exists a finitely generated $k$-subalgebra $R$ of $K$ such that $E_K$ is defined on $X_R$, i.e. there exists a divisor $E_R$ on $X_R$ such that the pullback to $X_K$ is $E_K$. 
Let $x$ be a closed $k$-point of $\text{Spec}(R)$, and $L_x$ be the restriction of $L_R$ to $X_{x} =  \{x\}\times_k X \subset X_R$. 
By upper semicontinuity, we have $h^0(X_{x},L_x)>0$. 

In this case, $L_x$ is effective, and its base change to $K$ is algebraically equivalent to $\mathcal{O}_{X_K}(E_K)$. So surjectivity holds. 

\emph{Ample cone}: 
$L$ is ample if and only if $L_K$ is ample since ampleness is invariant under base change of fields. 

The "only if" part means the induced map is well-defined. 
Note that being ample is invariant under numerical equivalence. 
So the "if" part implies Surjectivity. 

\emph{Closed movable cone}: By definition, $L$ being in the closed movable cone is equivalent to the following condition: For any ample $\mathbb{Q}$-divisor $A$, there exists a positive integer $N$ such that $L^N(NA)$ is generated by global sections up to a codimension $2$ set. But this condition is invariant under base change, so $L$ lies in the closed movable cone if and only if $L_K$ does. 

The "only if" part means the induced map is well-defined. 
Note that being in the closed movable cone is invariant under numerical equivalence. 
So the "if" part implies Surjectivity. 

(2) This can be deduced from the following commutative diagram: 
\[\begin{tikzcd}
	{{\rm Aut}(X,\Delta)} & {{\rm PsAut}(X,\Delta)} & {{\rm GL}(N^1(X))} \\
	{{\rm Aut}(X_K,\Delta_K)} & {{\rm PsAut}(X_K,\Delta_K)} & {{\rm GL}(N^1(X_K))}
	\arrow[from=1-2, to=1-3]
	\arrow[from=2-2, to=2-3]
	\arrow[from=1-2, to=2-2]
	\arrow["\simeq", from=1-3, to=2-3]
	\arrow[from=1-1, to=1-2]
	\arrow[from=1-1, to=2-1]
	\arrow[from=2-1, to=2-2]
\end{tikzcd}\]

(3) This result follows from (2) and  
Theorem \ref{Thm equivalent cone}.

\end{proof}

\begin{lemma}\label{Lem passage to limit}
Let $(X_K,\Delta_K)$ be a projective pair over a field $K$ of characteristic $0$. Then there exist a subfield $k$ of $K$, finitely generated
 over $\mathbb{Q}$, and a projective pair $(X,\Delta)$ over $k$ such that $(X,\Delta)\times_k K = (X_K,\Delta_K)$. 

Moreover, if for any finitely generated subfield $l$ of $K$ over $k$, the nef (resp. movable) cone conjecture holds for $(X_l,\Delta_l) : = (X,\Delta)\times_k l$, then the nef (resp. movable) cone conjecture holds for $(X_K,\Delta_K)$. 
\end{lemma}

\begin{proof}
By passage to limit, see \cite[Proposition 1.10.9]{Fu11}, there exists a finitely generated subfield $k$ and a pair $(X,\Delta)$ over $k$ such that $(X,\Delta)\times_k K = (X_K,\Delta_K)$. By faithfully flat descent of properness, see \cite[Corollery 1.7.12]{Fu11}, $X$ is proper over $k$. 
By \cite[Corollaire 6.6.5]{EGAII}, $X$ admits an ample divisor, so $X$ is projective over $k$. 

Pick line bundles $L_1,\dots,L_n$ on $X_K$ such that their images span $N^1(X_K)$ as a vector space. By passage to limit, we may choose a finitely generated extension $l/k$ such that $L_1,\dots,L_n$ are defined over $l$. By the fact that $[L_1],\dots,[L_n]$ are Galois invariant, the pullback map $\tau^*: N^1(X_l) \to N^1(X_K)$ is surjective. On the other hand, consider the commutative diagram: 
\[\begin{tikzcd}
	{N^1(X_l)} & {N^1(X_{\bar{l}})} \\
	{N^1(X_K)} & {N^1(X_{\bar{K}})}
	\arrow["{\tau^*}", from=1-1, to=2-1]
	\arrow[from=1-1, to=1-2]
	\arrow["\simeq", from=1-2, to=2-2]
	\arrow[from=2-1, to=2-2]
\end{tikzcd}\]
By Lemma \ref{Lemma cone of algebraic closure}, the horizontal maps are injective, so $\tau^*$ is injective. 

So we conclude that $\tau^*$ is an isomorphism, and the result follows from Lemma \ref{Lem cone up}.

\end{proof}

\section{Relations between nef and movable cone conjectures}

\label{Relations between cone conjectures}

First we note that this is a formal consequence of Li-Zhao's result \cite[Theorem 1.3 (3)]{LZ22} if we work over $\mathbb{C}$. The idea of the proof is essentially the same.  

\begin{theorem}\label{relation between nef and movable cone conjecture}

Let $(X,\Delta)$ be a projective klt Calabi-Yau pair over a field $k$ of characteristic $0$, where $X$ is $\mathbb{Q}$-factorial. Suppose Assumption \ref{Assumption main} or \ref{Assumption GLM} holds for $X$. Then the movable cone conjecture implies the nef cone conjecture.
\end{theorem}

\begin{proof}
By Corollary  \ref{strong GLM}, we may assume Assumption \ref{Assumption GLM} holds for $X$. 

Let $\Pi$ be a rational polyhedral fundamental domain for the action of ${\rm PsAut}(X,\Delta)$ on $\mathcal{M}^e(X)$. We have an equality $\mathcal{M}^e(X) = \bigcup_{g\in {\rm PsAut}(X,\Delta)} g_*\Pi$. Taking intersection with $\mathcal{A}^e(X)$, we have  $\mathcal{A}^e(X) = \bigcup_{g\in {\rm PsAut}(X,\Delta)} (g_*\Pi\cap \mathcal{A}^e(X))$. 

By Geography of log models, see Assumption  \ref{Assumption GLM}, $\Pi$ admit a decomposition into finitely many (not necessarily closed) rational polyhedral cones $\Pi = \bigcup_{i=1}^{n} \Pi_i$
such
that whenever
$B$, $D$ are effective divisors with $[B], [D] \in \Pi_i$
and $\delta \in \mathbb{R}_{>0}$ is sufficient small such that 
$(X, \Delta + \delta B)$, $(X, \Delta + \delta D)$ are lc,
then if $(Y, \Delta_Y +\delta B_Y )$ is a weak log canonical model of $(X, \Delta + \delta B)$, then $(Y, \Delta_Y +\delta D_Y )$ is a weak log canonical model of $(X, \Delta + \delta D)$. 

In particular, each $\Pi_i$ lies in a single chamber $g_*\mathcal{A}^e(X)$ for some $g\in {\rm PsAut}(X,\Delta)$. Let $\bar{\Pi}_i$ be the closure of $\Pi_i$, then $\bar{\Pi}_i$ still lies in a single chamber $g_*\mathcal{A}^e(X)$ for some $g\in {\rm PsAut}(X,\Delta)$, and we may assume $\{\bar{\Pi}_i\}_{i=1}^{m}$ are the cones among $\{\bar{\Pi}_i\}_{i=1}^{n}$ with maximal dimension $\rho(X)$. Then it is clear that $\Pi = \bigcup_{i=1}^{m} \bar{\Pi}_i$. 

Now we have equalities 

\begin{align*}
\mathcal{A}^e(X) &= \bigcup_{g\in {\rm PsAut}(X,\Delta)} g_*(\Pi\cap g^{-1}_*\mathcal{A}^e(X)) 
\\ &= \bigcup_{g\in {\rm PsAut}(X,\Delta)} g_*(\bigcup_{i=1}^{m} \bar{\Pi}_i \cap g^{-1}_*\mathcal{A}^e(X))
\\ &= \bigcup_{g\in {\rm PsAut}(X,\Delta)} g_*(\bigcup_{\bar{\Pi}_i \subset g^{-1}_*\mathcal{A}^e(X)} \bar{\Pi}_i )
\\ &= \bigcup_{g\in {\rm PsAut}(X,\Delta)} \bigcup_{\bar{\Pi}_i \subset g^{-1}_*\mathcal{A}^e(X)} g_*\bar{\Pi}_i 
\\ &= \bigcup_{g\in {\rm PsAut}(X,\Delta)} \bigcup_{g_*\bar{\Pi}_i \subset \mathcal{A}^e(X)} g_*\bar{\Pi}_i 
\\ \hfill
\end{align*}

We claim that such $g_*\bar{\Pi}_i $ fall into finitely many orbits under the action of ${\rm Aut}(X,\Delta)$. In fact, if $g_*\bar{\Pi}_i \subset \mathcal{A}^e(X)$, $h_*\bar{\Pi}_i \subset \mathcal{A}^e(X)$, since $g_*\bar{\Pi}_i$ is of dimension $\rho(X)$, we may pick any ample divisor $H$ on $g_*\bar{\Pi}_i$. Then $(h\circ g^{-1})_*[H] \in \mathcal{A}^e(X)$ is nef. Now we have a commutative diagram

\[\begin{tikzcd}
	& W \\
	{(X,H)} && {(X,(h\circ g^{-1})_*H)}
	\arrow["{h\circ g^{-1}}", dashed, from=2-1, to=2-3]
	\arrow["{\pi_1}"', from=1-2, to=2-1]
	\arrow["{\pi_2}", from=1-2, to=2-3]
\end{tikzcd}\]

We claim that  $\pi_1^*H = \pi_2^*((h\circ g^{-1})_*H)$. In fact, by negativity lemma for $\pi_1$, $\pi_1^*H - \pi_2^*((h\circ g^{-1})_*H) \geq 0$, while by negativity lemma for $\pi_2$, $\pi_1^*H - \pi_2^*((h\circ g^{-1})_*H) \leq 0$. 

Now all curves contracted by $\pi_2$ must be contracted by $\pi_1$, so $h\circ g^{-1}$ is a morphism. Interchanging $g$ and $h$, we can show $g\circ h^{-1}$ is also a morphism.

So $g\circ h^{-1}\in {\rm Aut}(X,\Delta)$, which means $g_*\bar{\Pi}_i$ and $h_*\bar{\Pi}_i$ lie in the same orbit, so such $g_*\bar{\Pi}_i$ falls into at most $m$ orbits under the action of ${\rm Aut}(X,\Delta)$.

\end{proof}

Moreover, we explain something about relations between the cone conjecture and the weak cone conjecture. See Remark \ref{remark cone} for our convention on weak cone conjectures. 

The following lemma is well-known, see also \cite{Tot09}: 

\begin{lemma}
For klt Calabi-Yau pairs $(X/S,\Delta)$, the nef cone conjecture implies the weak nef cone conjecture. 
\end{lemma}

\begin{proof}
Let $\Pi$ be a rational polyhedral fundamental domain for the action of ${\rm Aut}(X/S,\Delta)$ in $\mathcal{A}^e(X/S)$. Let $X\to Y$ be a contraction of projective varieties over $S$ given by
some semiample line bundle $L$ on $X/S$. 

By definition, there exists an element $g\in {\rm Aut}(X/S,\Delta)$ such that $L$ lies in $g_*\Pi$. Suppose $g_*^{-1}L$ lies in the interior of a face $F$ of $\Pi$. 

Now we have a map (depending on some choices) from the set of contractions starting from $X$ over $S$ (modulo ${\rm Aut}(X/S,\Delta)$) to the set of faces to $\Pi$. We claim that this map is injective. 

Indeed, suppose two contractions $X\to Y_1$ and $X \to Y_2$ over $S$ correspond to the same face $F$ of $\Pi$. Suppose $X\to Y_i$ correspond to semiample divisors $L_i$. Since $L_1$ and $L_2$ lie in the interior of the same face of $\Pi$, they lie in the interior of the same face of $\mathcal{A}^e(X/S)$. Hence these two contractions are the same up to ${\rm Aut}(X/S,\Delta)$.

\end{proof}

The following lemma is similar to Theorem \ref{relation between nef and movable cone conjecture}, see also \cite[Proposition 5.3(2)]{LZ22}: 

\begin{lemma}

Let $(X,\Delta)$ be a projective klt Calabi-Yau pair over a field $k$ of characteristic $0$, where $X$ is $\mathbb{Q}$-factorial. Suppose Assumption \ref{Assumption main} or \ref{Assumption GLM} holds for $X$. Then the movable cone conjecture implies the weak movable cone conjecture.
\end{lemma}

\begin{proof}
By Corollary  \ref{strong GLM}, we may assume Assumption \ref{Assumption GLM} holds for $X$. 

Let $\Pi$ be a rational polyhedral fundamental domain for the action of ${\rm PsAut}(X,\Delta)$ on $\mathcal{M}^e(X)$. We have an equality $\mathcal{M}^e(X) = \bigcup_{g\in {\rm PsAut}(X,\Delta)} g_*\Pi$.

By Geography of log models, see Assumption \ref{Assumption GLM}, $\Pi$ admit a decomposition into finitely many (not necessarily closed) rational polyhedral cones $\Pi = \bigcup_{i=1}^{n} \Pi_i$
such
that whenever
$B$, $D$ are effective divisors with $[B], [D] \in \Pi_i$
and $\delta \in \mathbb{R}_{>0}$ is sufficient small such that 
$(X, \Delta + \delta B)$, $(X, \Delta + \delta D)$ are lc,
then if $(Y, \Delta_Y +\delta B_Y )$ is a weak log canonical model of $(X, \Delta + \delta B)$, then $(Y, \Delta_Y +\delta D_Y )$ is a weak log canonical model of $(X, \Delta + \delta D)$. 

Let $\mathcal{A}^e(X_0,\alpha)$ be a chamber corresponding to a small $\mathbb{Q}$-factorial modification $\alpha: X_0 \to X$. Then there exists some $g\in {\rm PsAut}(X,\Delta)$ such that $g_*\mathcal{A}(X_0,\alpha) = \mathcal{A}(X_0,g\circ \alpha)$ intersects $\Pi$. Then $g_*\mathcal{A}(X_0,\alpha)\cap \Pi$ is a union of wlc equivalence classes. In particular, $g_*\mathcal{A}(X_0,\alpha)$ contains $\Pi_i$ for some $i$. 

Now we have a map (depending on some choices) from the set of small $\mathbb{Q}$-factorial modification $\alpha: X_0 \dashrightarrow X$ (modulo ${\rm PsAut}(X/S, \Delta)$) to the set $\{\Pi_i\}_{1\leq i\leq n}$. We claim that this map
is injective. 

Indeed, suppose two small $\mathbb{Q}$-factorial modifications $\alpha_i: X_i \dashrightarrow X$, $i=1,2$, correspond to the
same $\Pi_i$. By definition, there exists some $g_i\in {\rm PsAut}(X,\Delta)$, $i=1,2$, such that $g_{i*}\mathcal{A}(X_i,\alpha_i) = \mathcal{A}(X_i,g_i\circ \alpha_i)$, $i=1,2$, intersect $\Pi_i$. Then there exists a divisor $A$ such that the strict transforms on $X_i$ are ample. This means $(g_1\circ \alpha_1)^{-1}g_2\circ \alpha_2$ is an isomorphism, so $\alpha_1$ and $\alpha_2$ are the same up to ${\rm PsAut}(X,\Delta)$.

\end{proof}

\section{Proof of Theorem \ref{main thm}}
\label{Section reduction}

In this section, we prove Theorem \ref{main thm}.

\subsection{The pseudo-automorphisms}

Note that ${\rm PsAut}(X,\Delta_X)$ and ${\rm PsAut}(Y,\Delta_Y)$ can be identified as subgroups of ${\rm Bir}(X)$. 

In this part, we prove that ${\rm PsAut}(X,\Delta_X)$ and ${\rm PsAut}(Y,\Delta_Y)$ are the same up to finite index. More precisely, we prove: 

\begin{proposition} \label{Lifting birational}Notations as in Theorem \ref{main theorem 2}. 
$G = {\rm PsAut}(X,\Delta_X)\cap {\rm PsAut}(Y,\Delta_Y)$ is of finite index in each of these two groups, where the intersection is taken inside ${\rm Bir}(X)$. 
\end{proposition}

\begin{proof}

We note that it suffices to prove when $k$ is algebraically closed. Let $\bar{k}$ be the algebraic closure of $k$, and denote $(\bar{X},\Delta_{\bar{X}}) = (X,\Delta_X) \times_k \bar{k}$, $(\bar{Y},\Delta_{\bar{Y}}) = (Y,\Delta_Y) \times_k \bar{k}$ and $\bar{G} = {\rm PsAut}(\bar{X},\Delta_{\bar{X}})\cap {\rm PsAut}(\bar{Y},\Delta_{\bar{Y}})$. Indeed, we have the following commutative diagram, where each face is a pullback diagram: 
\[\begin{tikzcd}
	& {\bar{G}} && {{\rm PsAut}(\bar{X},\Delta_{\bar{X}})} \\
	G && {{\rm PsAut}(X,\Delta_X)} \\
	& {{\rm PsAut}(\bar{Y},\Delta_{\bar{Y}})} && {{\rm Bir}(\bar{X})} \\
	{{\rm PsAut}(Y,\Delta_X)} && {{\rm Bir}(X)}
	\arrow[from=2-3, to=4-3]
	\arrow[from=4-3, to=3-4]
	\arrow["i"'{pos=0.3}, from=2-1, to=4-1]
	\arrow["j"'{pos=0.2}, from=2-1, to=2-3]
	\arrow[from=4-1, to=4-3]
	\arrow[from=4-1, to=3-2]
	\arrow[from=3-2, to=3-4]
	\arrow[from=2-3, to=4-3]
	\arrow["{\bar{i}}"'{pos=0.3}, from=1-2, to=3-2]
	\arrow[from=1-4, to=3-4]
	\arrow[from=2-3, to=1-4]
	\arrow["{\bar{j}}"'{pos=0.2}, from=1-2, to=1-4]
	\arrow[from=2-1, to=1-2]
\end{tikzcd}\]
If $\bar{i}$, $\bar{j}$ are inclusions of subgroups of finite index, so are $i$, $j$. So we may assume $k$ is algebraically closed in the following.

Firstly, we show that $G$ is of finite index in ${\rm PsAut}(Y,\Delta_Y)$:

Set $H = {\rm PsAut}(Y,\Delta_Y)$. Let $g: (W,\Delta_W) \to (X,\Delta_X)$ be a terminal model of $(X,\Delta_X)$. Then $f\circ g:(W,\Delta_W) \to (Y,\Delta_Y)$ is a terminal model of $(Y,\Delta_Y)$. 

Let $\sigma$ be an element in ${\rm PsAut}(Y,\Delta_Y)$. Consider the following commutative diagram, where $Z$ is a common resolution: 
\[\begin{tikzcd}
	& Z \\
	W && W \\
	X && X \\
	Y && Y
	\arrow["f", from=3-1, to=4-1]
	\arrow["f", from=3-3, to=4-3]
	\arrow["g", from=2-1, to=3-1]
	\arrow["g", from=2-3, to=3-3]
	\arrow["\sigma", dashed, from=4-1, to=4-3]
	\arrow["{\sigma_X}", dashed, from=3-1, to=3-3]
	\arrow["{\sigma_W}", dashed, from=2-1, to=2-3]
	\arrow["{\pi_1}"', from=1-2, to=2-1]
	\arrow["{\pi_2}", from=1-2, to=2-3]
\end{tikzcd}\]

Note that $\pi_1^*(K_W+\Delta_W) \equiv \pi_2^*(K_W+\Delta_W) \equiv 0$ by the Calabi-Yau condition. We set $\Delta_{Z_i}$ ($i=1,2$) as the crepant pullback of $\Delta_W$ under $\pi_i$ ($i=1,2$), i.e. $$K_Z+\Delta_{Z_1}=\pi_1^*(K_W+\Delta_W), K_Z+\Delta_{Z_2}=  \pi_2^*(K_W+\Delta_W)$$ Then $\Delta_{Z_1}-\Delta_{Z_2}$ is supported in the exceptional locus of $Z/Y$. (Note that this exceptional locus of $f\circ g\circ \pi_1$ and $f\circ g\circ \pi_2$ are the same.) Then by negativity lemma, we have $\Delta_{Z_1}-\Delta_{Z_2} \geq 0$ and $\Delta_{Z_1}-\Delta_{Z_2} \leq 0$. So $\Delta_{Z_1}=\Delta_{Z_2}$. Denote $\Delta_Z := \Delta_{Z_1}=\Delta_{Z_2}$.

This means the contracted divisors of $\pi_1$ and $\pi_2$ are precisely the divisors in $\Delta_Z$ with negative coefficients. In other words, this means $\sigma_W$ is small. 

So $H$ also lies in  ${\rm PsAut}(W)$. 

This group $H$ acts on the exceptional divisors of $W/Y$, hence a finite index subgroup $H^\prime$ acts by trivial permutation. Then this subgroup lies in  ${\rm PsAut}(X)$. Moreover, the exceptional divisors of $X/Y$ is preserved by $H^\prime$, so $\Delta_X$ is preserved by $H^\prime$. Hence $H^\prime$ lies in ${\rm PsAut}(X,\Delta_X)$.

Now $H^\prime\subset G \subset {\rm PsAut}(Y,\Delta_Y)$, so $G$ is of finite index in ${\rm PsAut}(Y,\Delta_Y)$. 

Secondly, we show that $G$ is of finite index in ${\rm PsAut}(X,\Delta_X)$: 

In fact, ${\rm PsAut}(X,\Delta_X)$ acts on the components of the support of $\Delta_X$, hence a finite index subgroup $H^{\prime\prime}$ acts on them by identity permutation. In particular, $H^{\prime\prime}$ preserves the exceptional divisors of $X/Y$, so $H^{\prime\prime}$ lies in ${\rm PsAut}(Y,\Delta_Y)$. 

Now $H^{\prime\prime}\subset G \subset {\rm PsAut}(X,\Delta_X)$, so $G$ is of finite index in ${\rm PsAut}(X,\Delta_X)$.

\end{proof}

\begin{remark}
From the proof, we may observe that the condition that the support of $\Delta_X$ contains the exceptional divisors is only used in the second part of the proof. 

In fact, the second part may fail without this condition. There are many counterexamples, see for example \cite{PS12}. 
\end{remark}

\subsection{Proof of Theorem \ref{main theorem 2}}


\begin{proof}(of Theorem \ref{main theorem 2})

Suppose the cone conjecture holds for $(X,\Delta_X)$. Let $\Pi$ be a rational polyhedral fundamental domain of $\mathcal{M}^e(X)$ for the action of ${\rm PsAut}(X,\Delta_X)$. Let ${\rm PsAut}(X,\Delta_X) = \sqcup_i G\sigma_i$ be a right coset decomposition. Let $P$ be the cone generated by $\sigma_i\Pi$. Since $G$ is of finite index, $P$ is a rational polyhedral cone. 

Now $f_*P$ is a polyhedral cone inside $\mathcal{M}^e(Y)$. We claim that $f_*P$ guarantees the cone conjecture for $(Y,\Delta_Y)$ in the sense of Theorem \ref{Thm equivalent cone}, i.e. $\bigcup_{g\in {\rm PsAut}(Y,\Delta_Y)}
g_* f_*P \supset \mathcal{M}^e(Y)$. Note that $$ \bigcup_{g\in {\rm PsAut}(Y,\Delta_Y)}
g_* f_*P \supset \bigcup_{g\in G}
g_* f_*P = f_*(\bigcup_{g\in G}
g_* P) = f_*\mathcal{M}^e(X)$$ Now it suffices to show $f_*: \mathcal{M}^e(X) \to \mathcal{M}^e(Y)$ is surjective. 

This follows from \cite[Chapter 3 Proposition 1.8, 1.14]{Nak04}. Indeed, for any $B\in \mathcal{M}^e(X)$, we have $$f^*B-\sum_{D} \sigma_{D}(f^*B) \in \mathcal{M}^e(Y) \text{ and } f_*(f^*B-\sum_{D} \sigma_{D}(f^*B)) = B$$ 

Thus the cone conjecture holds for $(Y,\Delta_Y)$. 

On the other hand, suppose the cone conjecture holds for $(Y,\Delta_Y)$. Let $\Pi$ be a rational polyhedral fundamental domain of $\mathcal{M}^e(Y)$. Let ${\rm PsAut}(Y,\Delta_Y) = \sqcup_i G\sigma_i$ be a right coset decomposition. Let $P$ be the cone generated by $\sigma_i\Pi$. 

Consider $f_*: \mathcal{M}^e(X) \to \mathcal{M}^e(Y)$. We claim that $f_*^{-1}(P)$ is (closed) rational polyhedral. 

Let $\{[G_i]\}_{i\in I}$ be a finite set of generators of $P$. Let $\{D_j\}_{j\in J}$ be all prime divisors appearing in $\{G_i\}_{i\in I}$. Let $\tilde{D}_j$ be the strict transform of $D_j$. Write $\Delta = \sum_{i=1}^n d_i\Delta_i$ as a sum of prime divisors. Let $S = \{\tilde{D}_j\}_{j\in J} \cup \{\Delta_i\}_{1\leq i\leq n}$. 

Let $F = \sum_{D\in S} D$, and recall $B_F = \oplus_{D\in S}[0,1]D$. Write $$M_F = \{D\in B_F| K_X + D\in \mathcal{M}^e(X)\}$$
Now the subset $M_F$ of $B_F$  is precisely the union of the wlc equivalence classes that admit small wlc models, say $M_F = \sqcup_i C_i$. 

Hence by strong GLM, $\bar{M}_F$ is closed rational polyhedral. Moreover, we claim that $\bar{M}_F = M_F$. In fact, by non-vanishing, all divisors in $K_X+\bar{M}_F$ are effective, and they clearly lie in $\bar{\mathcal{M}}(X)$. 

We fix notations by the following commutative diagram: 
\[\begin{tikzcd}
	{\mathcal{M}^e(X)} & {N^1(X)} \\
	{\mathcal{M}^e(Y)} & {N^1(Y)}
	\arrow[hook, from=1-1, to=1-2]
	\arrow["{f_*}"', from=1-1, to=2-1]
	\arrow["{\bar{f}_*}", from=1-2, to=2-2]
	\arrow[hook, from=2-1, to=2-2]
\end{tikzcd}\]

Then by definition we have $f_*^{-1}(P) = \bar{f}_*^{-1}(P) \cap \mathcal{M}^e(X)$. 

Let $C$ be the cone in $N^1(X)$ generated by the image of $K_X+\bar{M}_F$. Then $C$ is generated by finitely many divisors as $\bar{M}_F$ is a rational polytope. We claim that $f_*^{-1}(P) = \bar{f}_*^{-1}(P) \cap C$. Clearly RHS is contained in LHS since $C\subset \mathcal{M}^e(X)$. 

For the converse inclusion, let $[D]$ be any element in $f_*^{-1}(P)$, where $D$ is an effective divisor on $X$. Then $f_*[D]\in P$, so we may write $[f_*D] = \sum_{j\in J} a_j[D_j]$ with $a_j\geq 0$ since $P$ is contained in the cone generated by $\{[D_j]\}_{j\in J}$. So 
$$[D]= \sum_{j\in J} a_j[\tilde{D}_j] + \sum b_i E_i = K_X+\Delta_X+ \sum_{j\in J} a_j[\tilde{D}_j] + \sum b_i [E_i]$$ 
where $E_i$ are exceptional divisors. So by multiplying $\varepsilon$, we have $$\varepsilon [D]= \sum_{j\in J} \varepsilon a_j[\tilde{D}_j] + [K_X+\Delta_X]+ \varepsilon \sum b_i [E_i]$$
When $\varepsilon$ is sufficiently small, all coefficients in $\sum_{j\in J} \varepsilon a_j\tilde{D}_j + \Delta_X+ \varepsilon \sum b_i E_i$ lies in $[0,1]$, which means $\varepsilon [D]$ lies in the image of $K_X+\bar{M}_F$. Hence $[D]\in C$. So the converse inclusion holds. 

So we have shown that $f_*^{-1}(P) = \bar{f}_*^{-1}(P) \cap C$ is an intersection of two rational polyhedral cones, hence it is itself rational polyhedral. 

Now $f_*^{-1}(P)$ guarantees the cone conjecture for $(X,\Delta_X)$ in the sense of Theorem \ref{Thm equivalent cone}, i.e. $\bigcup_{g\in {\rm PsAut}(X,\Delta_X)}
g_* f_*^{-1}(P) \supset \mathcal{M}^e(X)$. 
In fact, we have inclusions 
$$ \bigcup_{g\in {\rm PsAut}(X,\Delta_X)}
g_* f_*^{-1}(P) \supset \bigcup_{g\in G}
g_* f_*^{-1}(P) = f_*^{-1}(\bigcup_{g\in G}
g_* P) = f_*^{-1}\mathcal{M}^e(Y) = \mathcal{M}^e(X)$$ 

Thus the cone conjecture holds for $(X,\Delta_X)$.

\end{proof}

\section{The cone conjecture for surfaces and a generalisation}

\subsection{The cone conjecture for klt surfaces over non-closed fields}

We prove the cone conjecture for klt Calabi-Yau surfaces over an arbitrary field of characteristic $0$. This should be known in some sense, but we include a proof for completeness. The ideas essentially come from Kawamata \cite{Kaw97} and Totaro \cite{Tot09}, with techniques in Oguiso-Sakurai \cite{OS01} and Bright-Logan-van Luijk \cite{BLvL19}. 

Our main result here is: 

\begin{theorem}\label{Cone non-closed}
Let $X$ be a klt Calabi-Yau surface over a field $k$ of characteristic $0$.

Then there exists a closed rational polyhedral cone $P \subset N^1(X)$ such that ${\rm Aut}(X)P =
\mathcal{A}^e(X)$. Or equivalently, the cone conjecture holds for $X$.
\end{theorem}

As the first step, we observe that the proof of \cite[Lemma 3.4]{Tot09} is still valid for non-algebraically closed fields: 

\begin{lemma}(\cite[Lemma 3.4]{Tot09})\label{Tot09 Lemma 3.4}
Let $f: X \to Y$ be a proper birational morphism of projective klt surfaces over a field $k$ of characteristic $0$. Let $\Delta_X$
be an $\mathbb{R}$-divisor on $X$ and $\Delta_Y$ its pushforward to $Y$. If ${\rm Aut}(X, \Delta_X)$ admits a rational
polyhedral fundamental domain on the $\mathcal{A}^e(X)$, then ${\rm Aut}(Y, \Delta_Y)$ admits a rational
polyhedral fundamental domain on the $\mathcal{A}^e(Y)$.

\end{lemma}

\begin{proof}
Consider the natural map $f^*: N^1(Y) \to N^1(X)$. This is clearly injective, so we may identify $N^1(Y)$ as a subspace of $N^1(X)$. Under this identification, we have $\mathcal{B}^e(Y) = \mathcal{B}^e(X) \cap N^1(Y)$ and $\bar{\mathcal{A}}(Y) = \bar{\mathcal{A}}(X) \cap N^1(Y)$. So $\mathcal{A}^e(Y) = \mathcal{A}^e(X) \cap N^1(Y)$. 

The subgroup $H$ of $G = {\rm Aut}(X, \Delta)$ that maps the subspace $N^1(Y)$ into itself is a subgroup of ${\rm Aut}(Y, \Delta_Y)$. If we prove the cone conjecture
for this subgroup of ${\rm Aut}(Y, \Delta_Y )$, the statement for the whole group ${\rm Aut}(Y, \Delta_Y )$
follows from Theorem \ref{Thm equivalent cone}. 

We know that there is a rational polyhedral cone $\Pi$ for $G$ acting on $\mathcal{A}^e(X)$,
i.e. 
$\mathcal{A}^e(X)= \bigcup_{g\in G} g_*\Pi$. It follows that $\mathcal{A}^e(Y)= \bigcup_{g\in G} g_*\Pi\cap N^1(Y)$. Here each set 
$g_*\Pi\cap N^1(Y)$ is a rational polyhedral cone contained in $\mathcal{A}^e(Y)$. 

It remains to show that the cones $g_*\Pi\cap N^1(Y)$ lies in finitely many orbits under the action of $H$. 

For an element $g$ of $G$, $g_*\Pi\cap N^1(Y)$ is a face of $g_*\Pi$ (possibly just $\{0\}$). So we can divide
the nonzero intersections $g_*\Pi\cap N^1(Y)$
into finitely many classes corresponding to the
faces $\Pi_i$ of $\Pi$ such that $g_*\Pi\cap N^1(Y) = g_*\Pi_i$. Fix a face $\Pi_i$ of $\Pi$ (there are only finitely many possibilities). If $\Pi_i = \{0\}$, then all the
cones $g_*\Pi_i$ are equal to $\{0\}$ and so they form a single $H$-orbit. So we can assume that
the face $\Pi_i$ of $\Pi$ is not $\{0\}$. 

Consider the contraction defined by an element in the interior of $\Pi_i$, say $g: X \to Z$. Since $\Pi_i$ contains no ample divisors, $Z$ is nonzero, so the fiber dimension of $g$ is at most $1$. In particular, there are only finitely
many numerical equivalence classes of curves in $X$ contracted by $g$. Therefore, the number of intermediate contractions $X \to Y \to Z$ is finite, say $X \to Y_j \to Z$, $1\leq j\leq r$. Now we fix one
value of $j$, and consider only those cones $g_*\Pi_i$ such that $g$ moves the contraction
$X \to Y_j$ to $f: X \to Y$. We claim that these cones form only a single orbit under the
group $H$. Indeed, for two cones $g_{1*}\Pi_i$ and $g_{2*}\Pi_i$ with the fixed data $(i,j)$, the
automorphism $g_2g^{-1}_1$
of $X$ moves the cone $g_{1*}\Pi_i$ to the cone $g_{2*}\Pi_i$, and it preserves
the contraction $f: X \to Y$, which means that it belongs to the subgroup $H$ of $G$.

Thus the cones $g_*\Pi\cap N^1(Y)$ lies in finitely many orbits under the action of $H$. So the result follows.

\end{proof}

\begin{proof}(of Theorem \ref{Cone non-closed})
When $X$ is smooth, this follows from \cite{Kaw97}. 

When $X$ is canonical, by considering the minimal resolution, it follows from Lemma \ref{Tot09 Lemma 3.4} and the smooth case. The minimal resolution is defined over $k$ by uniqueness of the minimal resolution. 

When $X$ is not canonical, let $Y$ be the global index one cover of $X$, then $Y_{\bar{k}}$ is the global index one cover of $X_{\bar{k}}$. Let $g$ be a generator of ${\rm Aut}(Y/X) = {\rm Aut}(Y_{\bar{k}}/X_{\bar{k}})$. 

When $Y_{\bar{k}}$ is an abelian surface, we have the following commutative diagram by passage to limit (we may assume $g$ also descends to $g_l:Y_l \to Y_l$): 
\[\begin{tikzcd}
	{Y_{\bar{k}}} & {Y_l} \\
	{X_{\bar{k}}} & {X_l} \\
	{\bar{k}} & l
	\arrow[from=3-1, to=3-2]
	\arrow[from=2-2, to=3-2]
	\arrow[from=2-1, to=3-1]
	\arrow["\ulcorner"{anchor=center, pos=0.125}, draw=none, from=2-1, to=3-2]
	\arrow[from=1-1, to=2-1]
	\arrow[from=1-1, to=1-2]
	\arrow[from=1-2, to=2-2]
	\arrow[from=2-1, to=2-2]
	\arrow["\ulcorner"{anchor=center, pos=0.125}, draw=none, from=1-1, to=2-2]
\end{tikzcd}\]
where $l$ is a finitely generated extension over $\mathbb{Q}$. Then we choose an embedding $l \to \mathbb{C}$, giving (denote the base change of $g_l$ by $g_\mathbb{C}:X_\mathbb{C} \to X_\mathbb{C}$) : 
\[\begin{tikzcd}
	{Y_{\bar{k}}} & {Y_l} & {Y_\mathbb{C}} \\
	{X_{\bar{k}}} & {X_l} & {X_\mathbb{C}} \\
	{\bar{k}} & l & {\mathbb{C}}
	\arrow[from=3-1, to=3-2]
	\arrow[from=2-2, to=3-2]
	\arrow[from=2-1, to=3-1]
	\arrow["\ulcorner"{anchor=center, pos=0.125}, draw=none, from=2-1, to=3-2]
	\arrow[from=3-3, to=3-2]
	\arrow[from=2-3, to=2-2]
	\arrow[from=2-3, to=3-3]
	\arrow["\ulcorner"{anchor=center, pos=0.125, rotate=-90}, draw=none, from=2-3, to=3-2]
	\arrow[from=1-1, to=2-1]
	\arrow[from=1-1, to=1-2]
	\arrow[from=1-2, to=2-2]
	\arrow[from=1-3, to=2-3]
	\arrow[from=1-3, to=1-2]
	\arrow[from=2-1, to=2-2]
	\arrow["\ulcorner"{anchor=center, pos=0.125}, draw=none, from=1-1, to=2-2]
	\arrow["\ulcorner"{anchor=center, pos=0.125, rotate=-90}, draw=none, from=1-3, to=2-2]
\end{tikzcd}\]
Then by universal property of algebraic closure, we have: 
\[\begin{tikzcd}
	{Y_{\bar{k}}} & {Y_{\bar{l}}} & {Y_\mathbb{C}} \\
	{X_{\bar{k}}} & {X_{\bar{l}}} & {X_\mathbb{C}} \\
	{\bar{k}} & {{\bar{l}}} & {\mathbb{C}}
	\arrow[from=3-1, to=3-2]
	\arrow[from=2-2, to=3-2]
	\arrow[from=2-1, to=3-1]
	\arrow["\ulcorner"{anchor=center, pos=0.125}, draw=none, from=2-1, to=3-2]
	\arrow[from=3-3, to=3-2]
	\arrow[from=2-3, to=2-2]
	\arrow[from=2-3, to=3-3]
	\arrow["\ulcorner"{anchor=center, pos=0.125, rotate=-90}, draw=none, from=2-3, to=3-2]
	\arrow[from=1-1, to=2-1]
	\arrow[from=1-1, to=1-2]
	\arrow[from=1-2, to=2-2]
	\arrow[from=1-3, to=2-3]
	\arrow[from=1-3, to=1-2]
	\arrow[from=2-1, to=2-2]
	\arrow["\ulcorner"{anchor=center, pos=0.125}, draw=none, from=1-1, to=2-2]
	\arrow["\ulcorner"{anchor=center, pos=0.125, rotate=-90}, draw=none, from=1-3, to=2-2]
\end{tikzcd}\]

By the arguments after \cite[Theorem 1.6]{Suz01}, the action of $g_\mathbb{C}$ on $N^1(Y_{\mathbb{C}})$ is trivial, hence by Lemma \ref{Lem change of closed fields}, the action of $g$ on $N^1(Y_{\bar{k}})$ is trivial. On the other hand, the natural map $N^1(Y)\to N^1(Y_{\bar{k}})$ is injective, so the action of $g$ on $N^1(Y)$ is trivial. So it suffices to show that the image of $H^\prime = \{h\in {\rm Aut}(Y) | hg = gh\}$ is of finite index in the image of ${\rm Aut}(Y)$ in ${\rm GL}(N^1(Y))$, which reduces to show that the image of $H = \{h\in {\rm Aut}(Y_{\bar{k}}) | hg = gh\}$ is of finite index in ${\rm Aut}(Y_{\bar{k}})/Aut_0(Y_{\bar{k}})$. 

 Take a closed point $p$ in $Y_{\bar{k}}$, then consider the subgroup $K$ of ${\rm Aut}(Y_{\bar{k}})$ fixing $p$. Then the natural map $K \to {\rm Aut}(Y_{\bar{k}})/Aut_0(Y_{\bar{k}})$ is surjective, so it suffices to show $K\cap H$ is of finite index in $K$. Indeed, we show that $J = \{h^{-1}gh|h\in K\}$ is finite. If this were true, then we have an exact sequence $0 \to K\cap H \to K \to Sym(J)$, so $K\cap H$ is of finite index in $K$. 

Now it remains to show that $J = \{h^{-1}gh|h\in K\}$ is finite. Indeed, $J$ acts on $N^1(Y_{\bar{k}})$ trivially, so it lies in only finitely many components in ${\rm Aut}(Y_{\bar{k}})$. Moreover, the intersection with each component contains at most one element, so $J$ is finite. 

Otherwise (i.e. when $X$ is not canonical and the global index $1$ cover is not an abelian surface), as in \cite[Theorem 3.2]{Tot09}, $Y_{\bar{k}}$ is a canonical surface with trivial canonical bundle, so is $Y$. Let $Z$ be the minimal resolution of $Y$, $Z$ is defined on $k$ by uniqueness of the minimal resolution. Then $Z$ is a K3 surface, and we have a commutative diagram by uniqueness of the minimal resolution: 
\[\begin{tikzcd}
	{Z} & {W} \\
	{Y} & {X}
	\arrow["{\mathbb{Z}/I}"', from=1-1, to=1-2]
	\arrow[from=1-2, to=2-2]
	\arrow["{\mathbb{Z}/I}"', from=2-1, to=2-2]
	\arrow[from=1-1, to=2-1]
\end{tikzcd}\]

By Lemma \ref{Tot09 Lemma 3.4} again, it suffices to prove the cone conjecture for $W$, i.e. we may assume $Y$ is itself smooth. 
We prove this in the next lemma, see Lemma \ref{quotient of K3}. 

\end{proof}

We note that the following lemma is proved when $k = \mathbb{C}$ by Oguiso-Sakurai \cite{OS01} or when $H$ is trivial by Bright-Logan-van Luijk \cite{BLvL19}. The same technique works for the general case. We essentially follow the proof of \cite{OS01} and \cite{BLvL19}, making necessary adjustments:

\begin{lemma}\label{quotient of K3}

Let $X$ be a K3 surface over a field $k$ of characteristic 0, $H$ be a finite group acting faithfully on $X$. 
Then the action of ${\rm Aut}(X/H)$ on $\mathcal{A}^e(X/H)$ admits a rational polyhedral fundamental domain. 

\end{lemma}

\begin{proof}
Let $G$ be the absolute Galois group of $\eta$. Then the actions of $G$ and $H$ on $\bar{X} = X\times_\eta\bar{\eta}$ commute. 

We have natural identifications $N^1(X/H) \simeq N^1(\bar{X})^{G\times H}$, $\mathcal{A}^e(X/H) \simeq \mathcal{A}^e(\bar{X})^{G\times H}$. 

The natural action of $G\times H$ on $NS(\bar{X})$ gives a representation $G\times H \to O(NS(\bar{X}))$. We let $G\times H$ act on $O(NS(\bar{X}))$ by conjugation. Note that a conjugation of a $(-2)$-class is still a $(-2)$-class, so this action restricts to an action of $G\times H$ on $W(NS(\bar{X}))$. 
We define the equivariant reflection group $R_X$ to be $W(NS(\bar{X}))^{G\times H}$.

Let $C_X$ be the rational hull of the positive cone in $N^1(\bar{X})$, and $C_{X/H} = C_X \cap N^1(X/H)$. 

Consider the action of ${\rm Aut}(\bar{X})$ on $R_X$ by conjugation, defining the semi-direct product ${\rm Aut}(\bar{X})^{G\times H} \ltimes R_X$. 

We separate the proof into 4 steps. In Step 1, we give a concrete description of $R_X$, which shows that $R_X$ is generated by reflections of $G\times H$-orbits of $(-2)$-curves. In Step 2, we show that the effective nef cone is a fundamental domain for the faithful action of $R_X$ on $C_X$. In Step 3, we prove that the natural map
${\rm Aut}(\bar{X})^{G\times H} \ltimes R_X \to O(NS (\bar{X})^{G\times H})$
has finite kernel and image of finite index. In Step 4, we deduce our main result. 

\hfill

\noindent \textbf{Step 1. } Equivariant reflection group (\cite[Definition 3.3, Proposition 3.6]{BLvL19})

In this step, we prove that $R_X$ is generated by reflections of $G\times H$-orbits of $(-2)$-curves. It turns out that we can give a concrete description of $R_X$ as a Coxter group. 

Following \cite[Proposition 3.6]{BLvL19}, we describe $R_X = W(NS(\bar{X}))^{G\times H}$ as a Coxeter system: 

By \cite[Theorem 3.5]{BLvL19}, let $F$ be the set of orbits $I$ of $(-2)$-curves (under the action of $G\times H$) for which the subgroup $W_I$ generated by reflections in $I$ is finite, and let $r_I$ be the longest element in $(W_I,I)$, then $(R_X, \{r_I:I\in F\})$ is a Coxeter system. 

Moreover, we know more about all possible $I$ and $w_{I,0}$: 

(1) There are two possibilities of $I$: 

\text{ }\text{ }\text{ }\text{ }(i) $I$ consists of disjoint $(-2)$-curves;

\text{ }\text{ }\text{ }\text{ }(ii) $I$ consists of disjoint pairs of $(-2)$-curves, each pair having intersection number $1$. 

(2) For each $I \in F$, let $C_I \in {\rm Pic} (\bar{X})^{G\times H}$ be the sum of the classes in $I$. Then
$r_I$ acts on $NS (\bar{X})^{G\times H}$ as reflection in the class $C_I$. 

Now we prove (1) and (2): 

Proof of (1): Let $I$ be a $G\times H$-orbit of $(-2)$-curves, such that the subgroup $W_I$
is finite. Firstly, observe that no two
$(-2)$-curves in $I$ have intersection number at least $2$, otherwise the corresponding
reflections of these two curves would generate an infinite dihedral subgroup. Since $W_I$ is finite,
its Coxeter–Dynkin diagram would be a finite union of trees by \cite[Exercise 1.4]{BB05}. In particular,
it contains a vertex of degree $\leq 1$. However, the group $G\times H$ acts transitively
on the vertices of this diagram, so we conclude that either every vertex has degree $0$, or every
vertex has degree $1$. Then these two cases correspond to (i) and (ii). 

Proof of (2): In Case (i), we have
$I = \{E_1, \dots, E_r\}$. The reflections in the $E_i$ all commute, so $W_I$ is isomorphic to the
Coxeter group $A^r_1 = (\mathbb{Z}/2\mathbb{Z})^r$. The longest element is $$r_I = r_{E_1} \circ \dots \circ r_{E_r}$$
For $D \in NS(\bar{X})^{G\times H}$,
the intersection numbers $D \cdot E_i$ are all equal, and one calculates that $r_I$ coincides with the reflection in the class $C_I = E_1 + \cdots E_r$ for $D$.

In Case (ii), we write $I = {E_1, E^\prime_1
, \dots , E_r, E^\prime_r}$, where $E_i\cdot E^\prime_i = 1$ and all the
other intersections are zero. The two reflections $r_{E_i}$ and $r_{E^\prime_i}$
generate a
subgroup isomorphic to the Coxeter group $A_2 = S_3$, in which the longest element is
$$r_I = r_{E_i} \circ r_{E^\prime_i} 
\circ r_{E_i} = r_{E^\prime_i} \circ r_{E_i} 
\circ r_{E^\prime_i} = r_{E_i+E^\prime_i}$$
Thus we have $ W_I
= A^r_2 = (S_3)^r$ and the longest element is $$r_I = r_{E_1+E^\prime_1} \circ \dots \circ r_{E_r+E^\prime_r}$$
For $D \in NS(\bar{X})^{G\times H}$,
the intersection numbers $D \cdot E_i$ (resp. $D \cdot E^\prime_i$) are all equal, and one calculates that $r_I$ coincides with the reflection in the class $C_I = E_1 + \cdots E_r$ for $D$. 

In particular, (2) shows that $R_X$ is generated by reflections of $G\times H$-orbits of $(-2)$-curves. 

\hfill

\noindent \textbf{Step 2. } Identification of the nef cone (\cite[Proposition 3.7]{BLvL19})

In this step, we prove that the cone $\mathcal{A}^e(X/H) \cap C_{X/H}$ is a
fundamental domain for the action of $R_X$ on $C_{X/H}$, and this action
is faithful. 

More precisely, we need to prove two things: (1) every class in $C_{X/H}$ is $R_X$-equivalent to
an element of $\mathcal{A}^e(X/H) \cap C_{X/H}$; (2) the translates of $\mathcal{A}^e(X/H) \cap C_{X/H}$ by
two distinct elements of $R_X$ meet only along their boundaries. 

Proof of (1): Let $D$ be any class in $C_{X/H}$. Suppose first that $D$ has trivial
stabilizer for the action of the Weyl group $W(NS(\bar{X}))$. Then, by the corresponding result for K3 surfaces over algebraically closed fields, see for example \cite[Corollary 8.2.11]{Huy16}, there exists a unique $g \in W(NS(\bar{X}))$ such that $gD$ lies in the interior of $\mathcal{A}^e(\bar{X} )\cap  C_X$. We claim
that $g$ lies in $R_X$. Indeed, for any $\sigma \in G\times H$, since the $G\times H$ action preserves the properties of being nef and positive, we have
$$(\sigma \cdot g)D = \sigma(g(\sigma^{-1}D)) = \sigma(g(D)) \in \mathcal{A}^e(\bar{X} )\cap  C_X$$
By
uniqueness of $g$, we conclude that $g = \sigma \cdot g$, i.e. $g$ lies in the equivariant Weyl group $R_X$. Then $g$ preserves the subspace $N^1(X/H)$, so $gD$ lies in $\mathcal{A}^e(X/H)\cap  C_{X/H}$. 

Now suppose that $D$ has non-trivial stabilizer. Then $D$ lies on at least one
of the walls defined by the action of $W(NS(\bar{X}))$ on $N^1(\bar{X})$.
The chamber structure of this group action is locally polyhedral inside $C_X$, so a small enough neighbourhood of $D$ meets only finitely
many chambers. Also note that $NS(X/H)$ is not contained in any of the walls, since
$X/H$ admits an ample divisor. We may choose a sequence $(D_i)_{i=1}^{\infty}$ of
elements of $C_{X/H}$, tending to $D$ and all have trivial stabilizer for the action of $W(NS(\bar{X}))$. By passing to a subsequence, we may assume all $D_i$ lie in the interior of the same chamber
of $C_X$. As in the previous paragraph, there is a unique $g_i \in R_X$ satisfying $g_iD_i \in \mathcal{A}^e(X/H)\cap  C_{X/H}$ for all $i$. All $g_i$ coincide since all $D_i$ lie in the same chamber. By continuity, $gD$ also lies in $\mathcal{A}^e(X/H)\cap  C_{X/H}$.

Proof of (2): Suppose that $x \in C_{X/H}$
lies in the intersection $\mathcal{A}^e(X/H) \cap g \mathcal{A}^e(X/H)$, for some non-trivial $g \in R_X$. By the corresponding result for K3 surfaces over algebraically closed fields, see for example \cite[Corollary 8.2.11]{Huy16}
, we see that $x$ lies in the boundary of $\mathcal{A}^e(\bar{X})$. Then 
\cite[Lemma 3.8]{BLvL19} shows that in fact $x$ lies in the boundary of $\mathcal{A}^e(X/H)$. 

\hfill

\noindent \textbf{Step 3. } Equivariant Torelli (\cite[Proposition 3.10]{BLvL19})

In this step, we prove that the natural map
${\rm Aut}(\bar{X})^{G\times H} \ltimes R_X \to O(NS (\bar{X})^{G\times H})$
has finite kernel and image of finite index. 

This follows from abstract group theory and the usual version of Torelli theorem. Indeed, by the usual version of Torelli theorem, the natural map
$${\rm Aut}(\bar{X}) \to O(NS(\bar{X}))/W(NS (\bar{X}))$$ has finite kernel and image of finite index. Note that the action of $G\times H$ on all of these groups factors
through a finite quotient. So by \cite[Lemma 3.12]{BLvL19}, the induced homomorphism
$${\rm Aut}(\bar{X})^{G\times H} \to (O(NS(\bar{X}))/W(NS (\bar{X})))^{G\times H}$$
also has finite kernel and image of finite index. 

On the other hand, there is an exact sequence
$$1 \to R_X \to O(NS (\bar{X}))^{G\times H}
\to  (O(NS(\bar{X}))/W(NS (\bar{X})))^{G\times H}$$
and the homomorphism from ${\rm Aut}(\bar{X})$ to $(O(NS(\bar{X}))/W(NS (\bar{X})))^{G\times H}$ factors through $O(NS (\bar{X}))^{G\times H}$, and hence through the
injective map
$$O(NS (\bar{X}))^{G\times H}/R_X
\to  (O(NS(\bar{X}))/W(NS (\bar{X})))^{G\times H}$$
Therefore, by \cite[Lemma 3.13]{BLvL19}, the map
$${\rm Aut}(\bar{X})^{G\times H} \to O(NS (\bar{X}))^{G\times H}/R_X$$ has finite kernel and image of finite index. Since the image of ${\rm Aut}(\bar{X})^{G\times H}$ in $O(NS (\bar{X}))^{G\times H}$
meets $R_X$ only in the identity element, this shows that the natural map
$${\rm Aut}(\bar{X})^{G\times H} \ltimes R_X \to O(NS (\bar{X}))^{G\times H}$$
also has finite kernel and image of finite index.

We apply \cite[Proposition 2.2]{BLvL19} with $\Lambda = NS (\bar{X})^{G\times H}$ and $H^\prime$ being the image of $G\times H$
in $O(NS( \bar{X}))$.  The centralizer $Z_{O(\Lambda)}(H^\prime)$ is $O(NS (\bar{X}))^{G\times H}$, so part (2) of
\cite[Proposition 2.2]{BLvL19} shows that $O(NS (\bar{X}))^{G\times H}$ is of finite index in $O(NS( \bar{X})
,NS (\bar{X})^{G\times H})$.
Parts (1) and (2) of \cite[Proposition 2.2]{BLvL19} combined show that the natural map
$$O(NS( \bar{X})
,NS (\bar{X})^{G\times H}) \to  O(NS (\bar{X})^{G\times H})$$
has finite kernel and image of finite index; by \cite[Lemma 3.13]{BLvL19} so does the map
$$O(NS (\bar{X}))^{G\times H} \to O(NS (\bar{X})^{G\times H})$$

Combining these results, the natural map
$${\rm Aut}(\bar{X})^{G\times H} \ltimes R_X \to O(NS (\bar{X})^{G\times H})$$
also has finite kernel and image of finite index. 

\hfill

\noindent \textbf{Step 4. } Existence of a fundamental domain (\cite[Corollary 3.15]{BLvL19})

In this step, we finally deduce our main result: the action of ${\rm Aut}(X/H)$ on $\mathcal{A}^e(X/H)$ admits a rational polyhedral
fundamental domain. 

In the argument of \cite[Corollary 3.15]{BLvL19}, we take $\Lambda = NS(\bar{X})^{G\times H}$, and $\Gamma$ be the image of $${\rm Aut}(\bar{X})^{G\times H} \ltimes R_X \to O(NS (\bar{X})^{G\times H})$$ and $C = C_{X/H}$. Let $C_+$ be the rational hull of $C$. Pick any $y \in C \cap \Lambda$, then the set
$$\Pi = \{x \in C_+ | (\gamma x \cdot y) \geq (x \cdot y) \text{ for all } \gamma \in \Gamma\}$$
is rational polyhedral. In particular, $\Pi$ gives a rational polyhedral fundamental domain for the action of ${\rm Aut}(\bar{X})^{G\times H} \ltimes R_X$ on $C_+$. 

If we take $y$ to be an ample class, we claim that $\Pi$ is contained in the nef cone, following  \cite[Remark 3.9]{BLvL19}. Indeed, the nefness condition is given by positive intersection with all elements of the form $C_I$ as defined in Step $1$, and taking $\gamma = r_I$ shows that elements in $\Pi$ have positive intersection with $C_I$. 

Then we conclude that $\Pi$ is also a rational polyhedral fundamental domain for the action of ${\rm Aut}(\bar{X})^{G\times H}$ on $\mathcal{A}^e(X/H) = \mathcal{A}^e(\bar{X})^{G\times H}$. 
Indeed, if $x$ is a class in $\mathcal{A}^e
(X)$, since $\Pi$ is a fundamental
domain for the action of ${\rm Aut}(\bar{X})^{G\times H} \ltimes R_X$ on $C_+$, we can find $\phi \in {\rm Aut}(\bar{X})^{G\times H}$ and $r \in R_X$
such that $r\phi(x)$ lies in $\Pi$. But now $\phi(x)$ and $r\phi(x)$ both lie in $\mathcal{A}^e(X)$, so they are
equal and lie in $\Pi$. This shows that $\Pi$ is a fundamental domain for the action of
${\rm Aut}(\bar{X})^{G\times H}$ on $\mathcal{A}^e(X/H) = \mathcal{A}^e(\bar{X})^{G\times H}$.

But we have a commutative diagram: 
\[\begin{tikzcd}
	{{\rm Aut}(\bar{X})^{G\times H}} & {O(NS(\bar{X})^{G\times H})} & {O(N^1(\bar{X})^{G\times H})} \\
	{{\rm Aut}(X/H)} & {O(NS(X/H))} & {O(N^1(X/H))}
	\arrow[from=1-1, to=2-1]
	\arrow[from=2-1, to=2-2]
	\arrow[from=2-2, to=2-3]
	\arrow[from=1-1, to=1-2]
	\arrow[from=1-2, to=1-3]
	\arrow["\simeq", from=1-3, to=2-3]
\end{tikzcd}\]
Since $(\mathcal{A}^e(\bar{X})^{G\times H}, {{\rm Aut}(\bar{X})^{G\times H}})$ is of polyhedral type, $(\mathcal{A}^e(X/H), {{\rm Aut}(X/H)})$ is also of polyhedral type. Now the result follows from  \cite[Application 4.14]{Loo14}.

\end{proof}

\subsection{A generalised cone conjecture}

\label{Subsection generalised cone}

In this part, we prove Theorem \ref{Thm generalised cone}, which might be considered as a generalised version of the cone conjecture for klt Calabi-Yau varieties over finitely generated fields over $\mathbb{Q}$. We refer to Definition \ref{Def cones in Pic} for notations on cones in ${\rm Pic}(X)\otimes_\mathbb{Z} \mathbb{R}$. 

\subsubsection{Dimension 1 case}

\begin{lemma}\label{Lem Cone of elliptic}
Let $X$ be a genus $1$ curve over a finitely generated field $k$ over $\mathbb{Q}$. 

Then there exists a closed rational polyhedral cone $P\subset {\rm Pic}(X)\otimes_\mathbb{Z} \mathbb{R}$ such that ${\rm Aut}(X) P = {\rm Nef}^e(X)$. 
\end{lemma}

\begin{remark}
In this case, we have ${\rm Nef}^e(X) = \{0\}\cup \{D\in {\rm Pic}(X)\otimes_\mathbb{Z} \mathbb{R} |\deg D>0\}$. 
\end{remark}

\begin{proof}

Recall that ${\rm Pic}^0(X)$ acts on $X$ by translation, so ${\rm Pic}^0(X)$ acts on ${\rm Pic}(X)$ by pullback of line bundles. By \cite{Tot09}, we know that the translation by an element $x$ of ${\rm Pic}^0
(X)$ acts on
${\rm Pic}(X)$ by $\phi_x(y) = y + \deg(y)x$.

Pick $z_1,\dots,z_n\in {\rm Pic}^0(X)$ such that their images in ${\rm Pic}^0(X)\otimes_\mathbb{Z} \mathbb{R}$ form an $\mathbb{R}$-basis. Let $z$ be any degree $1$ element in ${\rm Pic}(X)\otimes_\mathbb{Z} \mathbb{Q}$, and $P$ be the cone in ${\rm Pic}(X)\otimes_\mathbb{Z} \mathbb{R}$ generated by the finite set $\{z+\sum_{i\in I}z_i|I\subset \{1,2,\dots,n\}\}$.  

Let $a\in {\rm Nef}^e(X)$ be any class. By scaling, we may assume $\deg a = 1$, say $a = z + \sum_{i=1}^n a_iz_i$. Taking translations by a integral combination of $z_1,\dots,z_n$, we may assume $0\leq a_i <1$, so $a$ lies in $P$. 

\end{proof}

\subsubsection{Abelian and bielliptic surfaces}

For abelian or bielliptic surfaces, we proceed as Kawamata's proof in \cite[Lemma 4.3]{Kaw97}. We note that Li's result \cite{Li23} essentially gives a systematic way to pass from the cone conjecture to this generalised version. 

\begin{theorem}\label{Thm generalised abelian}
Let $X$ be a surface defined on a finitely generated field extension of $\mathbb{Q}$, say $k$. Assume $X_{\bar{k}}$ is an abelian surface. 

Then there exists a closed rational polyhedral cone $P\subset {\rm Pic}(X)\otimes_\mathbb{Z} \mathbb{R}$ such that ${\rm Aut}(X)P = {\rm Nef}^e(X)$. 
\end{theorem}

\begin{proof}
We have a natural projection $p:{\rm Pic}(X)\otimes_\mathbb{Z} \mathbb{R} \to N^1(X)$. 

Since the cone conjecture holds by \cite[Theorem 2.1, Remark 2.2]{Kaw97}, let $\Pi \subset \mathcal{A}^e(X)$ be a closed rational polyhedral fundamental domain for the action of ${\rm Aut}(X)$ on $\mathcal{A}^e(X)$. 

Let $\{v_i\}_{i\in I}$ be a finite set of rational vectors that generate $\Pi$. 
Take a set of rational liftings $\{\tilde{v}_i\}_{i\in I}$ such that $\tilde{v}_i\in {\rm Nef}^e(X)$. 
Let $\tilde{\Pi}$ be the cone generated by $\{\tilde{v}_i\}_{i\in I}$. 

By assumption, $X$ is a torsor of the abelian surface $A = \mathcal{A}lb^0(X)$. 
Let $A(k)$ be the group of rational points in $A$, $x_1,\dots, x_n\in A(k)$ be a basis of the vector space $A(k)\otimes_\mathbb{Z} \mathbb{R}$. 
By the definition of torsors, we have an action $A \times X \to X$, hence we have a homomorphism  $A(k) \to {\rm Aut}(X)$.

Consider the cone $P$ generated by finitely many rational polyhedral cones $\left(\sum_{I\subset [n]}x_i\right) _*\tilde{\Pi}$, where $I$ run through all subsets of $[n] = \{1,2,\dots,n\}$. 
We claim that ${\rm Aut}(X)P = {\rm Nef}^e(X)$. Since ${\rm Nef}^e(X)$ is generated by integral divisors, it suffices to show any integral class in $C$ is contained in ${\rm Aut}(X)P$. 

Indeed, let $D \in {\rm Pic}(X)\otimes_\mathbb{Z} \mathbb{R}$ be any integral class in ${\rm Nef}^e(X)$. Let $[D]$ be the image of $D$ in $N^1(X)$. By construction, $g_*[D] \in \Pi$ for some $g\in {\rm Aut}(X)$. Now there exists some rational class $\tilde{v}\in \tilde{\Pi}$ such that $\tilde{v}-g_*D\in \ker p$. Replacing $D$ and $\tilde{v}$ by a suitable multiple, we may assume $\tilde{v}$ is an integral divisor. To show $D\in {\rm Aut}(X)P$, it suffices to show that $\tilde{v}-g_*D$ lies in the subspace spanned by $\{x_{i*}\tilde{v} - \tilde{v}\}_{ 1\leq i\leq n}$. 

If so, we have an equality in ${\rm Pic}(X)\otimes_\mathbb{Z} \mathbb{R}$: 
$$\tilde{v} - g_*D = \sum_i a_i(x_{i*}\tilde{v} - \tilde{v}) = \sum_i (n_i-c_i)(x_{i*}\tilde{v} - \tilde{v})$$ where $n_i\in \mathbb{Z}$, $0\leq c_i\leq 1$. Then this equility can be rewritten as $$g_*D +\sum_i n_i(x_{i*}\tilde{v} - \tilde{v})= \tilde{v} +\sum_i c_i(x_{i*}\tilde{v} - \tilde{v})$$
Since the translation action depends only on the numerical equivalence class, we have 
$$\left(\sum_i n_ix_{i}\right)_*g_*D = \tilde{v} +\sum_i c_i(x_{i*}\tilde{v} - \tilde{v})$$ 
Now it remains to show $\tilde{v} +\sum_i c_i(x_{i*}\tilde{v} - \tilde{v})$ lies in $P$. To see this, without loss of generality, we may assume $c_i$ is non-decreasing, otherwise we may change the order. Denote $c_0= 0$, then we have 
\begin{align*}
    \tilde{v} +\sum_i c_i(x_{i*}\tilde{v} - \tilde{v}) =& \tilde{v} +\sum_i \left(c_i-c_{i-1}\right)\left(\left(\sum_{1\leq j\leq i}x_{j}\right)_*\tilde{v} - \tilde{v}\right) 
    \\=& (1-c_n)\tilde{v} + \sum_i (c_i-c_{i-1})\left(\sum_{1\leq j\leq i}x_{j}\right)_*\tilde{v} \in P
\end{align*}

Now it remains to show that $\tilde{v}-g_*D$ lies in the subspace spanned by $\{x_{i*}\tilde{v} - \tilde{v}\}_{ 1\leq i\leq n}$.

When $\tilde{v}$ is ample, as in \cite[Lemma 4.3]{Kaw97}, the map $\phi_{\tilde{v}}: A \to \mathcal{P}ic^0(X)$ given by translation of divisors $\phi_{\tilde{v}}(x) = T_{x}(\tilde{v}) - \tilde{v}$ is finite \'etale of degree $d$. Now $\tilde{v}-D$ is a rational point in $\mathcal{P}ic^0(X)$. Consider the sum of the inverse images $y=\sum_{\phi_{\bar{v}}(x)=\tilde{v}-D} x$. It is an $\eta$-rational point since it is invariant under the Galois action. So $$\phi_{\tilde{v}}(y) = y_*(\tilde{v}) - \tilde{v} =d(\tilde{v}-D)$$ By the definition of $x_1,\dots,x_n$, we have an equality $my=\sum_{i=1}^n a_ix_i$ for some $m\in \mathbb{N}_{>0}, a_i\in \mathbb{Z}$. 

Combining the above equalities, and applying $\phi_{\tilde{v}}$ on both sides of $my=\sum_{i=1}^n a_ix_i$, we have $md(\tilde{v}-D) = \sum_{i=1}^n a_i(x_{i*}(\tilde{v})-\tilde{v})$, hence $\tilde{v}-D$ lies in the subspace spanned by $\{x_{i*}\tilde{v} - \tilde{v}\}_{ 1\leq i\leq n}$. 

In general, $\tilde{v}$ is effective, hence semiample, so it defines a contraction to another variety $Y$ over $k$, say $f: X\to Y$. Possibly replacing $\tilde{v}$ by a multiple, we have $\tilde{v} = f^*(H)$ for some ample divisor $H$ on $Y$. By the same reason, $D$ defines a contraction $f: X\to Y$, and $D = f^*(H^\prime)$ for some ample divisor $H^\prime$ on $Y$. Here we may choose the same $f$ since $\tilde{v}$ and $D$ are numerically equivalent. We note that in this case $Y_{\bar{k}}$ is an abelian variety since all contractions of abelian varieties are still abelian varieties. (It is in fact an elliptic curve when $\tilde{v}$ is nonzero and not ample.)

On the other hand, consider the map $\phi_{\tilde{v}}: A \to  \mathcal{P}ic^0(X)$ given by translation of divisors $\phi_{\tilde{v}}(x) = T_{x}(\tilde{v}) - \tilde{v}$. Now this is not necessarily finite \'etale, but we have a factorisation 
\[\begin{tikzcd}
	{A} & {\mathcal{P}ic^0(X)} \\
	{\mathcal{A}lb^0(Y)} & {\mathcal{P}ic^0(Y)}
	\arrow["{\mathcal{A}lb^0(f)}"', from=1-1, to=2-1]
	\arrow["{\phi_H}", from=2-1, to=2-2]
	\arrow["{\phi_{\tilde{v}}}", from=1-1, to=1-2]
	\arrow["{f^*}"', from=2-2, to=1-2]
\end{tikzcd}\]

Now $\phi_H\circ \mathcal{A}lb^0(f)$ is surjective, since $\mathcal{A}lb^0(f)$ is a contraction and $\phi_H$ is finite \'etale. Let $y\in (\phi_H\circ \mathcal{A}lb^0(f))^{-1}(H-H^\prime)$ be any preimage of $H-H^\prime$. Let $z$ be the sum of all Galois conjugates of $y$, then $z$ is an $\eta$-rational point and $\phi_H\circ \mathcal{A}lb^0(f)(z) = d (H-H^\prime)$, where $d$ is the number of the Galois conjugates of $y$. Now by the commutative diagram, $$\phi_{\tilde{v}}(z) = z_*(\tilde{v})-\tilde{v} = k(\tilde{v}-D)$$

Then we may proceed as in the ample case. By the definition of $x_1,\dots,x_n$, we have an equality $mz=\sum_{i=1}^n a_ix_i$ for some $m\in \mathbb{N}_{>0}, a_i\in \mathbb{Z}$. 

Combining the above equalities, we have $md(\tilde{v}-D) = \sum_{i=1}^n a_i(x_{i*}(\tilde{v})-\tilde{v})$, hence $\tilde{v}-D$ lies in the subspace spanned by $\{x_{i*}\tilde{v} - \tilde{v}\}_{ 1\leq i\leq n}$.

\end{proof}

Then we treat bielliptic surfaces by reducing to abelian surfaces: 

\begin{theorem}\label{Thm cone of bielliptic}
Let $X$ be a surface defined on a finitely generated field extension of $\mathbb{Q}$, say $k$. Assume $X_{\bar{k}}$ is a bielliptic surface. 

Then there exists a closed rational polyhedral cone $P\subset {\rm Pic}(X)\otimes_\mathbb{Z} \mathbb{R}$ such that ${\rm Aut}(X)P = {\rm Nef}^e(X)$. 
\end{theorem}

\begin{proof}
We have a natural projection $p:{\rm Pic}(X)\otimes_\mathbb{Z} \mathbb{R} \to N^1(X)$. 

Since the cone conjecture holds as before, let $\Pi \subset \mathcal{A}^e(X)$ be a closed rational polyhedral fundamental domain for the action of ${\rm Aut}(X)$ on $\mathcal{A}^e(X)$. 

Let $\{v_i\}_{i\in I}$ be a finite set of rational vectors that generates $\Pi$. Take a set of rational liftings $\{\tilde{v}_i\}_{i\in I}$ such that $\tilde{v}_i\in {\rm Nef}^e(X)$. 
Let $\tilde{\Pi}$ be the cone generated by $\{\tilde{v}_i\}_{i\in I}$. 

By assumption, a finite \'etale cover of $X$, say $Y$, with projection $\pi: Y \to X$, is a torsor of an abelian surface $A = \mathcal{A}lb(X)$. 
Let $A(k)$ be the group of rational points in $A$, $x_1,\dots, x_n\in A(k)$ be a basis of the vector space $A(k)\otimes_\mathbb{Z} \mathbb{R}$. By the definition of torsors, we have an action $A \times Y \to Y$, hence we have a homomorphism  $\iota: A(k) \to {\rm Aut}(Y)$. Note that for $g\in {\rm Aut}(Y/X)$, the action $g^{-1}\iota(x_i)g$ is still given by a translation, say by $gx_i$. 

Formally, we may expand $gx_i = \sum_j c_jx_j$ in $A(k)\otimes_\mathbb{Z} \mathbb{R}$ with $c_j\in \mathbb{Q}$. Let $V = \sum_{i=1}^n \mathbb{Q}x_i$ be the abstract $\mathbb{Q}$-vector space with basis $x_1,\dots,x_n$. Then this defines an action of ${\rm Aut}(Y/X)$ on $V$. Consider the fixed subspace $V^{{\rm Aut}(Y/X)}$, which is has a basis $y_1,\dots, y_k$. Replacing $y_i$ by a suitable multiple, we may assume they are indeed invariant under the conjugate action of ${\rm Aut}(Y/X)$ as elements in $A(k)$. So translations by $y_i$ descend to automorphisms of $X$, this defines a map $\iota^\prime: G = \oplus_{i=1}^k\mathbb{Z}y_k \to {\rm Aut}(X)$. 

Consider the cone $P$ generated by the finitely many rational polyhedral cones $\left(\sum_{I\subset [k]}y_i\right)_*\tilde{\Pi}$, where $I$ run through all subsets of $[k] = \{1,2,\dots,k\}$. 
We claim that ${\rm Aut}(X)P = {\rm Nef}^e(X)$. 

Let $D \in {\rm Pic}(X)\otimes_\mathbb{Z} \mathbb{R}$ be any integral class in ${\rm Nef}^e(X)$. (It suffices to check for integral classes as ${\rm Nef}^e(X)$ is generated by integral divisors.) Let $[D]$ be the image of $D$ in $N^1(X)$. By construction, $g_*[D] \in \Pi$ for some $g\in {\rm Aut}(X)$. Now there exists some rational class $\tilde{v}\in \tilde{\Pi}$ such that $\tilde{v}-D\in \ker p$. Replacing $D$ and $\tilde{v}$ by a suitable multiple, we may assume that $\tilde{v}$ is an integral divisor. To show $D\in {\rm Aut}(X)P$, it suffices to show $\tilde{v}-D$ lies in the subspace spanned by $\{y_{i*}\tilde{v} - \tilde{v}\}_{ 1\leq i\leq k}$. 

But by Theorem \ref{Thm generalised abelian}, we have shown that $\pi^*(\tilde{v}) - \pi^*(D)$ lies in the subspace of ${\rm Pic}(Y)$ spanned by $\{x_{i*}\pi^*(\tilde{v}) - \pi^*(\tilde{v})\}_{ 1\leq i\leq n}$, say $$\pi^*(\tilde{v}) - \pi^*(D) = \sum_{ 1\leq i\leq n} a_i(x_{i*}\pi^*(\tilde{v}) - \pi^*(\tilde{v}))$$
Since $\pi^*(\tilde{v}) - \pi^*(D)$ is invariant under the action of ${\rm Aut}(Y/X)$, we have \begin{align*}
    \sum_{ 1\leq i\leq n} a_i(x_{i*}\pi^*(\tilde{v}) - \pi^*(\tilde{v})) = & \frac{1}{|{\rm Aut}(Y/X)|}\sum_{g\in {\rm Aut}(Y/X)} \sum_{ 1\leq i\leq n} a_i((gx_{i})_*\pi^*(\tilde{v}) - \pi^*(\tilde{v})) \\=& \sum_j b_j (y_{j*}\pi^*(\tilde{v})-\pi^*(\tilde{v}))
\end{align*}
This gives $$\pi^*(\tilde{v}-D) = \pi^*\left(\sum_j b_j (y_{j*}\tilde{v}-\tilde{v})\right)$$
But $\pi^*: {\rm Pic}(X)\otimes_\mathbb{Z} \mathbb{R} \to {\rm Pic}(Y)\otimes_\mathbb{Z} \mathbb{R}$ is injective, so $\tilde{v}-D = \sum_j b_j (y_{j*}\tilde{v}-\tilde{v})$, which means $\tilde{v}-D$ lies in the subspace spanned by $\{y_{i*}\tilde{v} - \tilde{v}\}_{ 1\leq i\leq k}$.

\end{proof}

\subsubsection{Other cases in dimension 2}

The following lemma essentially shows that other cases of Theorem \ref{Thm generalised cone} follow directly from the cone conjecture. 

\begin{lemma}\label{Lem irregularity}
Let $X$ be a klt Calabi-Yau surface over $\mathbb{C}$ that is not smooth. Then $H^1(X,\mathcal{O}_X) = 0$. In particular, $N^1(X) = {\rm Pic}(X)\otimes_\mathbb{Z} \mathbb{R}$. 
\end{lemma}

\begin{proof}
Let $\tilde{X}$ be the minimal resolution of $X$. 

If $X$ is canonical, then  $\tilde{X}$ is a smooth Calabi-Yau surface. But the exceptional curves are rational, while abelian surfaces do not contain rational curves. Also, if $\tilde{X}$ were a bielliptic surface, let $\pi: Y\to \tilde{X}$ be a finite \'etale cover where $Y$ is an abelian surface. Let $E$ be an exceptional divisor of $\tilde{X}$ over $X$, then $\pi^*E$ is nef since effective divisors on an abelian surfaces is nef. But this then implies that $E$ is nef, which contradicts to negativity lemma. 
So $\tilde{X}$ is a K3 or Enriques surface, and $H^1(\tilde{X},\mathcal{O}_{\tilde{X}}) = 0$, which implies $H^1(X,\mathcal{O}_{X}) = 0$ since $X$ has rational singularities. 

If $X$ is not canonical, by \cite[Theorem 3.3]{Tot09}, let $\pi:Z \to X$ be the global index $1$ cover of $X$, then $Z$ is a surface with canonical singularities and trivial canonical bundle. If $Z$ is not smooth, $H^1(Z,\mathcal{O}_Z) = 0$, so $H^1(X,\pi_*\mathcal{O}_Z)=0$ since $\pi$ is finite. But the natural map $\mathcal{O}_X \to \pi_*\mathcal{O}_Z$ admits a section by the (normalised) trace map $\frac{1}{\deg f}tr: \pi_*\mathcal{O}_Z \to \mathcal{O}_X$, by \cite[Proposition 5.7]{KM98}. So $H^1(X,\mathcal{O}_X)$ is a direct summand of $H^1(X,\pi_*\mathcal{O}_Z)=0$, which implies $H^1(X,\mathcal{O}_X)=0$. 

When $Z$ is smooth, we may assume $H^1(Z,\mathcal{O}_Z) >0$. Otherwise $H^1(X,\mathcal{O}_X)$ is a direct summand of $H^1(X,\pi_*\mathcal{O}_Z)=H^1(Z,\mathcal{O}_Z)=0$, which implies $H^1(X,\mathcal{O}_X)=0$. 

In this case, $Z$ is an abelian surface or bielliptic surface, hence $X$ is a quotient of an abelian surface $A$ by a finite group $G$. 

For any variety $X$, denote the irregularity $\dim H^1(X,\mathcal{O}_X) = h^1(X,\mathcal{O}_X)=q(X)$. 

It remains to show $q(X)=0$ in this case. Suppose $q(X)>0$. Since $H^1(X,\mathcal{O}_X)$ is a direct summand of $H^1(Z,\mathcal{O}_Z)$, $q(X)\leq 2$. If $q(X)=2$, then $X$ is again an abelian surface by \cite{Yos95}, which contradicts to the fact that $X$ is singular. Moreover, if $q(X)=1$, by \cite[Lemma 3.4]{Yos95},  $H^0(X,mK_X) = 0$ for all $m>0$, which shows that $X$ is not a Calabi-Yau surface.

\end{proof}

\begin{proof}(of Theorem \ref{Thm generalised cone})

The case when $d=1$ is solved in Lemma \ref{Lem Cone of elliptic}. 

Now we consider $d=2$ case. Fix an embedding $k \to \mathbb{C}$. When $X_{\mathbb{C}}$ is an abelian or bielliptic surface, we proved the generalised cone conjecture in Theorem \ref{Thm generalised abelian} and Theorem \ref{Thm cone of bielliptic}. Otherwise, by classifiaction of smooth Calabi-Yau surfaces and Lemma \ref{Lem irregularity}, $H^1(X_{\mathbb{C}},\mathcal{O}_{X_{\mathbb{C}}})=0$. By flat base change, $H^1(X,\mathcal{O}_{X})=0$. So $\mathcal{P}ic(X)$ is an \'etale group scheme. 

This means that there are only finitely many linear equivalence classes that are numerically equivalent to $0$, so ${\rm Pic}(X)\otimes_\mathbb{Z} \mathbb{R}$ coincides with $N^1(X)$. Now Theorem \ref{Thm generalised cone} follows from the usual cone conjecture over $k$, which is proved in Subsection \ref{Cone non-closed}.

\end{proof}

\section{Proof of Theorem \ref{Lem Main} and Theorem \ref{The Main Theorem}}

In this part, we prove Theorem \ref{Lem Main} and Theorem \ref{The Main Theorem}.

\subsection{Essential finiteness}

For our purpose, we need to introduce certain finiteness conditions for groups. 

\begin{definition}\label{Def finiteness condition}
A group $G$ is said to be essentially finite if any finitely generated subgroup of $G$ is finite. 

Suppose $H \subset G$ is a normal subgroup. $H$ is said to be a normal subgroup of essentially finite index (or an essentially finite index normal subgroup) if $G/H$ is essentially finite. 

In general, suppose $G^\prime \subset G$ is a subgroup. $G^\prime$ is said to be a subgroup of essentially finite index (or an essentially finite index subgroup) if there exists a normal subgroup $H$ of $G$ such that $H$ is contained in $G^\prime$ and $G/H$ is essentially finite. 
\end{definition}

\begin{remark}
    We will see that essentially finite index subgroups behave like the usual finite index subgroups. 
\end{remark}

Recall a well-known result in group theory: 

\begin{lemma}\label{Lem gt}
Any finite index subgroup $H$ of a finitely generated group $G$ is finitely generated. 
\end{lemma}

\begin{proof}
We may assume $G$ is a free group. Then $G$ is the fundamental group of $\lor^n \mathbb{S}^1$, which is the  one point union of $n$ circles. In this case, by the Galois correspondence between fundamental groups and covering spaces, $H$ is the fundamental group of a finite cover of $\lor^n \mathbb{S}^1$, which is a finite CW complex. So $H$ is finitely generated.

\end{proof}

\begin{proposition}\label{Prop extension of essen fin}
Suppose we have an exact sequence $0\to F \to G \to H \to 0$ of groups. Then $G$ is essentially finite if and only if $F$ and $H$ are both essentially finite. 
\end{proposition}

\begin{proof}
(1) Suppose $G$ is essentially finite. Then any subgroup of $F$ is automatically a subgroup of $G$, so $F$ is essentially finite. 

Consider any finitely generated subgroup $S$ of $H$, say generated by $h_1,\dots, h_n$. Choose $g_1,\dots,g_n$ such that the image of $g_i$ is $h_i$. Then by our condition the subgroup $T$ of $G$ generated by $g_1,\dots,g_n$ is finite. But $T$ maps surjectively to $S$, so $S$ is finite. 

(2) Suppose $F$ and $H$ are essentially finite. Let $T$ be a finitely generated subgroup of $G$. Then the image $S$ of $T$ in $H$ is finite. Consider the kernel $K$ of $T\to S$. Then $K$ is finitely generated by Lemma \ref{Lem gt}. But $K$ is a finitely generated subgroup of $F$, so it is finite. Now $T$ is an extension of two finite groups $K$ and $S$, so $T$ is finite.

\end{proof}

\begin{remark}\label{Rem essen fin}
We have the following easy observations: 

$\bullet$ Finite groups are essentially finite.

$\bullet$ Torsion abelian groups are essentially finite. 

$\bullet$ Essentially finite groups with a finite generating set are finite. 
\end{remark}

\subsection{Lifting pseudo-automorphisms}

Let $(X/S,\Delta)$ be a klt Calabi-Yau fiber space. For our purpose, we would like to ask when a pseudo-automorphism of the generic fiber $(X_\eta,\Delta_\eta)$ can be lifted to a pseudo-automorphism of the total space $(X,\Delta)$. We only need the case when $\Delta_\eta = 0$, but we include the general case. 

We need the following easy lemma: 

\begin{lemma}\label{Lem small}
Suppose we have the following diagram
\[\begin{tikzcd}
	X & Y \\
	Z & W
	\arrow["f", dashed, from=1-1, to=1-2]
	\arrow["i"', from=1-1, to=2-1]
	\arrow["j"', from=1-2, to=2-2]
	\arrow["g", dashed, from=2-1, to=2-2]
\end{tikzcd}\]
where $X,Y,Z,W$ are normal varieties over a field $k$ of characteristic $0$, $i,j$ are finite surjective morphisms and $f,g$ are birational maps. 

Then $f$ is small if and only if $g$ is small. 
\end{lemma}

\begin{proof}
If $f$ is not small, then $f$ or $f^{-1}$ contracts some divisor. Without loss of generality, we may assume $f$ contracts a divisor $D$ on $X$. Then $g$ contracts $i(D)$, so $g$ is not small. 

Conversely, if $g$ is not small, then $g$ or $g^{-1}$ contracts some divisor. Without loss of generality, we may assume $g$ contracts a divisor $E$ on $Z$. Then $f$ contracts $i^{-1}(E)$, so $f$ is not small. 

\end{proof}

We also need the following lemma to preserve boundary divisors: 

\begin{lemma}\label{Lem easy lift}
Let $f:X\to S$ be a projective morphism of varieties over a field $k$ of characteristic $0$ and $\Delta$ be a $\mathbb{Q}$-divisor on $X$ such that $(X,\Delta)$ is a pair. Let $(X_\eta,\Delta_\eta)$ be the generic fiber. 

Then ${\rm PsAut}(X/S,\Delta)$  is a finite index subgroup of ${\rm PsAut}(X/S)\cap {\rm PsAut}(X_\eta,\Delta_\eta)$. 
\end{lemma}

\begin{proof}
Let $\Delta = \sum_{1\leq i\leq n} a_i\Delta_i$, $a_i>0$, where $\Delta_i$ maps surjectively to $S$ for $1\leq i\leq k$, and $\Delta_i$ maps to a proper subset of $S$ for $k+1\leq i\leq n$. 
Let $G_1, \dots, G_k$ be all the divisors contained in $f^{-1}(\bigcup_{i=k+1}^{n}f(\Delta_i))$. 

Since elements in ${\rm PsAut}(X/S)\cap K$ are pseudo-automorphisms of $(X_\eta,\Delta_\eta)$, they act by the trivial permutation on $\Delta_i$ for $1\leq i\leq k$. 
Since elements in ${\rm PsAut}(X/S)\cap {\rm PsAut}(X_\eta,\Delta_\eta)$ preserve the fiber space structure $f$, they induce permutations of $G_1, \dots, G_k$. 

Let $K$ be the subgroup of ${\rm PsAut}(X/S)\cap {\rm PsAut}(X_\eta,\Delta_\eta)$ given by 
$$K = \{g\in {\rm PsAut}(X/S)\cap {\rm PsAut}(X_\eta,\Delta_\eta)| g_*G_i = G_i, 1\leq i\leq k\}$$

Then $K \subset {\rm PsAut}(X/S)\cap {\rm PsAut}(X_\eta,\Delta_\eta)$ is of finite index. Then the result holds since we have $K \subset {\rm PsAut}(X/S,\Delta)$. 

\end{proof}

Then we prove the lifting lemma for varieties with canonical singularities: 

\begin{lemma}\label{Lem lifting 1}
Let $f:X\to S$ be a projective morphism of varieties over a field $k$ of characteristic $0$, such that $X$ is canonical and $K_X$ is $\mathbb{Q}$-linearly trivial over $S$. Let $X_\eta$ be the generic fiber, then $X_\eta$ is a (possibly not geometrically irreducible) canonical variety. 

Then ${\rm PsAut}(X/S)$ is a finite index subgroup of ${\rm PsAut}(X_\eta)$. 
\end{lemma}

\begin{proof}
It is clear that ${\rm PsAut}(X/S)$ is a subgroup of ${\rm PsAut}(X_\eta)$, so we only need to show that it is of finite index. 

First we claim that it suffices to prove this after base change to the algebraic closure of $k$. Indeed, consider the closure of the graph of $g$, say $\Gamma \subset X\times_k X$. Then the condition that the projections $p_1 : \Gamma \to X$ and $p_2 : \Gamma \to X$ are isomorphisms is represented by open sets by \cite[Theorem 5.22]{FAG}. In particular, the maximal set such that $g$ restricts to an isomorphism is defined over $k$ and is invariant under base change. 

Also, it suffices to prove this result locally on $S$. In fact, if we have a finite open covering $U_i$ of $S$, and the result holds for the base changes $X_i = X\times_S U_i$, then ${\rm PsAut}(X/S) = \bigcap_i {\rm PsAut}(X_i/U_i)$ is also of finite index in ${\rm PsAut}(X_\eta)$. Possibly shrinking $S$, we may assume $K_X$ is $\mathbb{Q}$-linearly equivalent to $0$. 

Let $q: Z \to X$ be a terminal model of $X$. Then $K_Z = q^*K_X$. 

Let $g$ be an element in ${\rm PsAut}(X_\eta)$, and $g_Z: Z \dashrightarrow Z$ be a lifting to $Z$. 
Since $Z$ is terminal and $K_Z$ is nef, $h := f\circ q: Z \to S$ and $g_Z$ satisfies the conditions of \cite[Corollary 3.54]{KM98}. This means that $g_Z$ is small. 

We summarize the notations in the following diagram: 
\[\begin{tikzcd}
	Z & Z \\
	X & X \\
	& S
	\arrow["{g_Z}", dashed, from=1-1, to=1-2]
	\arrow["q"', from=1-1, to=2-1]
	\arrow["q", from=1-2, to=2-2]
	\arrow["g", dashed, from=2-1, to=2-2]
	\arrow["f"', from=2-1, to=3-2]
	\arrow["f", from=2-2, to=3-2]
\end{tikzcd}\]

We claim that the subgroup of ${\rm PsAut}(X/S)$
$$ H:=\{g \in {\rm PsAut}(X_\eta)|g_Z:=q^{-1}gq \in {\rm PsAut}(Z/S), g_Z(E) = E, \forall E \in Exc(Z/X)\}$$
has finite index in ${\rm PsAut}(X_\eta)$. 

Clearly $H$ is a subgroup of ${\rm PsAut}(X/S)$ since $g_Z = q^{-1}gq$ preserve the exceptional divisors of $Z/X$, which means that $g:X \dashrightarrow X$ is small. 

Let $E_1,\dots, E_m$, $F_1, \dots, F_n$ be all the exceptional divisors of $Z/X$, where $E_1, \dots, E_m$ map surjectively to $S$, and $F_1, \dots, F_n$ maps to a proper subset of $S$. Let $G_1, \dots, G_k$ be all the divisors contained in $h^{-1}(\bigcup_{i=1}^{n}h(F_i))$. Since $g_Z$ is a lifting of a pseudo-automorphism of $X_\eta$, $g_Z$ induces a permutation of the exceptional divisors $E_1, \dots, E_m$ of $Z_\eta/X_\eta$. Since $g_Z$ preserves the fiber space structure $h$, it induces a permutation of $G_1, \dots, G_k$. 

This gives an action of ${\rm PsAut}(X_\eta)$ on the finite set $\{E_1,\dots,E_m, G_1, \dots, G_k\}$ of divisors on $Z$. So the subgroup $H$, which is defined by fixing the subset $\{E_1,\dots,E_m, F_1, \dots, F_m\}$, is of finite index in ${\rm PsAut}(X_\eta)$. So our claim holds. 

Since we have $H \subset {\rm PsAut}(X/S)\subset {\rm PsAut}(X_\eta)$, ${\rm PsAut}(X/S)$ is also a finite index subgroup of ${\rm PsAut}(X_\eta)$. 

\end{proof}

Next we deal with klt pairs, which is much subtler: 

\begin{lemma}\label{Lem lifting 2}
Let $f:X\to S$ be a projective contraction of varieties over a field $k$ of characteristic $0$ and $\Delta$ be a $\mathbb{Q}$-divisor on $X$, such that $(X,\Delta)$ is klt, $K_X+\Delta$ is $\mathbb{Q}$-linearly trivial over $S$ and the coefficients of $\Delta$ are of the form $1-1/n$ for some integer $n$. Let $(X_\eta,\Delta_\eta)$ be the generic fiber, then $(X_\eta,\Delta_\eta)$ is a klt pair. 

Then ${\rm PsAut}(X/S,\Delta)$ is an essentially finite index subgroup of ${\rm PsAut}(X_\eta,\Delta_\eta)$. 
\end{lemma}

\begin{proof}
As in Lemma \ref{Lem lifting 1}, we may assume $k$ is algebraically closed, and $K_X+\Delta$ is $\mathbb{Q}$-linearly equivalent to $0$. 

Fix a section 
$\sigma \in H^0(X,\omega_X^{[m]}(m\Delta))$
, where $m$ is the global index, i.e. the minimal positive integer $n$ such that $n(K_X+\Delta)$ is Cartier and linearly equivalent to $0$. Here $\omega_X$ is the canonical sheaf, and $\omega_X^{[m]}$ means the reflexive hull of the $m$-th power of $\omega_X$. Consider the index $1$ cover of $(X,\Delta)$ corresponding to $\sigma$, say $p: Y \to X$, as defined in \cite[Definition 2.49]{Ko13}, (see also \cite[Definition 5.19]{KM98} for some discussions when $\Delta = 0$). Then $Y$ has canonical singularities, and $K_Y \sim 0$. 

We note that we cannot lift a pseudo-automorphism of $(X_\eta,\Delta_\eta)$ simply by functorial property of index $1$ cover. Indeed, the index $1$ cover also relies on the choice of $\sigma$. Nevertheless, a non-unique lifting is enough for our purpose. So we analyse how $\sigma$ changes under the action of ${\rm PsAut}(X_\eta,\Delta_\eta)$. 

Recall that $k(Y)$ can be described as $k(Y) = k(X)((\sigma \tau^{-m})^{1/m})$, where $\tau$ is any rational section of $\omega_X$. We note that $\tau^{-m}$ is a rational section of $\omega_X^{-[m]}$, and $\sigma \tau^{-m}$ is a rational function on $X$, so our notations make sense. Let $g$ be any element in ${\rm PsAut}(X_{\eta},\Delta_\eta)$, we consider the following diagram: 
\[\begin{tikzcd}
	Y & Y \\
	X & X
	\arrow["p", from=1-1, to=2-1]
	\arrow["p", from=1-2, to=2-2]
	\arrow["g", dashed, from=2-1, to=2-2]
\end{tikzcd}\]

Then the rational map $g\circ p$ corresponds to the field extension $k(X)((g^*\sigma(g^*\tau)^{-m})^{1/m})$, where $g^*$ here denotes the map induced by pullback of differential forms. 

Note that $g$ preserve the boundary divisor $\Delta$. So $g^*$ induces an isomorphism from $H^0(X_\eta, \omega_{X_\eta}^{[m]}(m\Delta_\eta))$ to $H^0(X_\eta, \omega_{X_\eta}^{[m]}(mg^*\Delta_\eta)) = H^0(X_\eta, \omega_{X_\eta}^{[m]}(m\Delta_\eta))$. This shows that $g^*\sigma$ and $\sigma$ are both global sections of $\omega_{X_\eta}^{[m]}(m\Delta_\eta)$. Since $H^0(X_\eta,\omega_{X_\eta}^{[m]}(m\Delta_\eta))$ is of dimension $1$ as a $k(S)$-vector space, $g^*\sigma$ and $\sigma$ differ by an invertible element of $k(S)$, denoted by $g^*\sigma/\sigma \in k(S)^*$. 
Moreover, since $g^*$ is a $k(S)$-linear map, the map $${\rm PsAut}(X_\eta,\Delta_\eta) \to k(S)^*/k(S)^{*m}, g \mapsto g^*\sigma/\sigma$$ is a group homomorphism. Indeed, we have $$(gh)^*\sigma/\sigma = h^*(g^*\sigma)/\sigma = h^*(\sigma \cdot g^*\sigma/\sigma)/\sigma = h^*\sigma/\sigma \cdot g^*\sigma/\sigma= g^*\sigma/\sigma \cdot h^*\sigma/\sigma$$ 

Let $K$ be the kernel of ${\rm PsAut}(X_\eta,\Delta_\eta) \to k(S)^*/k(S)^{*m}$. $K$ is a essentially finite index subgroup of ${\rm PsAut}(X_\eta,\Delta_\eta)$ since $k(S)^*/k(S)^{*m}$ is essentially finite by Remark \ref{Rem essen fin}. Then for $g\in K$, $g$ can be (non-uniquely) lifted to a birational automorphism $g_Y$ of $Y$ such that the following diagram commutes: 
\[\begin{tikzcd}
	Y & Y \\
	X & X
	\arrow["{g_Y}", dashed, from=1-1, to=1-2]
	\arrow["p", from=1-1, to=2-1]
	\arrow["p", from=1-2, to=2-2]
	\arrow["g", dashed, from=2-1, to=2-2]
\end{tikzcd}\]

Let ${G}$ be the subgroup of ${\rm PsAut}(Y_\eta) \subset {\rm Aut}(k(Y)/k(S))$ given by 
$$ G = \{g_Y \in {\rm PsAut}(Y_\eta)| \text{ there exists } g\in K \subset {\rm PsAut}(X_\eta,\Delta_\eta), \text{ such that } g\circ p = p\circ g_Y\}$$
Since $k(Y)/k(X)$ is Galois, we have an exact sequence (note that ${\rm Gal}(k(Y)/k(X)) = {\rm Aut}(Y/X)$): 
$$1 \to {\rm Aut}(Y/X) \to {G} \to K \to 1$$
Note that an element $g$ in ${G}$ is a pseudo-automorphism of $Y$ if and only if the image in $K \subset {\rm Bir}(X)$ is a pseudo-automorphism of $X$, by Lemma \ref{Lem small}. So we have a commutative diagram whose rows are exact: 
\[\begin{tikzcd}
	1 & {{\rm Aut}(Y/X)} & {H={G}\cap {\rm PsAut}(Y/S)} & {K\cap {\rm PsAut}(X/S)} & 1 \\
	1 & {{\rm Aut}(Y/X)} & {G={G}\cap {\rm PsAut}(Y_\eta)} & {K} & 1
	\arrow[from=1-1, to=1-2]
	\arrow[from=1-2, to=1-3]
	\arrow["\simeq", from=1-2, to=2-2]
	\arrow[from=1-3, to=1-4]
	\arrow[hook, from=1-3, to=2-3]
	\arrow[from=1-4, to=1-5]
	\arrow[hook, from=1-4, to=2-4]
	\arrow[from=2-1, to=2-2]
	\arrow[from=2-2, to=2-3]
	\arrow[from=2-3, to=2-4]
	\arrow[from=2-4, to=2-5]
\end{tikzcd}\]
To show $K\cap {\rm PsAut}(X/S)$ is a finite index subgroup of $K$, it suffices to show $H$ is a finite index subgroup of $G$, which reduces to show that ${\rm PsAut}(Y/S)$ is a finite index subgroup of ${\rm PsAut}(Y_\eta)$. This is shown in the Lemma \ref{Lem lifting 1}. 

So we know that $K \subset {\rm PsAut}(X_\eta,\Delta_\eta)$ is of essentially finite index, and $K\cap {\rm PsAut}(X/S)\subset K$ is of finite index. Hence ${\rm PsAut}(X/S)\cap K \subset {\rm PsAut}(X_\eta,\Delta_\eta)$ is of essentially finite index. In particular, ${\rm PsAut}(X/S)\cap {\rm PsAut}(X_\eta,\Delta_\eta) \subset {\rm PsAut}(X_\eta,\Delta_\eta)$ is of essentially finite index. 

Finally, by Lemma \ref{Lem easy lift},  ${\rm PsAut}(X/S,\Delta) \subset {\rm PsAut}(X/S)\cap {\rm PsAut}(X_\eta,\Delta_\eta)$ is of finite index. 
So we conclude that ${\rm PsAut}(X/S,\Delta) \subset {\rm PsAut}(X_\eta,\Delta_\eta)$ is of essentially finite index.

\end{proof}

We deduce the following lifting lemma:

\begin{lemma}\label{Lem Aut of elliptic}
Let $(X,\Delta)$ be a projective klt Calabi-Yau pair over a field $k$ of characteristic $0$, where $X$ is $\mathbb{Q}$-factorial. Assume $-K_X$ is semiample, inducing a fibration $f: X\to S$. Let $X_\eta$ be the generic fiber of $f$. 
Viewing ${\rm PsAut}(X_\eta)$ and ${\rm PsAut}(X,\Delta)$ as subgroups of ${\rm Bir}(X)$, we set $G = {\rm PsAut}(X_\eta) \cap {\rm PsAut}(X,\Delta)$. 

Then $G$ is of essentially finite index in ${\rm PsAut}(X_\eta)$. 
\end{lemma}

\begin{proof}

Since $K_X = f^*A$ for some ample $\mathbb{Q}$-divisor $A$ on $S$, $K_X$ is $\mathbb{Q}$-linearly trivial over $S$. In this case, the restriction of $\Delta$ to the generic fiber is $0$. Now the conditions of Lemma \ref{Lem lifting 2} are satisfied (we only use the case when $\Delta = 0$), so ${\rm PsAut}(X/S)$ is an essentially finite index subgroup of ${\rm PsAut}(X_\eta)$. 

By Lemma \ref{Lem easy lift} in the case $\Delta_\eta = 0$, ${\rm PsAut}(X/S,\Delta)  \subset {\rm PsAut}(X/S)$ is of finite index. So ${\rm PsAut}(X/S,\Delta)$ is an essentially finite index subgroup of ${\rm PsAut}(X_\eta)$. 

Then the result holds since ${\rm PsAut}(X/S,\Delta) \subset G$. 
    
\end{proof}

\subsection{Proof of Theorem \ref{Lem Main}}

\label{Sec fibration and cone}

Then we consider the implication of the cone conjecture from the generic fiber to the total space. 

We need the following observation: 

\begin{lemma}\label{Lem Pic restriction}
Let $X \to S$ be a projective morphism of normal quasi-projective varieties over a finitely generated extension of $\mathbb{Q}$, say $k$. Let $\eta$ be the generic fiber of $S$. 

Then the kernel of the restriction map $\iota^*: {\rm Pic}(X)\otimes_{\mathbb{Z}}\mathbb{R} \to {\rm Pic}(X_\eta)\otimes_{\mathbb{Z}}\mathbb{R}$ is generated by vertical divisors.

\end{lemma}

\begin{proof}

Let $D \in {\rm Pic}(X)\otimes_{\mathbb{Z}}\mathbb{R}$ be an $\mathbb{R}$-divisor on $X$ such that $\iota^*D = D_\eta = 0$. Then $D_\eta$ is $\mathbb{R}$-linearly equivalent to $0$, i.e. there exists rational functions $\varphi_1, \dots, \varphi_n$ on $X_\eta$ and real numbers $a_1, \dots, a_n$ such that $D_\eta = \sum_{i=1}^n a_i(\varphi_i)$, where $(\varphi_i)$ denotes the principal divisor defined by $\varphi_i$. 

Now we view $\varphi_i$ as rational functions on $X$. 
Then $D - \sum_{i=1}^n a_i(\varphi_i)$ is an $\mathbb{R}$-divisor on $X$, which is $0$ when restricted to $X_\eta$. So $D - \sum_{i=1}^n a_i(\varphi_i)$ is a vertical divisor. 

This shows that $D$ is represented by a vertical divisor in ${\rm Pic}(X)\otimes_{\mathbb{Z}}\mathbb{R}$. 
    
\end{proof}

Then we are ready to prove Theorem \ref{Lem Main}:

\begin{proof} (of Theorem \ref{Lem Main})

We devide the proof into 5 steps. 

In Step 1, we apply our lifting result to prove that an analog of the movable cone conjecture holds for $X_\eta$ with smaller automorphism group $G$. 
In Step 2, we introduce a boundary divisor with larger support, in order to apply the geography of log models. 
In Step 3, we apply the GLM for movable divisors. 
In Step 4, we prove that the inverse image of a rational polytope $P \subset \rm{Mov}^e(X_\eta)$ in $\rm{Mov}^e(X)$ is still rational polyhedral, using detailed analysis of convex geometry. 
In Step 5, we combine the previous steps and finish the proof. 

\hfill

\noindent \textbf{Step 1. } 
In this step, we prove that there exists a rational polyhedral cone $P$ such that $GP = {\rm Mov}^e(X_\eta)$. 

Let $G^*$ be the image of $G$ in ${\rm GL}({\rm Pic}(X_\eta) \otimes_\mathbb{Z} \mathbb{R})$, and ${\rm PsAut}(X_\eta)^*$ be the image of ${\rm PsAut}(X_\eta)$ in ${\rm GL}({\rm Pic}(X_\eta) \otimes_\mathbb{Z} \mathbb{R})$. 

We claim that $G^*$ is of finite index in ${\rm PsAut}(X_\eta)^*$. Indeed, $G^*$ is of essentially finite index in ${\rm PsAut}(X_\eta)^*$ by Lemma \ref{Lem Aut of elliptic} and Proposition \ref{Prop extension of essen fin}. Since $({\rm Mov}^e(X_\eta), {\rm PsAut}(X_\eta)^*)$ is of polyhedral type, by \cite[Corollary 4.15]{Loo14}, ${\rm PsAut}(X_\eta)^*$ is finitely generated.

By definition, there exists a normal subgroup $H$ of ${\rm PsAut}(X_\eta)^*$ contained in $G^*$ such that  ${\rm PsAut}(X_\eta)^*/H$ is essentially finite. But ${\rm PsAut}(X_\eta)^*/H$ is a finitely generated, so it is finite by Remark \ref{Rem essen fin}. This shows our claim.

Let ${\rm PsAut}(X_\eta)^* = \sqcup_{i\in I} G^*\sigma_i$ be a finite right coset decomposition. 
Replacing $P$ by the cone generated by $\{\sigma_i P\}_{i\in I}$, we may assume $GP = {\rm Mov}^e(X_\eta)$. 

\hfill

\noindent \textbf{Step 2. } 
In this step, we introduce a new boundary divisor $\Delta^\prime$ with larger support. 

Let $\{G_i\}_{i\in I}$ be a finite set of generators of $P$. Let $\{\tilde{G}_i\}_{i\in I}$ be a set of liftings to ${\rm Pic}(X)\otimes_\mathbb{Z} \mathbb{R}$. Indeed, write $G_i$ as a sum of prime divisors with positive coefficients on $X_\eta$, then taking the closures of the prime divisors, we get a sum of Weil divisors on $X$. Since $X$ is $\mathbb{Q}$-factorial, the resulting divisors $\tilde{G}_i$ are effective and $\mathbb{Q}$-Cartier. 

Let $\{D_j\}_{j\in J}$ be all prime divisors appearing in $\{\tilde{G}_i\}_{i\in I}$. By Lemma \ref{Lem Pic restriction}, we may pick vertical prime divisors $F_1,\dots,F_n\in {\rm Pic}(X)$ such that $\{D_j\}_{j\in J}$ and $F_1,\dots,F_n$ spans ${\rm Pic}(X)\otimes_\mathbb{Z} \mathbb{R}$. Pick a $\mathbb{Q}$-divisor $B\sim_{\mathbb{Q}} A$ such that the support of $B$ contains $f(F_i)$, $1\leq i\leq n$. 

Now we consider a new pair $(X,(1-\varepsilon_0)\Delta+\varepsilon_0 f^*B)$. For sufficiently small $\varepsilon_0>0$, $(X,(1-\varepsilon_0)\Delta+\varepsilon_0 f^*B)$ is still klt, and it is clearly a Calabi-Yau pair. We fix a sufficiently small rational number $\varepsilon_0>0$, and write $\Delta^{\prime} = (1-\varepsilon_0)\Delta+\varepsilon_0 f^*B$. 

\hfill

\noindent \textbf{Step 3. } 
In this step, we apply the GLM for some movable divisors. 

Let $S = \{D_j\}_{j\in J} \cup \{D\}_{D\subset\rm{Supp}(\Delta^\prime)}$, and $F = \sum_{D\in S}D$, and we apply strong GLM, see Conjecture \ref{Conj strong GLM}, which holds by Corollary \ref{strong GLM}. 

Recall that $B_F = \oplus_{D\in S}[0,1]D$. Define $M_F = \{D\in B_F| [K_X + D]\in \mathcal{M}^e(X),(X,D) \text{ is lc}\}$, $N_F = \{D\in B_F| [K_X + D]\in \mathcal{B}^e(X),(X,D) \text{ is lc}\}$. Now the subset $M_F$ of $B_F$ is precisely the union of the wlc equivalence classes that admit small wlc models, say $M_F = \sqcup_i C_i$. Hence its closure $\overline{M_F}$ is rational polyhedral. 

Moreover, we claim that $\overline{M_F} = M_F$. In fact, by strong GLM, see Conjecture \ref{Conj strong GLM}, $N_F$ is closed. So we have $\overline{M_F} = \overline{M_F\cap N_F} \subset \overline{M_F}\cap \overline{N_F} = \overline{M_F}\cap N_F = M_F$. The converse inclusion is clear. 

So we conclude that $M_F$ is rational polyhedral. 

\hfill

\noindent \textbf{Step 4. } 
In this step, we prove that $\iota^{*-1}(P) \cap {\rm Mov}^e(X)$ is rational polyhedral. Here $\iota^*: {\rm Pic}(X)\otimes_\mathbb{Z} \mathbb{R} \to {\rm Pic}(X_\eta)\otimes_\mathbb{Z} \mathbb{R}$ is the restriction map. 

We fix notations by the following commutative diagram: 
\[\begin{tikzcd}
	{{\rm Mov}^e(X)} & {{\rm Pic}(X)\otimes_\mathbb{Z} \mathbb{R}} \\
	{{\rm Mov}^e(X_\eta)} & {{\rm Pic}(X_\eta)\otimes_\mathbb{Z} \mathbb{R}}
	\arrow[hook, from=1-1, to=1-2]
	\arrow["{\iota^*|_\mathcal{M}}"', from=1-1, to=2-1]
	\arrow["{\iota^*}", from=1-2, to=2-2]
	\arrow[hook, from=2-1, to=2-2]
\end{tikzcd}\]

Let $W$ be the cone in ${\rm Pic}(X)\otimes_\mathbb{Z} \mathbb{R}$ generated by the image of $K_X+{M_F}$ via the natural map $K_X+{M_F} \to {\rm Pic}(X)\otimes_\mathbb{Z} \mathbb{R}$, $K_X+D \mapsto K_X+D$. Then $W$ is generated by finitely many divisors as ${M_F}$ is a rational polytope. 

We claim that $$\iota^{*-1}(P)\cap {\rm Mov}^e(X)= \iota^{*-1}(P) \cap W$$
Clearly RHS is contained in LHS since $W\subset {\rm Mov}^e(X)$. 

For the converse inclusion, let $D$ be any effective $\mathbb{R}$-divisor in $\iota^{*-1}(P)\cap {\rm Mov}^e(X) = (\iota|_{\mathcal{M}})^{*-1}(P)$. We will show $D\in W$. 

Then $\iota^{*}D\in P$, so we may write $\iota^{*}D = \sum_{j\in J} a_j\iota^{*}D_j$ with $a_j\geq 0$ since $P$ is contained in the cone generated by $\{\iota^{*}D_j\}_{j\in J}$. 
So $$D= \sum_{j\in J} a_jD_j + \sum_{i=1}^n b_i F_i$$ where $F_i$ are vertical divisors chosen before. 
Since $K_X+\Delta^\prime$ is $\mathbb{Q}$-linearly trivial, for any $\varepsilon>0$, we have $$\varepsilon D\sim_\mathbb{Q} \sum_{j\in J} \varepsilon a_j{D}_j + K_X+\Delta^\prime+ \varepsilon \sum_{i=1}^n b_i F_i$$ 
Note that the support of $\Delta^\prime$ contains $F_1, \dots,F_n$, so when $\varepsilon$ is sufficiently small, all coefficients in $$\sum_{j\in J} \varepsilon a_j{D}_j + \Delta^\prime+ \varepsilon \sum_{i=1}^n b_i F_i$$ lies in $[0,1]$, and $\left(X,\Delta^\prime+\varepsilon \left(\sum_{j\in J} a_jD_j + \sum_{i=1}^n b_i F_i\right)\right) $ is klt as $(X,\Delta^\prime)$ is klt. This means that $\varepsilon D$ lies in the image of $K_X+M_F$. Hence $D\in W$. So the converse inclusion holds. 

So we have shown that $\iota^{*-1}(P)\cap {\rm Mov}^e(X)= \iota^{*-1}(P) \cap W$ is an intersection of two rational polyhedral cones, hence it is itself rational polyhedral. 

\hfill

\noindent \textbf{Step 5. } In this step, we finish the proof. 

Let $\Pi = \iota^{*-1}(P)\cap {\rm Mov}^e(X)$, and we claim that $\bigcup_{g\in {\rm PsAut}(X,\Delta)}
g_* \Pi = {\rm Mov}^e(X)$. 

Indeed, it is clear that $\bigcup_{g\in {\rm PsAut}(X,\Delta)} g_* \Pi \subset {\rm Mov}^e(X)$. 
Also, we have the converse inclusion

\begin{align*}
\bigcup_{g\in {\rm PsAut}(X,\Delta)}
g_* \Pi & \supset \bigcup_{g\in G}
g_* (\iota^{*-1}(P)\cap {\rm Mov}^e(X)) \\ & = (\bigcup_{g\in G}
g_* \iota^{*-1}(P))\cap {\rm Mov}^e(X) \\ & =   \iota^{*-1}(\bigcup_{g\in G}
(g_*P))\cap {\rm Mov}^e(X)
\\ & =  \iota^{*-1}{\rm Mov}^e(X_\eta)\cap {\rm Mov}^e(X)  \\ & =  {\rm Mov}^e(X)
\end{align*} 

Finally, consider the projection $\pi: {\rm Pic}(X)\otimes_\mathbb{Z} \mathbb{R} \to N^1(X)$, and let $Q$ be the image of $\Pi$. Then ${\rm PsAut}(X,\Delta)Q = \mathcal{M}^e(X)$ by taking the image of $\pi$ of the previous equality.

Thus the movable cone conjecture holds for $(X,\Delta)$ by Theorem \ref{Thm equivalent cone}. Then the cone conjecture holds by Section \ref{Relations between cone conjectures}.

\end{proof}

\subsection{Proof of Theorem \ref{Thm main dimension 3}}
\label{Subsection conclusion}

\begin{proof}(of Theorem \ref{Thm main dimension 3})

First we note that Assumption \ref{Assumption main} holds by Proposition \ref{Prop Assumption}. 

By Section \ref{Relations between cone conjectures}, it suffices to prove the movable cone conjecture. 

By Lemma \ref{Lem passage to limit}, we may assume $k$ is finitely generated over $\mathbb{Q}$. 

Then we run $\mathbb{Q}$-factorial $K_X+\Delta+\varepsilon\Delta$-MMP on $X$ for small $\varepsilon$, in the sense of \cite[3.31]{KM98}. This MMP terminates with $(X^{min}, (1+\varepsilon )\Delta^{min})$ by Proposition \ref{Prop non-closed mmp}. Let $\phi_i : (X_i,\Delta_i) \dashrightarrow  (X_{i+1},\Delta_{i+1})$ be a step of the MMP. If $\phi$ is a divisorial contraction, by Theorem \ref{main thm}, the movable cone conjecture holds for $(X_i,\Delta_i)$ if and only if the movable cone conjecture holds for $(X_{i+1},\Delta_{i+1})$. If $\phi$ is a flip, by definition, the movable cone conjecture holds for $(X_i,\Delta_i)$ if and only if the movable cone conjecture holds for $(X_{i+1},\Delta_{i+1})$.  

So we conclude that it suffices to prove the movable cone conjecture for the minimal model $(X^{min}, \Delta^{min})$, i.e. we may assume $\Delta$ is nef by replacing $(X,\Delta)$ with $(X^{min}, \Delta^{min})$. We note that log abundance holds for klt pairs when the Kodaira dimension of $K_X+\Delta$ is at least $d-3$ by \cite[Lemma 5.6]{KMM94} (see also \cite[Theorem 3]{Fuk02}), so $\Delta$ is semiample. 

Then the result follows from Theorem \ref{Lem Main}, where the Condition (1) is verified by Theorem \ref{Thm generalised cone} and the Condition (2) is verified by Lemma \ref{Lem Aut of elliptic}.

\end{proof}

\appendix

\section{Geography of log models}

In this appendix, we establish the geography of log models for NQC generalised pairs. The ideas are the same as \cite{Choi08}, \cite{SC11}. 

We also note that this proof is self contained in the sense that we only use basic results for minimal model program of generalised pairs. In particular, when $M = 0$, this proof might be considered as a reformulation of the proof of \cite{Choi08}, \cite{SC11}. We state the assumption on existence of minimal models in a more precise way, which is needed for the main theorem in this paper. 

We also note that some part of the Shokurov polytope for NQC generalised pairs is already established in \cite[Subsection 3.3]{HL22}. 

The structure of this appendix is as follows: 

Subsection A.1 gives the proof of geography of log models for NQC generalised pairs over $\mathbb{C}$. Some lemmas used in Subsection A.1 are proved in Subsection A.2-A.5. 

Subsection A.6 gives the proof of geography of log models for NQC generalised pairs over any field of characteristic $0$, using Galois descent. Some lemmas used in Subsection A.6 are proved in Subsection A.7.

\subsection{GLM over $\mathbb{C}$}

In this subsection, we work over $\mathbb{C}$. 

Firstly, we recall the conjecture of the geography of log models for usual pairs. 

\begin{conjecture}
Let $X \to S$ be a projective morphism of normal quasi-projective varieties, and $F = \sum_i D_i$ be a reduced divisor on $X$. 
Let $B_F \subset \bigoplus_i [0,1]D_i$ be a rational polytope and $N_F$ be the subset of divisors $B$ in $B_F$ such that $(X,B)$ is lc and $K_X+B$ is pseudo-effective. Then: 

(1) $N_F$ is closed convex rational polyhedral. 

(2) There exists a finite decomposition $N_F = \bigsqcup_i \Pi_i^\circ$ such that each $\Pi_i$ is closed convex rational polyhedral and any divisors $D_1$ and $D_2$ in $\Pi_i$ admit the same wlc models. 

\end{conjecture}

We may state similar conjecture for generalised pairs, following \cite{HL22}. 

\begin{conjecture}\label{Conj main}
Let $X \to S$ be a projective morphism of normal quasi-projective varieties, and $F = \sum_i D_i$ be a reduced divisor on $X$. Let $\{M_j\}_{j\in J}$ be b-nef $\mathbb{Q}$-divisors on $X$. 
Let $B_{F+M} \subset \bigoplus_i [0,1]D_i \oplus \bigoplus_j [0,1]M_j$ be a rational polytope and $N_{F+M}$ be the subset of b-divisors $B+M$ in $B_{F+M}$ such that $(X,B+M)$ is g-lc and $K_X+B+M$ is pseudo-effective. Then: 

(1) $N_{F+M}$ is closed convex rational polyhedral. 

(2) There exists a finite decomposition $N_{F+M} = \bigsqcup_i \Pi_i^\circ$ such that each $\Pi_i$ is closed convex rational polyhedral and any divisors $D_1+N_1$ and $D_2+N_2$ in $\Pi_i^\circ$ admit the same wlc models. 

\end{conjecture}

\begin{remark}
We refer to \cite{HL22} for the necessary definitions for generalised pairs. 

More precisely, for the definition of generalised pair, we refer to \cite[Definition 2.1]{HL22}. For generalised dlt pairs, we refer to \cite[Definition 2.2]{HL22}. For NQC (nef $\mathbb{Q}$-Cartier combinations) generalised pair, we refer to \cite[Definition 2.13]{HL22}. For wlc (weak log canonical) model of $(X,B+M)$, we mean a model $(X,B+M) \dashrightarrow Y$ that satisfies  \cite[Definition 2.9(1)-(3)]{HL22} and satisfies \cite[Definition 2.9(5)]{HL22} with non-strict inequality. 

For the definition of log birational transform of a usual divisor $B$, we refer to \cite[Formula 2.5]{HL22}, which means the sum of the strict transform and the reduced exceptional divisor. By log birational transform of $B+M$, we mean the sum of log birational transform of $B$ and the trace of the b-divisor $M$. This is compatible with \cite[Definition 2.11(2)]{HL22}. 
\end{remark}

\begin{remark}
A list of notations: 

$\Pi^\circ$: relative interior of $\Pi$; 

$\mathbb{B}(p,r)$: ball with center $p$ and radius $r$; 

Let $(X,B+M)$ be a g-pair and $X \dashrightarrow Y$ be a birational map, we denote the log birational transform of $B+M$ on $Y$ by $B_Y+M_Y$. 
When $Y = X^\prime$, we may simplify the notation. Let $X \dashrightarrow X^\prime$ be a birational map, we may denote the log birational transform of $B+M$ on $X^\prime$ by $B^\prime+M^\prime$. 

\end{remark}

\begin{theorem}\label{Theorem A.4}
Existence of minimal models in dimension $d$ for generalised pairs implies Conjecture \ref{Conj main} in dimension $d$. 

More precisely, we may fix the data $(X/S,F,\{M_j\}_{j\in J},B_{F+M})$ in Conjecture \ref{Conj main}. Assume existence of wlc models for generalised pairs $(X/S,B+M)$, where $B+M \in N_{F+M}$. Then Conjecture \ref{Conj main} holds for $(X/S,F,\{M_j\}_{j\in J},B_{F+M})$. 
\end{theorem}


\begin{proof}

We proceed by induction on the dimension of $B_{F+M}$. 

\hfill

\noindent \textbf{Step} 0: Do some easy reductions. 

We may assume that $(X,B+M)$ is g-lc for all $B+M\in B_{F+M}$, since lc condition is rational polyhedral, which can be observed by taking a log resolution. 

By compactness, it suffices to show that for any $B_0+M_0\in B_{F+M}$, there exists $\epsilon>0$ and some rational polytope $B_{F+M,0} \subset B_{F+M}$ such that $B_{F+M}\cap \mathbb{B}(B_0,\epsilon) \subset B_{F+M,0}$ and the result holds for $B_{F+M,0}$. This means that we may shrink $B_{F+M}$ freely around $B_0+M_0$. 

We may assume $X$ is $\mathbb{Q}$-factorial by taking a $\mathbb{Q}$-factorial g-dlt modification, whose existence is established in \cite[Proposition 3.10]{HL22}. In fact, for each $B_0+M_0\in B_{F+M}$, we may find some $\mathbb{Q}$-factorial g-dlt modification $X^\prime \to X$ of $(X,B_0+M_0)$. 
We may assume the log discrepancies of the exceptional divisors are $0$. 
Then it is a $\mathbb{Q}$-factorial g-dlt modification of $(X,B+M)$ for $B+M$ in an open neighbourhood of $B_0+M_0$, using the criterion in \cite[Proposition 2.40]{KM98}. 
Then the result holds since we may shrink freely around $B_0+M_0$. 

If $K_X+B_0+M_0$ is not pseudo-effective over $S$, the result is clear as pseudo-effectivity is a closed condition. 

Otherwise, $K_X+B_0+M_0$ is pseudo-effective over $S$. 

\hfill

\noindent \textbf{Step} 1: Reduce to the case when $K_X+B_0+M_0$ is nef. 

We may run 
{\color{black} MMP} 
on $K_X+B_0+M_0$ over $S$, which exists by Proposition \ref{Prop special termination}, and end up with a model $\phi: X \dashrightarrow X^\prime/S$. Then $K_{X^\prime}+B^\prime_0+M^\prime_0$ is nef. Possibly shrinking $B_{F+M}$, we may assume $\phi$ is $K_X+B+M$-negative for all $B\in B_{F+M}$. 

Let $F^\prime = \phi_*F$, $M^\prime_j = \phi_*M_j$ and $B_{F^\prime+M^\prime} = \phi_*B_{F+M}$. 

Firstly, we need to verify existence of wlc models for generalised pairs $(X^\prime/S,B^\prime+M^\prime)$, where $B^\prime+M^\prime \in N_{F^\prime+M^\prime}$. Indeed, by assumption, we know existence of wlc models for generalised pairs $(X/S,B+M)$, $B+M \in B_{F+M}$, which are automatically a wlc models for generalised pairs $(X^\prime/S,B^\prime+M^\prime)$, $B^\prime+M^\prime \in B_{F^\prime+M^\prime}$ by Lemma \ref{lem wlc equiv}. 

Assume that our result holds for nef case, then there exists a closed convex rational polytope $N_{F^\prime+M^\prime}$ and a decomposition $N_{F^\prime+M^\prime} = \bigsqcup_i \Pi^{\prime \circ}_i$ satisfying the Conjecture \ref{Conj main}. Possibly shrinking $B_{F+M}$, we may assume each $\Pi_i^\prime$ contains $B_0^\prime+M_0^\prime$. 

Let $N_{F+M} = \phi_*^{-1} N_{F^\prime}$ and $\Pi_i =  \phi_*^{-1}\Pi^{\prime}_i$. It's clear that $B\in N_{F+M}$ iff $K_X+B+M$ is pseudo-effective by negativity of $\phi$. We verify that $\Pi_i$ gives the desired decomposition, i.e. 
suppose $B_1+M_1, B_2+M_2 \in \Pi_i^\circ$, and $Z$ is a wlc model of $(X,B_1+M_1)$, then $Z$ is a wlc model of $(X,B_2+M_2)$. 

Indeed, by Lemma \ref{lem wlc equiv}, $Z$ is a wlc model of $(X^\prime,B_1^\prime+M_1^\prime)$. By our assumption for the nef case, $Z$ is a wlc model of $(X^\prime,B_2^\prime+M_2^\prime)$. By Lemma \ref{lem wlc equiv} again, $Z$ is a wlc model of $(X,B_2+M_2)$.

\hfill

\noindent \textbf{Step} 2: Application of stability results. 

From now on, we assume $K_X+B_0+M_0$ is nef. 

By Proposition \ref{Prop stability of wlc models}, possibly shrinking $B_{F+M}$,  for any $B+M\in B_{F+M}$, any wlc model of $K_X+B+M$ is crepant with respect to $K_X+B_0+M_0$. 

By Corollary \ref{Cor Bir11 Prop 3.2}(1), there exists nonnegative 
real numbers $a_1,\dots,a_r$, b-divisors $B_{0,i} +M_{0,i} \in B_{F+M}$, $1\leq i\leq r$, 
and $ m\in \mathbb{N}$ 
depending only on $X/S$, $F$, $\{M_j\}_{1\leq j\leq l}$ and $(B_0+M_0,B_{F+M})$ so that we have $\sum_{i=1}^r a_i = 1$, $B_0+M_0 = \sum_{i=1}^r a_i(B_{0,i}+M_{0,i})$, $a_1, \dots,a_r$ are linearly independent over $\mathbb{Q}$, and each $m(K_X + B_{0,i}+M_{0,i})$ is Cartier and nef. 

By induction hypothesis, for any proper face $B_{F+M}^{k}$ such that $B_1+M_1$ is not contained in $B_{F+M}^k$, $N_{F+M}\cap B_{F+M}^k$ is rational polyhedral. Moreover, there exists a finite decomposition for the relative interior of $B_{F+M}^k$, say $N_{F+M}^{k} \cap B_{F+M}^{k\circ} = \bigsqcup_i \Sigma_i^{k\circ}$, such that any divisors $D_1+N_1$ and $D_2+N_2$ in $\Sigma_i^{k\circ}$ admit the same wlc models. 

We note that the following two claims correspond to the locally conical property of the decomposition in \cite{Choi08}, \cite{SC11}. 

\hfill

\noindent \textbf{Claim} 2.1:  We have $N_{F+M} \cap B_{F+M} = \{B_1+M_1\} \sqcup\bigsqcup_k(\bigsqcup_i \text{Conv}(B_1+M_1, \Sigma_i^{k\circ})^\circ)$. 

Indeed, the right hand side is contained in the left hand side since convex combination of pseudo-effective divisors are pseudo-effective. On the other hand, suppose that $B+M\in N_{F+M} \cap B_{F+M}$ and $B+M\neq B_1+M_1$, then there exists some positive real number $c$ such that $B_{(c)}+M_{(c)} = B_1+M_1+c(B+M-B_1-M_1)$ lies in $B_{F+M}^{k\circ}$ for some $k$. 

To prove our claim, it suffices to show that $K_X+B_{(c)}+M_{(c)}$ is pseudo-effective. Suppose that $K_X+B_{(c)}+M_{(c)}$ is not pseudo-effective, then by \cite[Lemma 4.4(1)]{BZ16}, we may run $K_X+B_{(c)}+M_{(c)}$-MMP which terminates with a Mori fiber space. By Corollary \ref{Cor Bir11 Prop 3.2}(3), this MMP is $K_X+B_0+M_0$-trivial, and by linear independence of $a_i$, this MMP is $K_X+B_{0,i}+M_{0,i}$-trivial for $1\leq i\leq r$. So this MMP is also a $K_X+B+M$-MMP, which leads to a contradiction since $K_X+B+M$ pseudo-effective. So Claim 2.1 holds. 

\hfill

\noindent \textbf{Claim} 2.2:  Any divisors $D_1+N_1$ and $D_2+N_2$ in $\text{Conv}(B_1+M_1, \Sigma_i^{k\circ})^\circ$ admit the same wlc models. 

Indeed, there exists some real number $e>1$ such that $D_e+N_e = B_1+M_1+e(D_1+N_1-B_1-M_1)$ lies in $\Sigma_i^{k\circ}$. It suffices to show that $D_1+N_1$ and $D_e+M_e$ admit the same wlc models. 

Suppose $X \dashrightarrow X^\prime$ gives a wlc model of $(X, D_e+N_e)$. Then by our assumption on $B_{F+M}$, $X \dashrightarrow X^\prime$ is crepant with respect to $K_X+B_1+M_1$. So $X \dashrightarrow X^\prime$ gives a wlc model of $(X, D_1+N_1)$ since $D_1+N_1$ is a convex combination of $B_1+M_1$ and $D_e+N_e$. 

On the other hand, suppose $X \dashrightarrow X^\prime$ gives a wlc model of $(X, D_1+N_1)$. By our assumption on $B_{F+M}$, $X \dashrightarrow X^\prime$ is crepant with respect to $K_X+B_1+M_1$. So it remains to show that the log birational transform of $K_X+D_e+N_e$ is nef on $X^\prime$. 

By Lemma \ref{Lem realization of wlc models}, there exists a $\mathbb{Q}$-factorial dlt modification $(Y,D_{Y,1}+N_{Y,1})$ of $(X,D_{1}+N_{1})$, a sequence of $K_Y+D_{Y,1}+N_{Y,1}$-MMP $(Y,D_{Y,1}+N_{Y,1}) \dashrightarrow (Z,D_{Z,1}+N_{Z,1})$ over $S$, a sequence of $K_Z+D_{Z,1}+N_{Z,1}$-flops $(Z,D_{Z,1}+N_{Z,1}) \dashrightarrow (Z^\prime,D_{Z^\prime,1}+N_{Z^\prime,1})$ over $S$, and a $K_{Z^\prime}+D_{Z^\prime,1}+N_{Z^\prime,1}$-trivial contraction $(Z^\prime,D_{Z^\prime,1}+N_{Z^\prime,1}) \to (X^\prime,D_1^\prime+N_1^\prime)$ over $S$. 

Since $(Z^\prime,D_{Z^\prime,e}+N_{Z^\prime,e}) \to (X^\prime,D_e^\prime+N_e^\prime)$ is $K_{Z^\prime}+D_{Z^\prime,e}+N_{Z^\prime,e}$-trivial, it suffices to show that $K_{Z^\prime}+D_{Z^\prime,e}+N_{Z^\prime,e}$ is nef. 
Suppose not, we may further run MMP on $(Z^\prime,D_{Z^\prime,e}+N_{Z^\prime,e})$ and end with a minimal model $(Z^{\prime\prime},D_{Z^{\prime\prime},e}+N_{Z^{\prime\prime},e})$. 

Since $Z^{\prime\prime}$ is a wlc model for $(X,D_e+N_e)$, it is also a wlc model for $(X,D_1+N_1)$. Also by our assumption on $B_{F+M}$, $X \dashrightarrow Z^{\prime\prime}$ is crepant with respect to $K_X+B_1+M_1$. This shows that $Z^\prime \dashrightarrow Z^{\prime\prime}$ is crepant for both $(Z^\prime, B_{Z^\prime,1}+M_{Z^\prime,1})$ and $(Z^\prime, D_{Z^\prime,1}+N_{Z^\prime,1})$. So $Z^\prime \dashrightarrow Z^{\prime\prime}$ is also crepant for $(Z^\prime, D_{Z^\prime,e}+N_{Z^\prime,e})$. But this contradicts the assumption that $K_{Z^\prime}+D_{Z^\prime,e}+N_{Z^\prime,e}$ is not nef. 

\hfill

Claim 2.1 and Claim 2.2 combine to give a desired decomposition in Conjecture A.2. 

\end{proof}

\subsection{Some easy lemmas}

Here we collect some easy lemmas, which are known to experts. 

\begin{lemma}\label{lem wlc equiv}
Let $(X,B+M)$ be a generalised pair, and $X \to X^\prime$ be a $K_X+B+M$-negative birational contraction. Let $X \dashrightarrow Z$ be a birational map. 

Then $Z$ is a wlc model of $(X,B+M)$ if and only if $Z$ is a wlc model of $(X^\prime,B^\prime+M^\prime)$. 
\end{lemma}

\begin{proof} 
Let $B_Z+M_Z$ be the log birational transform of $B+M$ on $Z$. 

Let $W$ be a common resolution of $X,X^\prime,Z$ as in the diagram. 

\[\begin{tikzcd}
	& W \\
	{X^\prime} & X & Z
	\arrow["q"', from=1-2, to=2-1]
	\arrow["p"', from=1-2, to=2-2]
	\arrow["r", from=1-2, to=2-3]
	\arrow["\phi"', dashed, from=2-2, to=2-1]
\end{tikzcd}\]

Firstly we show that $q^*(K_{X^\prime}+B^\prime+M^\prime) \geq r^*(K_Z+B_Z+M_Z)$. This inequality will be called $(*)$. 

Since $r^*(K_Z+B_Z+M_Z) - q^*(K_{X^\prime}+B^\prime+M^\prime)$ is nef over $X^\prime$, by negativity lemma, it suffices to show $K_{X^\prime}+B_1^\prime+M_1^\prime \geq q_*r^*(K_Z+B_Z+M_Z)$. By definition of wlc models, $p^*(K_{X}+B_1+M_1) \geq r^*(K_Z+B_Z+M_Z)$, and taking $q_*$ gives our desired inequality. 

Secondly we claim that $B_Z+M_Z$ is also the log birational transform of $B^\prime+M^\prime$ on $Z$. 
To prove this, it suffices to show that for any exceptional divisor $E$ of $X^\prime/X$ with $a(E,X,B+M)<a(E,X^\prime,B^\prime+M^\prime)$, the strict transform of $E$ on $W$ is contracted by $r$. 
Suppose not, then we may consider $D = r^*(K_Z+B_Z+M_Z)-q^*(K_{X^\prime}+B^\prime+M^\prime) $. 
Then we have $D\leq 0$ by $(*)$. But we have 
\begin{equation*}
\begin{split}
{\rm Coeff}_E(D) &= -a(E,Z,B_Z+M_Z)+a(E,X^\prime,B^\prime+M^\prime) \\
&= -a(E,X,B+M)+a(E,X^\prime,B^\prime+M^\prime)>0
\end{split}
\end{equation*}
This leads to a contradiction. 

Then we prove ``only if" part. This follows from $(*)$ and the fact that $B_Z+M_Z$ is also the log birational transform of $B^\prime+M^\prime$ on $Z$. 

Finally we prove ``if" part. 
We need to show $p^*(K_{X}+B_1+M_1) \geq r^*(K_Z+B_Z+M_Z)$. 

By the assumption that $X \to X^\prime$ is a $K_X+B+M$-negative birational contraction, we have $p^*(K_{X}+B_1+M_1) \geq q^*(K_{X^\prime}+B^\prime+M^\prime)$.
By definition of wlc models, we have $q^*(K_{X^\prime}+B^\prime+M^\prime) \geq r^*(K_Z+B_Z+M_Z)$. These two inequalities give $p^*(K_{X}+B_1+M_1) \geq r^*(K_Z+B_Z+M_Z)$. 

\end{proof}

\subsection{Stability of extremal rays}

In this subsection, we try to prove an analog of \cite[Proposition 3.2]{Bir11}, with essentially the same arguments. Such lemma is also established in \cite{HL22} for g-klt g-pairs. 

We work over $\mathbb{C}$. 

For an $\mathbb{R}$-b-divisor $D +M = \sum d_iD_i +\sum u_jM_j$ where the $D_i$ are the irreducible components of $D$, and $M_j$ belongs to a set of b-divisors, we define $||D|| := \max\{|d_i|, |u_j|\}$. 

\begin{lemma}\label{Lem Bir11 Prop 3.2}
Let $X \to S$ be a projective morphism of normal quasi-projective varieties, and $F = \sum_{i=1}^k D_i$ be a reduced divisor on $X$. Let $\{M_j\}_{1\leq j\leq l}$ be b-nef $\mathbb{Q}$-divisors on $X$. 
Let $B_{F+M} \subset \bigoplus_{i=1}^k [0,1]D_i \oplus \bigoplus_{j=1}^l [0,1]M_j$ be a rational polytope such that $(X,B+M)$ is g-lc for $B+M\in B_{F+M}$. 
Fix some $B_0 + M_0\in B_{F+M}$. 
Then, there are real numbers $\alpha, \delta > 0$, depending on 
$X/S$, $F$, $\{M_j\}_{1\leq j\leq l}$ and $(B_0+M_0,B_{F+M})$
such that 

(1) If $R$ is any extremal curve over $S$ and if we have $(K_X + B_0 + M_0) \cdot R > 0$, then we have
 $(K_X + B_0 + M_0) \cdot R > \alpha$;

(2) If  we have $B +M \in B_{F+M}$, $||B +M - B_0 - M_0|| < \delta$, then we have $B_0 +M_0 +(1+2\dim X/\alpha)(B +M - B_0 - M_0) \in B_{F+M}$; 

(3) If we have $B +M \in B_{F+M}$, $||B+M - B_0-M_0|| < \delta$ and $(K_X +B+M)\cdot R \leq 0$ for an extremal ray $R$ over $S$, then we have $(K_X +B_0+M_0)\cdot R \leq 0$; 

(4) Let $\{R_t\}_{t\in T}$ be a family of extremal rays of $\overline{NE}(X/Z)$. Define
$\mathcal{N}_T = \{B+M\in B_{F+M}|(K_X +B+M)\cdot R_t \geq 0 \text{ for any } t\in T\}$. Then we have $\mathcal{N}_T = \mathcal{N}_{T^\prime}$ for some finite subset $T^\prime$ of $T$. In particular, the set $\mathcal{N}_T$
is a rational polytope.  

\end{lemma}

\begin{proof}
(1) The set $B_{F+M}$ is rational polyhedral, so there exist $\mathbb{Q}$-b-divisors $B_{0,i}+M_{0,i}\in B_{F+M}$, $1\leq i\leq r$, and $a_1, \dots, a_r \in \mathbb{R}_{>0}$ such that we have $(K_X+B_0+M_0) = \sum_{i=1}^r a_i (K_X+B_{0,i}+M_{0,i})$, where $a_i$ and the coefficients of $B_{0,i}$ depend only on $(B_0+M_0,B_{F+M})$ as a pointed subset of $\bigoplus_{i=1}^k [0,1]D_i \oplus \bigoplus_{j=1}^l [0,1]M_j = [0,1]^{k+l}$. Then by the cone theorem, see \cite[Theorem 1.3 (2)]{HL21}, we have $(K_X+B_{0,i}+M_{0,i}) \cdot R \geq -2\dim X$. So the number $(K_X+B_0+M_0) \cdot R =\sum_{i=1}^r a_i (K_X+B_{0,i}+M_{0,i})\cdot R $ lies in a DCC set that depends only on $\dim X$, $(a_i)_{i\in I}$ and the Cartier indices of $K_X+B_{0,i}+M_{0,i}$. In particular, such $\alpha$ exists. 

(2) If $B_0+M_0$ lies in the interior of $B_F$, let $d(B_0+M_0, \partial B_{F+M})$ be the distance between $B_0+M$ and the boundary of $B_{F+M}$. Then  we have $d(B_0+M_0, \partial B_{F+M})>0$ by compactness of $\partial B_{F+M}$. We may choose $\delta = d(B_0+M_0, \partial B_{F+M})/(1+ 2\dim X/\alpha)$. 

If $B_0+M_0$ lies in the boundary of $B_{F+M}$, there exists a unique face $S$ of $B_{F+M}$ such that $B_0+M_0$ lies in the relative interior of $S$. Let $\mathcal{T}$ be the union of faces $T$ of $B_{F+M}$ that $S$ is not contained in $T$. Then we have $d(B_0+M_0, \mathcal{T})>0$ by compactness of $\partial B_F$. 

We claim that $\delta = d(B_0+M_0, \mathcal{T})/(1+ 2\dim X/\alpha)$ satisfies our conditions. Indeed, $B_{F+M}$ is defined by finitely many linear inequalities $f_i(D+N) \geq 0$, $1\leq i\leq N$, where each hyperplane $f_i(D+N)=0$ intersects $B_{F+M}$ at some face. Let $t$ be the maximal real number such that $B_0+M_0+t(B+M-B_0-M_0) \in B_F$. Then there exists some $f_i$ such that $f_i(B_0+M_0)>0$ but $f_i(B_0+M_0+t(B+M-B_0-M_0)) = 0$. Then $B_0+M_0+t(B+M-B_0-M_0)$ lies in the face of $B_{F+M}$ defined by $f_i(D+N) = 0$, so we have $||t(B +M -B_0 -M_0)||\geq d(B_0+M_0,\mathcal{T})$, which implies $t\geq d(B_0+M_0,\mathcal{T})/||B +M -B_0 -M_0|| > d(B_0,\mathcal{T})/\delta = 1+2\dim X/\alpha$. 

(3) Suppose that we have $(K_X+B_0+M_0)\cdot R>0$, then we have $(K_X+B_0+M_0)\cdot R>\alpha$ by (1). Let $t = 1+2\dim X /\alpha$. By (2) and the cone theorem, see \cite[Theorem 1.3 (2)]{HL21}, we have $(K_X+B_0 +M_0+t(B +M -B_0 -M_0))\cdot R\geq -2\dim X$. So, we may compute 
\begin{align*}
&(K_X+B+M)\cdot R \\
= &
\frac{2\dim X}{\alpha+ 2\dim X} (K_X+B_0+M_0)\cdot R + \frac{\alpha}{\alpha+ 2\dim X} (K_X+B_0+M_0+t(B + M-B_0 -M_0))\cdot R 
\\> & \frac{2\dim X}{\alpha+ 2\dim X} \alpha - \frac{\alpha}{\alpha+ 2\dim X} 2\dim X 
\\= & 0
\end{align*}

This contracts the assumption that $(K_X+B+M)\cdot R$ is non-positive. 

(4) We proceed by induction on the dimension of $B_{F+M}$. The result is trivial when $\dim B_{F+M} = 0$. 

The set $\mathcal{N}_T$ is clearly compact, so by (3) there exist $B_1+M_1,\dots,B_n+M_n \in \mathcal{N}_T$ and $\delta_1,\dots,\delta_n > 0$ such that $\mathcal{N}_T$ is covered by $\mathbb{B}_i = \{B+M\in  \mathcal{N}_T | \text{ } ||B+M-B_i-M_i|| < \delta_i\}$ and such that if $B+M \in \mathbb{B}_i$ and $(K_X +B+M)\cdot R_t \leq 0$ for an extremal ray $R$ over $S$, then $(K_X +B_{i}+M_{i})\cdot R = 0$. 

Let $T_i = \{t\in T| (K_X+B+M)\cdot R_t<0 \text{ for some } B+M\in \mathbb{B}_i\}$. Then $\mathcal{N}_T = \cap_i \mathcal{N}_{T_i} = \mathcal{N}_{\cup _iT_i}$. So it suffices to show $\mathcal{N}_{T_i}=\mathcal{N}_{T_i^\prime}$ for some finite subset $T_i^\prime$ of $T_i$ for each $1\leq i\leq n$. 

Let $B_F^k$, $1\leq k\leq m$ be the proper faces of $B_F$. By induction hypothesis, $\mathcal{N}_{T_i} \cap B_F^k = \mathcal{N}_{T_i^{k\prime}} \cap B_F^k$ for some finite subset $T_i^{k\prime}$ of $T_i$. 

Set $T_i^\prime$ to be $\cup_k T_i^{k\prime}$, which is a finite subset of $T_i$. We claim that we have $\mathcal{N}_{T_i} = \mathcal{N}_{T_i^\prime}$

Indeed, as $(K_X +B_{i}+M_{i})\cdot R_t = 0$ for each $t \in T_i$, $\mathcal{N}_{T_i}$ is the convex hull of $\{B_{i}+M_{i}\}$ and $\mathcal{N}_{T_i} \cap B_{F+M}^k$, $1\leq k\leq m$. 

Similarly, as $(K_X +B_{i}+M_{i})\cdot R_t = 0$ for each $t \in  T_i^\prime$, $\mathcal{N}_{T_i^\prime}$ is the convex hull of $\{B_{i}+M_{i}\}$ and $\mathcal{N}_{T_i^{k\prime}} \cap B_{F+M}^k$, $1\leq k\leq m$.

So we have $\mathcal{N}_{T_i}=\mathcal{N}_{T_i^\prime}$ since we know $\mathcal{N}_{T_i} \cap B_{F+M}^k = \mathcal{N}_{T_i^{k\prime}} \cap B_{F+M}^k$, and $T_i^\prime$ is a finite subset of $T_i$.

\end{proof}

\begin{corollary}\label{Cor Bir11 Prop 3.2}
Notations as Lemma \ref{Lem Bir11 Prop 3.2}. Assume further that $X$ is $\mathbb{Q}$-factorial and $K_X+B_0+M_0$ is nef. 
Then, there are real numbers $\alpha, \delta > 0$, depending on
$X/S$, $F$, $\{M_j\}_{1\leq j\leq l}$ and $(B_0+M_0,B_{F+M})$
such that Lemma \ref{Lem Bir11 Prop 3.2}(1)-(4) holds, and moreover 

(1) There exists nonnegative 
real numbers $a_1,\dots,a_r$, b-divisors $B_{0,i} +M_{0,i} \in B_{F+M}$, $1\leq i\leq r$, 
and $ m\in \mathbb{N}$ 
depending only on $X/S$, $F$, $\{M_j\}_{1\leq j\leq l}$ and $(B_0+M_0,B_{F+M})$ so that we have $\sum_{i=1}^r a_i = 1$, $B_0+M_0 = \sum_{i=1}^r a_i(B_{0,i}+M_{0,i})$, $a_1, \dots,a_r$ are linearly independent over $\mathbb{Q}$, and each $m(K_X + B_{0,i}+M_{0,i})$ is Cartier and nef; 

(2) Let $\phi: X \dashrightarrow X^\prime$ be a birational contraction, which is crepant with respect to $K_X+B_0+M_0$. Assume $m\phi_*(K_X + B_{0,i}+M_{0,i})$ is Cartier for $1\leq i\leq r$. Then 

a) if $R$ is any extremal curve on $X^\prime$ over $S$ and if $\phi_*(K_X + B_0 + M_0) \cdot R > 0$, then
 $ \phi_*(K_X + B_0 + M_0) \cdot R > \alpha$;

b) if we have $B+M \in B_{F+M}$, $||B+M - B_0-M_0|| < \delta$ and $\phi_*(K_X +B+M)\cdot R \leq 0$ for some extremal ray $R$ over $S$. 
Assume further that $\phi$ is $K_X+B+M$-non-positive. Then we have $\phi_*(K_X +B_0+M_0)\cdot R = 0$; 

c) if we have $B+M \in B_{F+M}$, $||B+M - B_0-M_0|| < \delta$, and $\psi: X^\prime \to Y$ is a $K_X+B+M$-trivial birational contraction.  
Assume further that $\phi$ is $K_X+B+M$-non-positive. Then $\psi\circ \phi$ is crepant with respect to $K_X+B_0+M_0$. 

(3) If we have $B+M \in B_{F+M}$, $||B+M - B_0-M_0|| < \delta$, and $\phi: X \dashrightarrow X^\prime$ is a sequence of $K_X+B+M$-negative divisorial contractions, $K_X+B+M$-flips or $K_X+B+M$-flops. Then $m\phi_*(K_X + B_{0,i}+M_{0,i})$ is Cartier for $1\leq i\leq r$ and $\phi$ is crepant with respect to $K_X+B_0+M_0$. 

\end{corollary}

\begin{proof}
(1) Since the set of $B+M \in B_{F+M}$ such that $K_X+B+M$ is nef is rational polyhedral by Lemma \ref{Lem Bir11 Prop 3.2} (4), there exists nonnegative real numbers $a_1,\dots,a_r$, and $ m\in \mathbb{N}$ so that $\sum_{i=1}^r a_i = 1$, $B = \sum_{i=1}^r a_iB_i$, and each $m(K_X + B_{0,i}+M_{0,i})$ is Cartier. 

We may further assume $a_1, \dots,a_r$ are linearly independent over $\mathbb{Q}$. Indeed, once we have a nontrivial linear relation $\sum_{i=1}^r \lambda_ia_i = 0$ with $\lambda_i\in \mathbb{Q}$. Then $P = \{(x_1,\dots, x_r)\in [0,1]^r|\sum_{i=1}^r a_i = 1, \sum_{i=1}^r \lambda_ia_i = 0\}$ is a rational polytope of dimension $r-2$ containing $a=(a_1,\dots,a_r)$. So there exists $n-1$ rational points $p_i = (p_{i1}, \dots p_{ir})\in P$, $1\leq i\leq r-1$, such that there exists nonnegative real numbers $b_1,\dots,b_{r-1}$ such that $\sum_{i=1}^{r-1} b_ip_i = a$. So we have $B_0+M_0 = \sum_{i=1}^{r-1} b_i(\sum_{j=1}^rp_{ij}(B_{0,j}+M_{0,j}))$. So we can decrease $r$ until there are no such linear relations. 

(2)
By linear independence, $\phi$ is crepant with respect to $K_X+B_{0,i}+M_{0,i}$ for $1\leq i\leq r$. 

a) By the cone theorem, see \cite[Theorem 1.3 (2)]{HL21}, we have $\phi_*(K_X+B_{0,i}+M_{0,i}) \cdot R \geq -2\dim X$. So $\phi_*(K_X+B_0+M_0) \cdot R = \sum_{i=1}^r a_i \phi_*(K_X+B_{0,i}+M_{0,i})\cdot R $ lies in a DCC set depending only on $\dim X$, $(a_i)_{i\in I}$ and $m$. In particular, such $\alpha$ exists. 

b) Note that by Lemma \ref{Lem Bir11 Prop 3.2} (2), there exists $\delta$ depending only on $X/S$, $F$, $\{M_j\}_{1\leq j\leq l}$ and $(B_0+M_0,B_{F+M})$ such that if we have $B+M \in B_{F+M}$, $||B + M - B_0 -M_0|| < \delta$, then we have $B_0 +M_0 +(1+2\dim X/\alpha)(B +M -B_0 -M_0) \in B_{F+M}$.

Suppose $\phi_*(K_X +B_0+M_0)\cdot R \neq 0$. Since $K_X+B_0+M_0$ is nef and $\phi$ is $K_X+B_0+M_0$-crepant, we have $\phi_*(K_X +B_0+M_0)\cdot R \geq 0$. By part a), we have $\phi_*(K_X +B_0+M_0)\cdot R > \alpha$. 

Let $t = 1+2\dim X/\alpha$. 

Since $\phi$ is $K_X+B+M$-non-positive and $K_X+B_0+M_0$-crepant, $\phi$ is $K_X+B_0+M_0+t(B+M-B_0-M_0)$-non-positive. In particular, we have 
$$\phi_*(K_X+B_0+M_0+t(B+M-B_0-M_0)) = K_{X^\prime} + \phi_*B_0+ M_{0,X^\prime}+ t\phi_*(B+M-B_0-M_0)$$
and $(X^\prime, \phi_*B_0+ M_{0,X^\prime}+ t\phi_*(B+M-B_0-M_0))$ is g-lc. By the cone theorem, see \cite[Theorem 1.3 (2)]{HL21}, we have $\phi_*(K_X+B_0+M_0+t(B+M-B_0-M_0))\cdot R\geq -2\dim X$. So 
\begin{align*}
&\phi_*(K_X+B+M)\cdot R \\
= &
\frac{2\dim X}{\alpha+ 2\dim X} \phi_*(K_X+B_0+M_0)\cdot R + \frac{\alpha}{\alpha+ 2\dim X} \phi_*(K_X+B_0+M_0+t(B+M-B_0-M_0))\cdot R
\\> & \frac{2\dim X}{\alpha+ 2\dim X} \alpha - \frac{\alpha}{\alpha+ 2\dim X} 2\dim X 
\\= & 0
\end{align*}

This contracts the assumption that $\phi_*(K_X+B+M)\cdot R$ is non-positive. 

c) 
Suppose that $\psi \circ \phi$ is not crepant with respect to $K_X+B_0+M_0$. 
Then there exists some curve $R$ on $X^\prime$ such that we have $\psi(R) = 0$ and $\phi_*(K_X+B_0+M_0) \cdot R \neq 0$. Since $\phi_*(K_X+B_0+M_0)$ is nef, we have $\phi_*(K_X+B_0+M_0) \cdot R > 0$. 

Let $B_\mu+M_\mu = B_0+M_0 + (1+\mu)(B+M-B_0+M_0)$, then we have $B_\mu+M_\mu \in B_{F+M}$ for $\mu>0$ sufficiently small by Lemma \ref{Lem Bir11 Prop 3.2}(2). 

Moreover, we have $\phi_*(K_X+B_\mu+M_\mu) \cdot R < 0$. So $\phi_*(K_X+B_\mu+M_\mu)$ is not nef over $Y$. By the cone and contraction theorem \cite[Theorem 1.3(4)]{HL21}, there exists a $\phi_*(K_X+B_\mu+M_\mu)$-negative extremal ray $R$ on $X^\prime$ over $Y$. Then $R$ is also a $\phi_*(K_X+B_\mu+M_\mu)$-negative extremal ray on $X^\prime$ over $S$. 

By b), we deduce that $\phi_*(K_X+B_0+M_0) \cdot R = 0$. By assumption, We have $\phi_*(K_X+B+M)\cdot R=0$. So we have $\phi_*(K_X+B_\mu+M_\mu) \cdot R = 0 $, which contradicts to the assumption that $\phi_*(K_X+B_\mu+M_\mu) \cdot R < 0 $. 

(3) We proceed by induction on the length of the sequence. When the length is $0$, the result is trivial. 

Let $\phi: X \dashrightarrow X^\prime$ be factored as $f\circ \phi$, where $\psi: X\dashrightarrow Y$ is a sequence  of smaller length and $f: Y \dashrightarrow X^\prime$ is a divisorial contraction, a flip or a flop. Let $R$ be the extremal ray corresponding to the divisorial contraction, the flipping contraction or the flopping contraction. 

By induction hypothesis, $m\psi_*(K_X + B_{0,i}+M_{0,i})$ is Cartier for $1\leq i\leq r$ and $\psi$ is crepant with respect to $K_X+B_0+M_0$. 

By (2) b) and (2) c), we hace $\phi_*(K_X +B_0+M_0)\cdot R = 0$. By linear independence of $a_i$, $\phi_*(K_X +B_{0,i}+M_{0,i})\cdot R = 0$ for $1\leq i\leq r$. 

If $f$ is a divisorial contraction or a flip, by the base point free theorem, see \cite[Theorem 1.5]{Xie24}, $mf_*\phi_*(K_X+B_{0,i}+M_{0,i})$ is Cartier for $1\leq i\leq r$ and $f$ is crepant with respect to $\psi_*(K_X+B_0+M_0)$. 

If $f$ is a flop, let $H$ be an ample divisor on $X^\prime$. We may view $H$ as a b-divisor, then $H$ is a b-nef b-divisor. So $f^*H$ is a b-nef b-divisor on $Y$. Applying \cite[Theorem 1.5]{Xie24} for the generalised pair $(Y,\phi_*(B_0+M_0)+f^*H)$, we have $mf_*\phi_*(K_X+B_{0,i}+M_{0,i})$ is Cartier for $1\leq i\leq r$ and $f$ is crepant with respect to $\psi_*(K_X+B_0+M_0)$.

\end{proof}

\subsection{Realization of weak log canonical models by MMP}

The following proposition might be considered as an analog of \cite[Proposition 6.2]{HH19}. We restrict to $\mathbb{Q}$-factorial case
to avoid some technical issue. 

We work over $\mathbb{C}$. 

\begin{proposition}\label{Prop special termination}
Let $(X/S,B+M)$ be an $\mathbb{Q}$-factorial g-dlt g-pair. Suppose $(X/S,B+M)$ admits a wlc model over $S$. 

Then there exists a sequence of MMP, which terminates with a minimal model. 
\end{proposition}

\begin{proof}
Let $\mathbb{K}$ be the field generated by $\mathbb{Q}$ and the coefficients of $B+M$.
Let $\rho = \rho(X/Z)$ be the relative Picard number of $X/S$. Let $\alpha_1,\dots,\alpha_\rho$ be real numbers which are linearly independent over $\mathbb{K}$, and $A_1,\dots,A_\rho$ be sufficiently ample divisors that generate $N^1(X/Z)$ as an $\mathbb{R}$-vector space. Put $A = \sum_i \alpha_iA_i$. Since $A_i$ are sufficiently ample, we may assume $K_X +B +A+M$ is nef over $S$. 

We run MMP with scaling of $A$, in the sense of \cite[Definition 3.20]{HL22}. 
$$(X,B+M) = (X_1,B_1+M_1) \dashrightarrow (X_2,B_2+M_2) \dashrightarrow \cdots \dashrightarrow (X_i,B_{i}+M_{i}) \dashrightarrow$$
and set $$\lambda_i = \inf\{\mu \in \mathbb{R}_{\geq 0} | K_{X_i} + B_i + \mu A_i + M_i \text{ is nef over } Z\}$$ where $A_i$
is the birational
transform of $A$ on $X_i$. 

By the definition of MMP with scaling, see \cite[Definition 3.20]{HL22}, in the $i$-th step we choose a $K_{X_i}+B_i+M_i$-negative extremal ray $R_i$ on $X_i$, such that $(K_{X_i}+B_i+\lambda_iA_i)\cdot R_i = 0$. Such an extremal ray exists by \cite[Lemma 3.19]{HL22}. Then we do a divisorial contraction or a flip with respect to $R_i$ to get $X_{i+1}$. In this case, $K_{X_{i+1}}+B_{i+1}+\lambda_i A_{i+1}$ is nef, so $\lambda_{i+1}\leq \lambda_i$. 

We claim that this MMP terminates. If not, since the relative Picard number is non-increasing, by taking a subsequence, we may assume this MMP consists of flips. 

As in \cite[Proposition 6.2]{HH19}, we claim that $\lambda_i$ are strictly decreasing. Otherwise, suppose $\lambda_i = \lambda_{i+1}$, and consider the following diagram: 
\[\begin{tikzcd}
	{X_i} && {X_{i+1}} && {X_{i+2}} \\
	& {V_i} && {V_{i+1}}
	\arrow[dashed, from=1-1, to=1-3]
	\arrow["{f_i}"', from=1-1, to=2-2]
	\arrow[dashed, from=1-3, to=1-5]
	\arrow["{f_i^+}", from=1-3, to=2-2]
	\arrow["{f_{i+1}}"', from=1-3, to=2-4]
	\arrow["{f_{i+1}^+}", from=1-5, to=2-4]
\end{tikzcd}\]
where $f_i$, $f_{i+1}$ are flipping contractions, and $f_i^+$, $f_{i+1}^+$ are the corresponding flips. 

Let $C$ be a curve contracted by $f_i^+$, and $C^\prime$ be a curve contracted by $f_{i+1}$. Then by the definition of MMP with scaling, $$C \cdot (K_{X_{i+1}}+B_{i+1}+\lambda_iA_{i+1}+M_{i+1}) = C^\prime\cdot (K_{X_{i+1}}+B_{i+1}+\lambda_iA_{i+1}+M_{i+1}) = 0$$
We may assume $\lambda_i>0$, otherwise the MMP already terminates. 
So we have $$C\cdot A_{i+1} = \lambda^{-1}_iC\cdot (K_{X_{i+1}}+B_{i+1}+M_{i+1}), \text{ and } C^\prime\cdot A_{i+1} = \lambda^{-1}_iC^\prime\cdot (K_{X_{i+1}}+B_{i+1}+M_{i+1})$$ which shows that there exists a nonzero number $q\in \mathbb{K}$ such that $C\cdot A_{i+1} = qC^\prime\cdot A_{i+1}$. But since $A= \sum_i \alpha_iA_i$ with $A_i$ generate $N^1(X)$ and $\alpha_i$ linearly independent over $\mathbb{K}$, we deduce that $C \equiv qC^\prime$. But this forces $f_i^+ = f_{i+1}$, which is a contradiction. This shows that $\lambda_i$ are strictly decreasing.

Let $\mu = \lim_{i\to \infty} \mu_i$. If $\mu>0$, $(X/Z,B+\mu A+M)$ admits a minimal model by \cite[Lemma 4.4]{BZ16}. If $\mu=0$, by \cite[Proposition 3.9]{HL22}, we can take a $\mathbb{Q}$-factorial dlt modification of the wlc model, which becomes a g-log minimal model in the sense of \cite[Definition 2.9]{HL22}. 

By \cite[Theorem 4.1]{HL22}, this MMP terminates.

\end{proof}

The following lemma should be useful. We note that this lemma might be considered as a generalization of the fact that flops connect minimal models, which is proved by \cite{Kaw07}. Such a decomposition for usual pairs appears in \cite[Theorem 4.1]{Bir12}. 

\begin{lemma}\label{Lem realization of wlc models}
Let $(X/S,B+M)$ be a NQC g-lc g-pair, and $(X^\prime,B^\prime+M^\prime)$ be a wlc model of $(X,B+M)$ over $S$. 

Then there exists a $\mathbb{Q}$-factorial dlt modification $(Y,B_Y+M_Y)$ of $(X,B+M)$, a sequence of $K_Y+B_Y+M_Y$-MMP $(Y,B_Y+M_Y) \dashrightarrow (Z,B_Z+M_Z)$ over $S$, a sequence of $K_Z+B_Z+M_Z$-flops $(Z,B_Z+M_Z) \dashrightarrow (Z^\prime,B_{Z^\prime}+M_{Z^\prime})$ over $S$, and a $K_{Z^\prime}+B_{Z^\prime}+M_{Z^\prime}$-trivial contraction $(Z^\prime,B_{Z^\prime}+M_{Z^\prime}) \to (X^\prime,B^\prime+M^\prime)$ over $S$. 

Moreover, we may assume the exceptional divisors of $Y/X$ have log discrepancy $0$ with respect to $(X,B+M)$. 
\end{lemma}

\begin{proof}
Let $(W,B_W+M_W)$ be a common log resolution of $(X,B+M)$ and $(X^\prime,B^\prime+M^\prime)$. 

Let $B_W^{log}$ be the log birational transform of $B$, i.e. $B_W^{log}$ is the sum of the strict transform of $B$ and the reduced exceptional divisor. 
We may run $K_W+B_W^{log}+M_W$-MMP on $W$ over $X$, which terminates with a minimal model $(Y,B_Y+M_Y)$ by \cite[Proposition 3.8]{HL22}. 

Similarly, we may run $K_W+B_W^{log}+M_W$-MMP on $W$ over $X^\prime$, which terminates with a minimal model $(Z^\prime,B_{Z^\prime}+M_{Z^\prime})$ by \cite[Proposition 3.8]{HL22}. Then $Z^\prime \to X^\prime$ is a $K_{Z^\prime}+B_{Z^\prime}+M_{Z^\prime}$-trivial contraction over $S$ by negativity lemma. 

Then $\phi: (Y,B_Y+M_Y) \dashrightarrow (X^\prime,B^\prime+M^\prime)$ is a birational contraction and $(X^\prime,B^\prime+M^\prime)$ is a wlc model of $(Y,B_Y+M_Y)$. 
We may run $K_Y+B_Y+M_Y$-MMP on $Y$ over $S$, which terminates with a minimal model $(Z,B_Z+M_Z)$ by Proposition \ref{Prop special termination}. 

We claim that the natural map $\mu: (Z,B_Z+M_Z) \dashrightarrow (Z^\prime,B_{Z^\prime}+M_{Z^\prime})$ is small. Indeed, the contracted divisors of $p_{Z^\prime}: W \dashrightarrow Z^\prime$ are precisely the support of $K_W+B_W^{log}+M_W - p_{Z^\prime}^*(K_{Z^\prime}+B_{Z^\prime}+M_{Z^\prime})$. On the other hand, the contracted divisors of $p_Y: W \dashrightarrow Y$ are precisely the support of $K_W+B_W^{log}+M_W - p_{Y}^*(K_{Y}+B_{Y}+M_{Y})$, and the contracted divisors of $q: Y \dashrightarrow Z$ are precisely the support of $K_Y+B_Y+M_Y - p_{Z}^*(K_{Z}+B_{Z}+M_{Z})$. So the contracted divisors of $p_{Z}: W \dashrightarrow Z$ are precisely the support of $K_W+B_W^{log}+M_W - p_{Z}^*(K_{Z}+B_{Z}+M_{Z})$. Thus, our claim holds. 

It remains to show that the natural map $\mu: (Z,B_Z+M_Z) \dashrightarrow (Z^\prime,B_{Z^\prime}+M_{Z^\prime})$ can be decomposed into a sequence of $K_Z+B_Z+M_Z$-flops over $S$. 

By \cite[Lemma 5]{Cha22}, there exists an open subset $U$ of $Z$ such that $\mu|_U$ is an isomorphism and all glc centers of $(Z,B_Z+M_Z)$ intersect $U$. 

Then by \cite[Theorem 1]{Cha22}, can be decomposed into a sequence of $K_Z+B_Z+M_Z$-flops over $S$. Note that we are working on $\mathbb{Q}$-factorial case, so small modifications and switching of $Z$ and $Z^\prime$ are not needed.

\end{proof}

\subsection{Stability of weak log canonical models}

\begin{proposition}\label{Prop stability of wlc models}
Let $X \to S$ be a projective morphism of normal quasi-projective varieties, and $F = \sum_{i=1}^k D_i$ be a reduced divisor on $X$. Let $\{M_j\}_{1\leq j\leq l}$ be b-nef $\mathbb{Q}$-divisors on $X$. 
Let $B_{F+M} \subset \bigoplus_{i=1}^k [0,1]D_i \oplus \bigoplus_{j=1}^l [0,1]M_j$ be a rational polytope such that $(X,B+M)$ is g-lc for $B+M\in B_{F+M}$. 
Fix some $B_0 + M_0\in B_{F+M}$. 
Assume further that $X$ is $\mathbb{Q}$-factorial and $K_X+B_0+M_0$ is nef. 

Then, there is a real number $\delta > 0$, depending on 
$X/S$, $F$, $\{M_j\}_{1\leq j\leq l}$ and $(B_0+M_0,B_{F+M})$
such that: 

If $B+M \in B_{F+M}$, $||B+M - B_0-M_0|| < \delta$ and $(X^\prime, B^\prime+M^\prime)$ is a wlc model of $(X, B+M)$ over $S$, then the rational map $(X,B_0+M_0) \dashrightarrow (X^\prime,B_0^{\prime}+M_0^{\prime})$ is crepant, where $B_0^\prime+M_0^\prime$ is the log birational transform of $B_0+M_0$. 
\end{proposition}

\begin{proof}
We choose $\delta$ as in Lemma \ref{Lem Bir11 Prop 3.2} and Corollary \ref{Cor Bir11 Prop 3.2}. 

By Lemma \ref{Lem realization of wlc models}, there exists a $\mathbb{Q}$-factorial dlt modifacation $(Y,B_Y+M_Y)$ of $(X,B+M)$, a sequence of $K_Y+B_Y+M_Y$-MMP $(Y,B_Y+M_Y) \dashrightarrow (Z,B_Z+M_Z)$ over $S$, a sequence of $K_Z+B_Z+M_Z$-flops $(Z,B_Z+M_Z) \dashrightarrow (Z^\prime,B_{Z^\prime}+M_{Z^\prime})$ over $S$, and a $K_{Z^\prime}+B_{Z^\prime}+M_{Z^\prime}$-trivial contraction $(Z^\prime,B_{Z^\prime}+M_{Z^\prime}) \to (X^\prime,B^\prime+M^\prime)$ over $S$. 
Moreover, the exceptional divisors of $Y/X$ have log discrepancy $0$ with respect to $(X,B+M)$. 

Let $E_1,\dots, E_n$ be the exceptional divisors of $Y/X$. We will show that $B+M$ and $B_0+M_0$ are still close after taking crepant pullback to $Y$, and then we may apply Corollary \ref{Cor Bir11 Prop 3.2}. 

Since $B+M \in B_{F+M}$, $||B+M - B_0-M_0|| < \delta$, the generalised pair $(X,B+M+\epsilon(B+M-B_0-M_0))$ is g-lc for some small $\epsilon>0$ by Lemma \ref{Lem Bir11 Prop 3.2}(2). But for any exceptional divisor over $X$, $a(X,B+M,E)$ is a nontrivial convex combination of $a(X,B_0+M_0,E)$ and $a(X,B+M+\epsilon(B+M-B_0-M_0),E)$. This means that the exceptional divisors of $Y/X$ have log discrepancy $0$ with respect to $(X,B_0+M_0)$. 
In particular, we have $||B_Y+M_Y - B_{0,Y}-M_{0,Y}|| < \delta$. 

Let $F_Y$ be the sum of the strict transform of $F$ and the reduced exceptional divisor, say $F_Y = \sum_{i=1}^kD_{Y,i}+ \sum_{u=1}^nE_u$, and $M_{Y,j}$ be the $b$-divisor on $Y$ corresponding to $M_j$. Consider the affine map $$\bigoplus_{i=1}^k [0,1]D_i \oplus \bigoplus_{j=1}^l [0,1]M_j \to \bigoplus_{i=1}^k [0,1]D_{Y,i} \oplus \bigoplus_{u=1}^n [0,1]E_u\oplus \bigoplus_{j=1}^l [0,1]M_{Y,j}$$
$$\sum_{i=1}^k d_iD_i +\sum_{j=1}^l m_jM_j \mapsto \sum_{i=1}^k d_iD_{Y,i} + \sum_{u=1}^n E_u+\sum_{j=1}^l m_jM_{Y,j}$$
and denote the image of $B_{F+M}$ by $B_{F_Y+M_Y}$. 

Let $K_Y+B_{Y,0}+M_{Y,0}$ be the crepant pullback of $K_X+B_0+M_0$. One can check directly that the same $\alpha,\delta$ as $(X/Z, F, \{M_{j}\}_{1\leq j\leq l}, (B_{0}+M_{0},B_{F+M}))$ works for $(Y/Z, F_Y, \{M_{Y,j}\}_{1\leq j\leq l}, (B_{Y,0}+M_{Y,0},B_{F_Y+M_Y}))$ in Lemma \ref{Lem Bir11 Prop 3.2} and Corollary \ref{Cor Bir11 Prop 3.2}. 

Then by Corollary \ref{Cor Bir11 Prop 3.2}(2)c) and Corollary \ref{Cor Bir11 Prop 3.2}(3), the birational map $(Y,B_{Y,0}+M_{Y,0}) \dashrightarrow (X^\prime,B_0^\prime+M_0^\prime)$ is crepant. So the rational map $(X,B_0+M_0) \dashrightarrow (X^\prime,B_0^{\prime}+M_0^{\prime})$ is crepant.

\end{proof}

\subsection{Galois descent}

In this subsection, we work on an arbitrary field $k$ of characteristic $0$. 

\begin{conjecture}\label{Conj main 2}
Let $X \to S$ be a projective morphism of normal quasi-projective varieties, and $F = \sum_i D_i$ be a reduced divisor on $X$. Let $\{M_j\}_{j\in J}$ be b-nef $\mathbb{Q}$-divisors on $X$. 
Let $B_{F+M} \subset \bigoplus_i [0,1]D_i \oplus \bigoplus_j [0,1]M_j$ be a rational polytope and $N_{F+M}$ be the subset of b-divisors $B+M$ in $B_{F+M}$ such that $(X,B+M)$ is g-lc and $K_X+B+M$ is pseudo-effective. Then: 

(1) $N_{F+M}$ is closed convex rational polyhedral. 

(2) There exists a finite decomposition $N_{F+M} = \bigsqcup_i \Pi_i^\circ$ such that each $\Pi_i$ is closed convex rational polyhedral and any divisors $D_1+N_1$ and $D_2+N_2$ in $\Pi_i^\circ$ admit the same wlc models. 

\end{conjecture}

\begin{theorem}\label{Thm Galois descent}
Conjecture \ref{Conj main} in dimension $d$ implies Conjecture \ref{Conj main 2} in dimension $d$. 

More precisely, we may fix the data $(X/S,F,\{M_j\}_{j\in J},B_{F+M})$ in Conjecture \ref{Conj main 2}. Then there exist some data $(X_{\mathbb{C}}/S_{\mathbb{C}},F_{\mathbb{C}},\{M_{{\mathbb{C}},j}\}_{j\in J},B_{{\mathbb{C}},F+M})$ such that: 

(1) There exist a subfield $l \subset k$, some data $(X_l/S_{l},F_{l},\{M_{{l},j}\}_{j\in J},B_{{l},F+M})$ over $l$, and an embedding $l \to \mathbb{C}$ such that: 

The data
$(X/S,F,\{M_j\}_{j\in J},B_{F+M})$ and $(X_{\mathbb{C}}/S_{\mathbb{C}},F_{\mathbb{C}},\{M_{{\mathbb{C}},j}\}_{j\in J},B_{{\mathbb{C}},F+M})$ are given by taking base change from $(X_l/S_{l},F_{l},\{M_{{l},j}\}_{j\in J},B_{{l},F+M})$ by the inclusions $l \subset k$ and $l \to \mathbb{C}$. 

(2) Assume \ref{Conj main} holds for the data $(X_{\mathbb{C}}/S_{\mathbb{C}},F_{\mathbb{C}},\{M_{{\mathbb{C}},j}\}_{j\in J},B_{{\mathbb{C}},F+M})$. Then Conjecture \ref{Conj main 2} holds for $(X/S,F,\{M_j\}_{j\in J},B_{F+M})$. 
\end{theorem}

\begin{proof}
Consider the set $S$ of $b$-divisors $B+M \in B_{F+M}$ such that $(X,B+M)$ is generalised lc. 
For $B+M \in S$, let $W_{B+M}$ be the set of wlc models of $(X,B+M)$. 

By Proposition \ref{Prop countable}, for fixed $B+M \in S$, the set $W_{B+M}$ is countable. Since the cardinality of $S$ is at most the cardinality of $\mathbb{C}$, the cardinality of set $\bigsqcup_{B+M\in S} W_{B+M}$ is at most the cardinality of $\mathbb{C}$. 

So there exists a subfield $l\subset k$ generated by a set of cardinality at most the cardinality of $\mathbb{C}$ such that all the birational maps $(X,B+M) \dashrightarrow Y$ that exhibit $Y$ as a wlc model of $(X,B+M)$ for some $B+M\in S$ are defined on $l$. 

Choose an embedding $l \to \mathbb{C}$. By Conjecture \ref{Conj main}, $N_{F_{\mathbb{C}}+M_{\mathbb{C}}}$ is closed rational polyhedral, and so is $N_{F+M}$, since g-lc property and pseudo-effectivity are preserved under change of base fields. 

Moreover, we have a finite decomposition $N_{F_{\mathbb{C}}+M_{\mathbb{C}}} = \bigsqcup_i (\Pi_{\mathbb{C}})_i^\circ$ such that each $(\Pi_{\mathbb{C}})_i$ is closed convex rational polyhedral and any divisors $D_1+N_1$ and $D_2+N_2$ in $(\Pi_{\mathbb{C}})_i^\circ$ admit the same wlc models. 

We claim that the corresponding decomposition $N_{F+M} = \bigsqcup_i \Pi_i^\circ$ satisfies the condition that any divisors $D_1+N_1$ and $D_2+N_2$ in $\Pi_i^\circ$ admit the same wlc models. 

Let $\phi: (X,D_1+N_1) \dashrightarrow Y$ be a wlc model for $(X,D_1+N_1)$, then by our choice of $l$, $\phi: (X,D_1+N_1) \dashrightarrow Y$ is defined on $l$, say $\phi = (\phi_l) \times_l k$ with $\phi_l: (X_l,(D_l)_1+(N_l)_1) \dashrightarrow Y_l$. 

Let $\phi_{\mathbb{C}} = (\phi_l) \times_l \mathbb{C}$ with $\phi_{\mathbb{C}}: (X_{\mathbb{C}},(D_{\mathbb{C}})_1+(N_{\mathbb{C}})_1) \dashrightarrow Y_{\mathbb{C}}$. 

By our choice of $\Pi_i^\circ$, $\phi_{\mathbb{C}}$ is a wlc model for $(X_{\mathbb{C}},(D_{\mathbb{C}})_2+(N_{\mathbb{C}})_2)$. Since being a wlc model is invariant under base change, $\phi$ is a wlc model for $(X,D_2+N_2)$. 

We may switch $D_1+N_1$ and $D_2+N_2$, and we deduce that any divisors $D_1+N_1$ and $D_2+N_2$ in $(\Pi_{\mathbb{C}})_i^\circ$ admit the same wlc models. 

\end{proof}

\subsection{Countability of wlc models}

In this subsection, we work on an arbitrary field $k$ of characteristic $0$. 

\begin{lemma}\label{Lem countable2}
Let $X$ be a normal $\mathbb{Q}$-factorial variety, then the set of birational contractions $\phi: X \dashrightarrow Y$ is countable. 
\end{lemma}

\begin{proof}
Let $W$ be the normalization of the closure of the graph of $\phi$ with natural maps $p: W\to X$ and $q: W\to Y$. Let $A$ be an ample divisor on $Y$, $A_X = p_*q^*A$, and $A_W = q^*A$. 

We claim that $\phi$ coincides with the natural rational map $$X \dashrightarrow \text{Proj}\oplus_{n=0}^{\infty}H^0(X,\mathcal{O}_X(nA_X))$$

To see this, it suffices to show that the natural map $ H^0(Y,\mathcal{O}_Y(nA_Y)) \to H^0(X,\mathcal{O}_X(nA_X))$ is an isomorphism for sufficiently divisible $n$. Injectivity is clear, so it suffices to show surjectivity. 

Let $f_X$ be a rational function on $X$ such that $(f_X)+nA_X \geq 0$. Let $f_W$, $f_Y$ be the corresponding rational function on $W$, $Y$. To show surjectivity, it suffices to show that $(f_W)+nA_W \geq 0$. Since $A_W = q^*A$, $(f_W)+nA_W \geq 0$ if and only if $(f_Y)+nA_Y \geq 0$. But we indeed have $(f_Y)+nA_Y \geq 0$ since $(f_Y)+nA_Y = \phi_*((f_X)+nA_X)$. 

Then countability of $\phi$ follows from countability of $\mathbb{Q}$-divisor classes on $X$.

\end{proof}

\begin{lemma}\label{Lem countable3}
Let $(X,B)$ be a log smooth lc pair over an algebraically closed field of characteristic $0$, then there are only countably many divisors over $(X,B)$ whose log discrepancy is less than $1$. 
\end{lemma}

\begin{remark}
The result is not true if we allow log discrepancy $1$. 
\end{remark}

\begin{proof}
Let $E$ be a divisor over $(X,B)$ with log discrepancy less than $1$. 
Then there exists $n>0$ such that $E$ be a divisor over $(X,(1-1/n)B)$ with log discrepancy less than $1$. 
But $(X,(1-1/n)B)$ is klt, there are finitely many divisors over $(X,(1-1/n)B)$ with log discrepancy less than $1$ for fixed $n$. 

Taking union with respect to all $n\in \mathbb{N}_{>0}$, we deduce that there are only countably many divisors over $(X,B)$ whose log discrepancy is less than $1$.

\end{proof}

\begin{lemma}\label{Lem countable1}
Let $(X,B+M)$ be a g-lc pair over an algebraically closed field over a field of characteristic $0$, then the set of $\mathbb{Q}$-factorial dlt modifications $\phi: (Y,B_Y+M_Y) \to (X,B+M)$ of $(X,B+M)$ is countable. 
\end{lemma}

\begin{proof}
By Lemma \ref{Lem countable2}, it suffices to show that there are only countably many divisors over $(X,B+M)$ whose log discrepancy is at most $1$. 

Let $f: W \to (X,B+M)$ be a log resolution where $M$ descends, and $K_W+E+M_W = f^*(K_X+B+M)+F$, where $E$ and $F$ are effective without common components. By Lemma \ref{Lem countable3}, the set of divisors over $(W,E)$ of log discrepancy less than $1$ is countable, so the set of divisors over $(X,B+M)$ of log discrepancy less than $1$ is also countable.

\end{proof}

\begin{lemma}\label{Lem countable inf}
Let $(X,B+M)$ be a g-lc pair over an algebraically closed field of characteristic $0$, then the set of wlc models $\psi:  (X,B+M) \dashrightarrow (Y,B_Y+M_Y)$ of $(X,B+M)$ is countable. 
\end{lemma}

\begin{proof}
We may choose a $\mathbb{Q}$-factorial dlt modification $(W,B_W+M_W)$ of $(X,B+M)$ such that the natural map $W \dashrightarrow Y$ is a birational contraction. Then the result follows from Lemma \ref{Lem countable2} and Lemma \ref{Lem countable1}.

\end{proof}

\begin{proposition}\label{Prop countable}
Let $(X,B+M)$ be a g-lc pair over field $k$ of characteristic $0$, then the set of wlc models $\psi:  (X,B+M) \dashrightarrow (Y,B_Y+M_Y)$ of $(X,B+M)$ is countable. 
\end{proposition}

\begin{proof}
Observe that if $\psi: (X,B+M) \dashrightarrow (Y,B_Y+M_Y)$ is a wlc model of $(X,B+M)$, then the base change to $\bar{k}$, say $\psi_{\bar{k}}: (X_{\bar{k}},B_{\bar{k}}+M_{\bar{k}}) \dashrightarrow (Y_{\bar{k}},(B_Y)_{\bar{k}}+(M_Y)_{\bar{k}})$, is a wlc model of $(X_{\bar{k}},B_{\bar{k}}+M_{\bar{k}})$. 

Moreover, if $(\psi_1)_{\bar{k}}$ coincide with $(\psi_2)_{\bar{k}}$, then we must have $\psi_1$ coincide with $\psi_2$. 

So the result follows from Lemma \ref{Lem countable inf}.

\end{proof}

\vspace{1em}
 

\noindent\small{\textsc{Qiuzhen College, Tsinghua University, Beijing, China} }

\noindent\small{Email: \texttt{xufl23@mails.tsinghua.edu.cn}}



\begin{thebibliography}{100}


\bibitem[Bir09]{Bir09} C. Birkar, \emph{On existence of log minimal models}, Compositio Math., \textbf{145} (2009), 1442--1446.

\bibitem[Bir11]{Bir11} C. Birkar, \emph{On existence of log minimal models {II}}, Journal f\"ur die reine und angewandte
Mathematik, \textbf{2011} (2011), no. 658, 99--113.



\bibitem[Bir12]{Bir12} C. Birkar, \emph{Existence of log canonical flips and a special LMMP}, Publ. Math. de l'IH\'ES, 
\textbf{115} (2012), no. 1, 325--368.

\bibitem[BZ16]{BZ16} C. Birkar and D-Q. Zhang, \emph{Effectivity of Iitaka fibrations and pluricanonical systems of
polarized pairs}, Publications math\'ematiques de l’IH\'ES 123 (2016), 283--331.

\bibitem[BB05]{BB05} A. Bj\"orner and F. Brenti, \emph{Combinatorics of Coxeter groups}, volume 231 of Graduate Texts
in Mathematics. Springer-Verlag, 2005.

\bibitem[BLvL20]{BLvL19}
M. Bright, A. Logan, and R. van Luijk, \emph{Finiteness results for K3 surfaces over arbitrary fields},
Eur. J. Math., \textbf{6} (2020), no. 2, 336--366. 


\bibitem[Cha22]{Cha22} Priyankur Chaudhuri, \emph{Flops and minimal models for generalized pairs}, arXiv e-prints (2022),
arXiv:2211.11230.

\bibitem[CHL25]{CHL25} G. Chen, J. Han, and W. Liu, \emph{On the Iitaka volumes of log canonical surfaces and threefolds},
Journal of the London Mathematical Society, \textbf{111} (2025), no. 4, e70132. 

\bibitem[CT23]{CT23} G. Chen, N. Tsakanikas, \emph{On the termination of flips for log canonical generalized pairs}, Acta. Math. Sin.-English Ser. \textbf{39} (2023), 967--994. 


\bibitem[Choi08]{Choi08} S. R. Choi, \emph{Geography of log models and its applications}, Thesis. Johns Hopkins University, 2008.



\bibitem[FAG05]{FAG}
A. Vistoli et al., \emph{Grothendieck's FGA Explained}, Mathematical Surveys and
Monographs 123, American Mathematical Society, 2005.

\bibitem[FHS24]{FHS21}S. Filipazzi, C. D. Hacon, and R. Svaldi, \emph{Boundedness of elliptic Calabi-Yau
threefolds}, J. Eur. Math. Soc. (2024), published online first. 


\bibitem[Fu11]{Fu11}
L. Fu, \emph{\'Etale cohomology theory}, Nankai Tracts in Mathematics, 13.
World Scientific Publishing Co. Pte. Ltd., Hackensack, NJ, 2011.

\bibitem[Fuk02]{Fuk02}
S. Fukuda, \emph{On numerically effective log canonical divisors}, Int. J. Math.
Math. Sci., \textbf{30} (2002), no. 9, 521--531. 

\bibitem[Gac25]{Gac25} C. Gachet, \emph{Well-clipped cones under finite quotients and applications to the cone conjecture}, arXiv e-prints (2025), arXiv:2504.01753. 

\bibitem[GLSW24]{GLSW24} C. Gachet, H. Lin, I. Stenger, and L. Wang, \emph{The effective cone conjecture for Calabi--Yau pairs}, arXiv e-prints (2024), arXiv:2406.07307.

\bibitem[GLW24]{GLW24} C. Gachet, H. Lin, L. Wang, \emph{Nef cones of fiber products and an application to the cone conjecture}, Forum of Mathematics, Sigma., \textbf{12} (2024), e28, 1 -- 22. 

\bibitem[GKP16]{GKP16} D. Greb, S. Kebekus, and T. Peternell, \emph{\'Etale fundamental groups of Kawamata log terminal spaces, flat
sheaves, and quotients of abelian varieties}, Duke Math. J., \textbf{165} (2016), 1965--2004.

\bibitem[EGAII61]{EGAII} A. Grothendieck and J. Dieudonn\'e, \emph{El\'ements de g\'eom\'etrie alg\'ebrique: II. \'Etude globale \'el\'ementaire de quelques classes de morphismes}, Publ. Math. IH\'ES \textbf{8} (1961), 5--222.



\bibitem[HL23]{HL21} C. D. Hacon and J. Liu, \emph{Existence of flips for generalized lc pairs}, Camb. J. Math., \textbf{11} (2023), no. 4,
795--828. 

\bibitem[HMX18]{HMX18} C. D. Hacon, J. McKernan, C. Xu, \emph{Boundedness of moduli of varieties of general type}, J. Eur. Math. Soc., \textbf{20} (2018), no. 4, pp. 865--901

\bibitem[HL22]{HL22} J. Han and Z. Li, \emph{Weak Zariski decompositions and log terminal models for generalized
pairs}, Mathematische Zeitschrift, \textbf{302} (2022), 707--741.

\bibitem[HLQ23]{HLQ20} J. Han, Y. Liu, and L. Qi, \emph{ACC for local volumes and boundedness of singularities}, J. Algebraic Geom. \textbf{32} (2023), 519--583. 

\bibitem[HH20]{HH19} K. Hashizume, Z. Hu, \emph{On minimal model theory for log abundant lc pairs}, J. Reine Angew. Math.,
\textbf{767} (2020), 109--159. 

\bibitem[Huy16]{Huy16} D. Huybrechts, \emph{Lectures on K3 surfaces}, volume 158 of Cambridge Studies in Advanced
Mathematics. Cambridge University Press, Cambridge, 2016. 

\bibitem[Kaw97]{Kaw97}  Y. Kawamata, \emph{On the cone of divisors of Calabi-Yau fiber spaces}, Int. J. Math., \textbf{8} (1997), 665--687. 

\bibitem[Kaw07]{Kaw07} Y. Kawamata, \emph{Flops connect minimal models}, arXiv e-prints (2007), arXiv:0704.1013. 

\bibitem[KM98]{KM98}
J. Koll\'ar and S. Mori, \emph{Birational geometry of algebraic varieties},
volume 134 of Cambridge Tracts in Mathematics. Cambridge University
Press, Cambridge, 1998. With the collaboration of C. H. Clemens and A.
Corti, Translated from the 1998 Japanese original.

\bibitem[KMM94]{KMM94} S. Keel, K. Matsuki, and J. McKernan, \emph{Log abundance theorem for threefolds}, Duke Math.
J., \textbf{75} (1994), 99--119.

\bibitem[Ko13]{Ko13} J. Koll\'ar, \emph{Singularities of the minimal model program}, volume 200 of Cambridge Tracts in Mathematics. Cambridge University Press, Cambridge, 2013. 

\bibitem[LN59]{LN59} S. Lang and A. N\'eron, \emph{Rational points of abelian varieties over function fields}, Amer. J. Math., \textbf{81} (1959), 95--118.



\bibitem[LOP18]{LOP18} V. Lazi\'c, K. Oguiso, and T. Peternell, \emph{The Morrison-Kawamata cone conjecture and
 abundance on Ricci flat manifolds}, In Uniformization, Riemann-Hilbert correspondence, Calabi-Yau manifolds \& Picard-Fuchs equations, volume 42 of Adv. Lect. Math. (ALM), pages 157--185.
 Int. Press, Somerville, MA, 2018.

\bibitem[LOP20]{LOP20} V. Lazi\'c, K. Oguiso, and T. Peternell, \emph{Nef line bundles on Calabi-Yau threefolds, I}, Int. Math. Res. Not. \textbf{2020} (2020), no. 19, 6070--6119. 

\bibitem[LP13]{LP13} V. Lazi\'c and T. Peternell, \emph{On the cone conjecture for Calabi-Yau manifolds with
Picard number two}, Math. Res. Lett., \textbf{20} (2013), no. 6, 1103--1113.

\bibitem[Li23a]{Li23a} J. Li, \emph{On the cone conjecture for log Calabi-Yau mirrors of Fano 3-folds}, arXiv e-prints (2023), arXiv:2310.02962. 


\bibitem[Li23]{Li23}  Z. Li, \emph{On the relative Morrison-Kawamata cone conjecture (II)}, arXiv e-prints (2023), arXiv:2309.04673. 

\bibitem[LZ22]{LZ22} Z. Li and H. Zhao, \emph{On the relative Morrison-Kawamata cone conjecture}, arXiv e-prints (2022), arXiv:2206.13701. 


\bibitem[Liu02]{Liu02}
Q. Liu, \emph{Algebraic geometry and arithmetic curves}, Oxford Graduate Texts in Mathematics,
vol. 6, Oxford University Press, Oxford, 2002. Translated from the French by R. Ern\'e,
Oxford Science Publications. 

\bibitem[Loo14]{Loo14} E. Looijenga, \emph{Discrete automorphism groups of convex cones of finite type}, Compos. Math.,
\textbf{150} (2014), no. 11, 1939--1962.

\bibitem[Lut24]{Lut24} W. Lutz, \emph{The Morrison cone conjecture under deformation}, arXiv e-prints (2024), arXiv:2410.05949.


\bibitem[Mar11]{Mar11}  E. Markman, \emph{A survey of Torelli and monodromy results for holomorphic-symplectic varieties},
in Complex and differential geometry, volume 8 of Springer Proc. Math., pages 257--322. Springer,
Heidelberg, 2011.

\bibitem[Me19]{Me19} F. Meng, \emph{Some non-vanishing results on log canonical pairs of dimension 4}, arXiv e-prints (2019), arXiv:1906.11451. 

\bibitem[MZ23]{MZ23} F. Meng and Z. Zhuang. \emph{MMP for locally stable families and wall crossing for moduli of stable
pairs}, arXiv e-prints (2023), arXiv:2311.01319.

\bibitem[MQ24]{MQ24} M. Monti, A. Quedo, \emph{The Kawamata-Morrison cone conjecture for generalized hyperelliptic variety}, arXiv e-prints (2024), arXiv:2403.13156. 

\bibitem[Mor96]{Mor96} D. R. Morrison, \emph{Beyond the K\"ahler cone}, Proceedings of the Hirzebruch 65 conference
 on algebraic geometry (1996), 361--376. 




\bibitem[Nak04]{Nak04} N. Nakayama, \emph{Zariski-decomposition and abundance}, MSJ Memoirs, vol. 14,
Mathematical Society of Japan, Tokyo, 2004.

\bibitem[Ogu14]{Ogu14} K. Oguiso,  \emph{Automorphism groups of Calabi-Yau manifolds of Picard number 2}, J. Alg.
Geom., \textbf{23} (2014), no. 4, 775--795.

\bibitem[OS01]{OS01} K. Oguiso and J. Sakurai, \emph{Calabi-Yau threefolds of quotient type}, Asian J. Math., \textbf{5}
(2001), no. 1, 43--77.

\bibitem[Pro21]{Pro21} Y. G. Prokhorov, \emph{Equivariant minimal model program}, 
Russ. Math. Surv., \textbf{76} (2021), no. 3, 461--542. 

\bibitem[PS12]{PS12} A. Prendergast-Smith, \emph{The cone conjecture for some rational elliptic threefolds}, 
Math. Z., \textbf{272} (2012), 589--605.

\bibitem[PS12a]{PS12a} A. Prendergast-Smith, \emph{The cone conjecture for abelian varieties}, J. Math. Sci. Univ. Tokyo,
\textbf{19} (2012), 243--261.

\bibitem[PS15]{PS15}A. Prendergast-Smith, \emph{Finiteness results for 3-folds with semiample anticanonical bundle}, Math. Res. Lett., \textbf{22} (2015), 549--578.

\bibitem[PS26]{PS23} G. Pacienza, A. Sarti, \emph{On the cone conjecture for Enriques manifolds}, \'Epijournal de
G\'eom\'etrie Alg\'ebrique, Special volume in honour of Claire Voisin (2026). 

\bibitem[Sho96]{Sho96}  V.V. Shokurov, \emph{3-fold log models}, J. Math. Sci., \textbf{81}
(1996), no. 3, 2667--2699.


\bibitem[Sho09]{Sho09} V.V. Shokurov, \emph{Letters of a bi-rationalist. VII ordered termination}, Proc. Steklov Inst. Math., \textbf{264} (2009), 178--200. 

\bibitem[SC11]{SC11} V.V. Shokurov, S.R. Choi, \emph{Geography of log models: theory and applications}, Centr. Eur. J. Math., \textbf{9} (2011), 489--534. 

\bibitem[Suz01]{Suz01} 
K. Suzuki, \emph{On Morrison's cone conjecture for klt surfaces with $K_X$ $\equiv$ 0}, Comment. Math. Univ. St. Paul., \textbf{50} (2001), 173--180.

\bibitem[Tak25]{Tak25} T. Takamatsu, \emph{On the finiteness of twists of irreducible symplectic varieties}, Math. Ann. \textbf{392} (2025), 339--371.

\bibitem[Tot10]{Tot09} B. Totaro, \emph{The cone conjecture for Calabi-Yau pairs in dimension two}, Duke Math. J., \textbf{154} (2010),
241--263.


\bibitem[Xie24]{Xie24} L. Xie, \emph{Contraction theorem for generalized pairs}, Algebraic Geometry and Physics, \textbf{1}
(2024), no. 1, 101–124.


\bibitem[Yos95]{Yos95} H. Yoshihara, \emph{Quotients of abelian surfaces}, Publ. Res. Inst. Math. Sci., \textbf{31} (1995), no. 1, 135--143. 
\end{thebibliography}
\end{document}